\documentclass[11pt]{article}
\usepackage{amsmath,amssymb,amsthm}

\newcommand{\prf}{\noindent{\bf Proof.}\ }
\newcommand{\Tbar}{\ensuremath{\overline{\mathcal{T}}}}

\newcommand{\hess}{\operatorname{Hess}}
\newcommand{\grad}{\operatorname{grad}}
\newcommand{\dt}{\frac{d}{dt}}
\newcommand{\dr}{\frac{\partial}{\partial r}}
\newcommand{\dthet}{\frac{\partial}{\partial \vartheta}}

\newcommand{\mathH}{\mathbb{H}}
\newcommand{\mathR}{\mathbb{R}}
\newcommand{\mathS}{\mathbb{S}}

\newcommand{\mathZ}{\mathbb{Z}}
\newcommand{\bA}{\mathbf{A}}
\newcommand{\bB}{\mathbf{B}}
\newcommand{\bC}{\mathbf{C}}

\newcommand{\caF}{\mathcal{F}}

\newcommand{\caH}{\mathcal{H}}

\newcommand{\caQ}{\mathcal{Q}}
\newcommand{\caT}{\mathcal{T}}

\newtheorem{definition}{Definition}[section]
\newtheorem{lemma}[definition]{Lemma}
\newtheorem{theorem}[definition]{Theorem}
\newtheorem{corollary}[definition]{Corollary}

\newtheorem{example}[definition]{Example}
\newtheorem*{theorem'}{Theorem}
\newtheorem*{corollary'}{Corollary}
\newtheorem*{lemma'}{Lemma}

\begin{document}

\title{Behavior of geodesic-length functions on Teichm\"{u}ller space\footnote{2000 Mathematics Subject Classification Primary: 32G15; Secondary: 20H10, 30F60, 53C55.}} 
\author{Scott A. Wolpert}        
\date{September 28, 2007}          
\maketitle

\begin{abstract}
Let $\mathcal T$ be the Teichm\"{u}ller space of marked genus $g$, $n$ punctured Riemann surfaces with its bordification $\Tbar$  the {\em augmented Teichm\"{u}ller space} of marked Riemann surfaces with nodes, \cite{Abdegn, Bersdeg}.  Provided with the WP metric $\Tbar$ is a complete  $CAT(0)$ metric space, \cite{DW2, Wlcomp, Yam2}.  An invariant of a marked hyperbolic structure is the length $\ell_{\alpha}$ of the geodesic $\alpha$ in a free homotopy class.  A basic feature of Teichm\"{u}ller theory is the interplay of two-dimensional hyperbolic geometry,  Weil-Petersson (WP) geometry and the behavior of geodesic-length functions.  Our goal is to develop the understanding of the intrinsic local WP geometry through a study of the gradient and Hessian of geodesic-length functions.  Considerations include expansions for the WP pairing of gradients, expansions for the Hessian and covariant derivative, comparability models for the WP metric, as well as the behavior of WP geodesics including a description of the Alexandrov tangent cone at the augmentation. Approximations and applications for geodesics close to the augmentation are developed. An application for fixed points of group actions is described.  Bounding configurations and functions on the hyperbolic plane is basic to our approach.  Considerations include analyzing the orbit of a discrete group of isometries and bounding sums of the inverse square exponential-distance.  
\end{abstract}

\pagebreak

\tableofcontents

\section{Introduction}

\subsection{Background}
Collections of geodesic-length functions provide local coordinates for Teichm\"{u}ller space.  Gardiner provided a formula for the differential of geodesic-length and for the WP gradient of geodesic-length, \cite{Gardtheta}.  In \cite{WlFN} a geometric and analytic description of the deformation vector field for the infinitesimal Fenchel-Nielsen right twist $t_{\alpha}$ revealed the twist-length duality formula 
$2t_{\alpha}=i\grad\ell_{\alpha}$, \cite[Thrm. 2.10]{WlFN}.  The formula combines with the calculation of the twist variation of length to give a formula for the WP Riemannian pairing
\[
\langle\grad\ell_{\alpha},i\grad\ell_{\beta}\rangle=-2\sum_{p\in\alpha\cap\beta}\cos\theta_p
\]
for geodesics $\alpha,\,\beta$, \cite[Thrm. 3.3]{Wlsymp}.  Hyperbolic trigonometric formulas were also provided for the Lie derivative $t_{\alpha}t_{\beta}\ell_{\gamma}$ and the Lie bracket $[t_{\alpha},t_{\beta}]$, \cite[Thrm. 3.4, Thrm. 4.8]{Wlsymp}.  The convexity of geodesic-length functions along twist and earthquake paths is a consequence, \cite{Kerck, Wlsymp}.  In \cite{Wlnielsen} an analytic description of deformation vector fields was applied to find that geodesic-length functions are strictly convex along WP geodesics.  The WP geodesic convexity of Teichm\"{u}ller space is a consequence.

\subsection{The gradient and Hessian of geodesic-length}
Riera has provided a formula for the pairing of geodesic-length gradients, \cite[Theorem 2]{Rier}.  Let $R$ be a finite area Riemann surface $R$ with uniformization group $\Gamma\subset PSL(2;\mathbb R)$ acting on the upper half plane $\mathbb H$.   For closed geodesics $\alpha,\beta$ on $R$ with deck transformations $A, B$ with axes $\tilde\alpha,\tilde\beta$ on $\mathbb H$ the pairing formula is
\[
\langle\grad\ell_{\alpha},\grad\ell_{\beta}\rangle=\frac{2}{\pi}\bigl(\ell_{\alpha}\delta_{\alpha\beta}\ +\ 
\sum_{\langle A\rangle\backslash\Gamma/\langle B\rangle}(u\log\frac{u+1}{u-1}-2\bigr)\bigr)
\]
for the Kronecker delta $\delta_*$, where for $C\in\langle A\rangle\backslash\Gamma/\langle B\rangle$ then $u=u(\tilde\alpha,C(\tilde\beta))$ is the cosine of the intersection angle if $\tilde\alpha$ and $C(\tilde\beta)$ intersect and is otherwise $\cosh d(\tilde\alpha,C(\tilde\beta))$; for $\alpha=\beta$ the double-coset of the identity element is omitted from the sum.  In Section \ref{gradhess} we compare the sum to $\sum_{\langle A\rangle\backslash\Gamma/\langle B\rangle}e^{-2d}$, use a mean value estimate and analysis of the orbit of $\Gamma$ to find the following.
\begin{lemma'}
The WP pairing of geodesic-length gradients of disjoint geodesics $\alpha,\beta$ satisfies 
\[
0<\langle\grad\ell_{\alpha},\grad\ell_{\beta}\rangle-\frac{2}{\pi}\ell_{\alpha}\delta_{\alpha\beta}\quad\mbox{is}
\quad O(\ell_{\alpha}^2\ell_{\beta}^2)
\]
where for $c_0$ positive the remainder term constant is uniform for $\ell_{\alpha},\ell_{\beta}\le c_0$.  The shortest nontrivial segment on $\mathbb H/\Gamma$  connecting $\alpha$ and $\beta$ contributes a term of order $\ell_{\alpha}^2\ell_{\beta}^2$ to the pairing.
\end{lemma'}
Further considerations provide an almost sharp general bound for the norm of the gradient.  The expansion is one of several suggesting introduction of the root geodesic-lengths $\ell_{\alpha}^{1/2}$ and the root geodesic-length gradients $\lambda_{\alpha}=\grad\ell_{\alpha}^{1/2}$ for which $\langle\lambda_{\alpha},\lambda_{\alpha}\rangle=\frac{1}{2\pi}+O(\ell_{\alpha}^3)$.

A description of the deformation vector field on $\mathbb H$ was developed from the Eichler-Shimura isomorphism in \cite{Wlnielsen}.  The description enables the calculation of a second order deformation beginning from the harmonic Beltrami differential representing the infinitesimal deformation.  The resulting expressions \cite[Lemma 4.2, Coro. 4.3, Lemma 4.4, Coro. 4.5]{Wlnielsen} are complicated.  We show now that the infinite horizontal strip is better suited for the calculation.  The resulting expressions are simpler and the formulas enable Hessian bounds.  In Section \ref{defvf} we present a formula for the Hessian and find the following.
\begin{corollary'}
The complex and real Hessians of geodesic-length are uniformly comparable 
\[
\partial\overline\partial\ell_{\alpha}\le\hess\ell_{\alpha}\le3\partial\overline\partial\ell_{\alpha}.
\]  
The first and second derivatives of geodesic-length satisfy 
\[
2\ell_{\alpha}\Ddot\ell_{\alpha}[\mu,\mu]-\dot\ell_{\alpha}^2[\mu]-3\dot\ell_{\alpha}^2[i\mu]\ge0\quad\mbox{and}\quad\ell_{\alpha}\partial\overline\partial\ell_{\alpha}-2\partial\ell_{\alpha}\overline\partial\ell_{\alpha}\ge 0
\] 
for a harmonic Beltrami differential $\mu=\overline{\varphi}(ds^2)^{-1}$ with equality only for the elementary
 $\varphi=a(\frac{dz}{z})^2$.  
\end{corollary'}

We apply the formulas to study geodesic-length and the length of compactly supported geodesic laminations.  A geodesic lamination is a closed subset of a surface which is foliated by complete geodesics.  A measured geodesic lamination with compact support is a compact geodesic lamination with a transverse measure with support the lamination.  A transverse measure $\upsilon$ is alternately specified on the space $\mathcal G$ of complete geodesics on the upper half plane.  For a point $p$ and a complete geodesic $\tilde\beta$ on $\mathbb H$ the distance $d(p,\tilde\beta)$ is defined.  We consider the inverse square exponential-distance $\mathbb P_{\upsilon}(p)=\int_{\mathcal G}e^{-2d(p,\tilde\beta)}\,d\upsilon(\tilde\beta)$.  The integral compares to $\ell_{\upsilon}\int_{\mathbb H}e^{-2d(p,q)}\,dA$.      The map $\upsilon\rightarrow \mathbb P_{\upsilon}$ from $\mathcal{ML}_0$ measured geodesic laminations with compact support 
 to $C(R)$ is continuous in the compact-open topology; see Lemma \ref{Pbound1}.  In Section \ref{bdhess} we bound the Hessian of geodesic-length as follows.

\begin{theorem'}
The Hessian of length of a measured geodesic lamination $\upsilon$ is bounded in terms of the WP pairing and the weight $\mathbb P_{\upsilon}$ as follows
\[
\langle\mu,\mu\mathbb{P}_{\upsilon}\rangle \le 3\pi\partial\overline{\partial}\ell_{\upsilon}[\mu,\mu]\le 16\langle\mu,\mu\mathbb{P}_{\upsilon}\rangle 
\]
and
\[
\langle\mu,\mu\mathbb{P}_{\upsilon}\rangle \le 3\pi\Ddot\ell_{\upsilon}[\mu,\mu]\le 48\langle\mu,\mu\mathbb{P}_{\upsilon}\rangle 
\]
for $\mu$ a harmonic Beltrami differential.  There are positive functions $c_1$ and $c_2$ such that
\[
c_1(inj_{red})\ell_{\upsilon}\langle\ ,\ \rangle 
\le\partial\overline{\partial}\ell_{\upsilon},\ \Ddot\ell_{\upsilon}\le c_2(inj_{red})\ell_{\upsilon}
\langle\ ,\ \rangle 
\]
for the reduced injectivity radius with $c_1(\rho)$ an increasing function vanishing at the origin and $c_2(\rho)$ a decreasing function tending to infinity at the origin.  
\end{theorem'}

In Section \ref{gradhess} we separate the contribution of the zeroth rotational term in the collar and apply the Cauchy Integral Formula and Schwarz Lemma to find the following elaboration.
\begin{theorem'}
The variation of geodesic-length $\ell_{\alpha}$ satisfies
\[
2\ell_{\alpha}\Ddot\ell_{\alpha}[\mu,\mu]-\dot\ell_{\alpha}^2[\mu]-3\dot\ell_{\alpha}^2[i\mu]\quad\mbox{is}\quad O(\ell_{\alpha}^3\|\mu\|^2_{WP})
\]
and
\[
\ell_{\alpha}\partial\overline\partial\ell_{\alpha}[\mu,\mu]-2\partial\ell_{\alpha}[\mu]\overline\partial\ell_{\alpha}[\mu]\quad\mbox{is}\quad O(\ell_{\alpha}^3\|\mu\|^2_{WP})
\]
for a harmonic Beltrami differential $\mu$ where for $c_0$ positive the remainder term constant is uniform for $\ell_{\alpha}\le c_0$.
\end{theorem'}
An intrinsic expansion for covariant differentiation follows from the direct relation with the Hessian; see Section \ref{gradhess}.  The almost complex structure of Teichm\"{u}ller space is parallel with respect to WP covariant differentiation.  

\begin{theorem'}
The WP covariant derivative $D$ of the root geodesic-length gradient $\lambda_{\alpha}$ satisfies
\[
D_U\lambda_{\alpha}\ =\ 3\ell_{\alpha}^{-1/2}\langle J\lambda_{\alpha},U\rangle J\lambda_{\alpha}\ +\ O(\ell_{\alpha}^{3/2}\|U\|)
\]
for $J$ the almost complex structure and $\|\ \|$ the WP norm.  For $c_0$ positive the remainder term constant is uniform for $\ell_{\alpha}\le c_0$.
\end{theorem'}

\subsection{The augmentation, WP comparisons and the Alexandrov tangent cone}
The augmentation $\Tbar$ is a bordification (a partial compactification) introduced by extending the range of Fenchel-Nielsen parameters, \cite{Abdegn, Bersdeg}.  The added points correspond to unions of marked hyperbolic surfaces with formal pairings of cusps.  For a geodesic-length $\ell_{\alpha}$ equal to zero in place of the geodesic $\alpha$ there appears a pair of cusps; the marking map is a homeomorphism from the complement of a curve to a union of hyperbolic surfaces.  In general for a collection $\sigma$ of homotopically non trivial, non peripheral simple closed curves the locus $\mathcal T(\sigma)=\{\ell_{\alpha}=0\mid\alpha\in\sigma\}\subset\Tbar$ has codimension twice the cardinality of $\sigma$.  We show in Corollary \ref{termbhv} that  on $\Tbar$ the distance $d_{\mathcal T(\sigma)}$ to the closure of $\mathcal T(\sigma)$ satisfies
$d_{\mathcal T(\sigma)}\le(2\pi\sum_{\alpha\in\sigma}\ell_{\alpha})^{1/2}\quad\mbox{and}\quad d_{\mathcal T(\sigma)}=(2\pi\sum_{\alpha\in\sigma}\ell_{\alpha})^{1/2}\, +\, O(\sum_{\alpha\in\sigma}\ell_{\alpha}^{5/2})$.  

As a means for understanding the geometry and the drop of dimension at the augmentation we describe a frame for the tangent bundle of $\mathcal T$ near $\mathcal T(\sigma)$ in terms of geodesic-lengths.  We first consider the description of the bordification in terms of representations of groups and the Chabauty topology.  We find that the pairing $\langle\grad\ell_{\beta},\grad\ell_{\gamma}\rangle$ is continuous at $\mathcal T(\sigma)$ for closed curves $\beta$ and $\gamma$ disjoint from $\sigma$.  Then for $p\in\mathcal T(\sigma)$ corresponding to a union of marked hyperbolic surfaces we define a {\em relative length basis} as a collection $\tau$ of simple closed curves disjoint from the elements of $\sigma$ such that at $p$ the geodesic-lengths $\ell_{\beta},\,\beta\in\tau,$ provide local coordinates for $\mathcal T(\sigma)$.  We consider the WP Riemannian pairing matrix $P$ for $\{\lambda_{\alpha},i\lambda_{\alpha},\grad\ell_{\beta}\}_{\alpha\in\sigma,\,\beta\in\tau}$ where pairings with $\lambda_{\alpha}$ on $\{\ell_{\alpha}=0\}$ are determined by continuity.  We show in Lemma \ref{paircont} that $P$ is the germ of a continuous map from $\Tbar$ to the general linear group $GL(\mathbb R)$ and in Corollary \ref{WPcont} that the pairing satisfies $P(q)=P(p)+O(d_{WP}(p,q))$. The local frames provide for a bundle extension over $\Tbar$ of the tangent bundle of $\mathcal T$.  

We combine techniques in Section \ref{WPcompar} to find comparisons and expansions for the WP metric.  Comparisons are given in terms of geodesic-lengths for a partition $\sigma$, a maximal collection of $3g-3+n$ disjoint distinct simple non trivial non peripheral free homotopy classes.  Comparisons are for the Bers regions $\mathcal B=\{\ell_{\alpha}\le c \mid \alpha\in\sigma \}$.  The complex differentials $\{\partial\ell_{\alpha}^{1/2}\}_{\alpha\in\sigma}$ provide a global frame for the tangent bundle of $\mathcal T$, \cite[Thrm. 3.7]{WlFN}.  
\begin{theorem'}
The WP metric is comparable to a sum of first and second derivatives of  geodesic-length functions for a partition $\sigma$ and $J$ the almost complex structure as follows
\[
\langle\ ,\ \rangle \ \asymp\ \sum_{\alpha\in\sigma} (d\ell_{\alpha}^{1/2})^2+(d\ell_{\alpha}^{1/2}\circ J)^2\ \asymp\ \sum_{\alpha\in\sigma}\hess\ell_{\alpha} 
\]
where given $c_0$ positive there are positive constants $c_1$ and $c_2$ for the comparability on the Bers region $\mathcal B=\{\ell_{\alpha}\le c_0\mid\alpha\in\sigma\}$.
The WP metric has the expansions
\begin{align*}
\langle\ ,\ \rangle \ =&\ 2\pi\sum_{\alpha\in\sigma}(d\ell_{\alpha}^{1/2})^2+(d\ell_{\alpha}^{1/2}\circ J)^2\ +\ O((\sum_{\alpha\in\sigma}\ell_{\alpha}^3)\,\langle\ ,\ \rangle )\\ =&\ \frac{\pi}{6}\sum_{\alpha\in\sigma}\frac{\hess\ell_{\alpha}^2}{\ell_{\alpha}} \ +\ O((\sum_{\alpha\in\sigma}\ell_{\alpha}^2)\,\langle\ ,\ \rangle )
\end{align*}
at the maximal frontier point of $\overline{\mathcal T}$ for the partition $\sigma$.
\end{theorem'}

The expansions are a counterpart to Masur's original local expansion, \cite{Msext}.  To date all expansions involve considerations of algebraic geometry and the local plumbing family $\{zw=t\}\rightarrow\{|t|<1\}$.  The present formulas are intrinsic to hyperbolic geometry and are valid on the {\em larger} Bers regions.  The total geodesic-length of a partition is uniformly convex on WP geodesics in the Bers region.  The above comparability on Bers regions provides for local considerations for Brock's quasi isometry of the pants complex $C_{\mathbf P}$ to Teichm\"{u}ller space with the WP metic, \cite{Brkwp}.   A Lipschitz model for Teichm\"{u}ller space can now be constructed from the pants complex $C_{\mathbf P}$ and corresponding Bers regions.

In \cite{Wlext} we use the current techniques of investigation to develop a normal from for the WP Riemanninan connection for surfaces with short geodesics.  The normal form is used to establish approximation of geodesics in strata of the bordification $\Tbar$.  Considerations are then combined to establish the convexity along WP geodesics of the functions the distance between horocycles for a hyperbolic metric with cusps.   

In Sections \ref{WPconn} and \ref{Alextgtcn} we combine techniques to investigate geodesics terminating at the augmentation. We consider in Corollary \ref{termbhv} the tangent field for a geodesic $\gamma(t)$ terminating at $\mathcal T(\sigma)$.  The geodesic-length functions for $\alpha\in\sigma$ have one-sided $C^2$ expansions $\ell_{\alpha}^{1/2}(\gamma(t))=a_{\alpha}t+O(t^5)$, while the tangent field pairings have one-sided $C^1$ expansions with $\langle\frac{d}{dt},J\lambda_{\alpha}\rangle$ bounded as $O(t^4)$.  Additionally on a geodesic a geodesic-length function either vanishes identically or has positive initial one-sided derivative.  
Earlier expansions indicate that the metric tangent cone at the augmentation is singular, \cite{DW2, Wlcomp, Yam2}.  We consider the Alexandrov tangent cone as a generalization of the tangent space.  See \cite[Chap. II.3]{BH} for the notion of Alexandrov angle for $CAT(0)$ metric spaces.  A triple of points $(p,q,r)$ has Euclidean comparison triangle with angle $\angle(p,q,r)$ at $p$ determined by the Law of Cosines.  For constant speed geodesics $\gamma_0(t),\gamma_1(t)$ with common initial point the comparison angle for $(\gamma_0(0),\gamma_1(t),\gamma_1(t'))$ is a non decreasing function of $t$ and $t'$.  The Alexandrov angle $\angle(\gamma_0,\gamma_1)$ is the limit for $t$ and $t'$ tending to zero.  The Alexandrov tangent cone $AC_p$ is the set of constant speed geodesics beginning at $p$ modulo the  equivalence relation same speed and at zero angle.  

We present for an augmentation point $p\in\mathcal T(\sigma)$ an isometry between the Alexandrov tangent cone $AC_p$ and the product $\mathbb R_{\ge 0}^{|\sigma|}\times T_p\mathcal T(\sigma)$ with the first factor the Euclidean orthant and the second factor the stratum tangent space with WP metric. The Alexandrov tangent cone is given the structure of an inner product by the formal relation 
$\langle\gamma_0,\gamma_1\rangle =\|\gamma_0'\|\|\gamma_1'\|\cos\angle (\gamma_0,\gamma_1)$.  For a relative length basis $\tau$ the mapping for a WP geodesic $\gamma(t)$ beginning at $p$ is given by associating
\[
\Lambda:\gamma\rightarrow(2\pi)^{1/2}\frac{d\mathcal L(\gamma)}{dt}(0)
\]
the initial one-sided derivative of $\mathcal L(\gamma(t))=(\ell_{\alpha}^{1/2},\ell_{\beta}^{1/2})_{\alpha\in\sigma,\,\beta\in\tau}(\gamma(t))$.  The tuple $(\ell_{\beta}^{1/2})_{\beta\in\tau}$ provides local coordinates for $\mathcal T(\sigma)$ and thus  $(2\pi)^{1/2}\bigl(\frac{d\ell_{\beta}^{1/2}(\gamma)}{dt}(0)\bigr)_{\beta\in\tau}$ defines a tangent vector in $T_p\mathcal T(\sigma)$.  The following is presented in Section \ref{Alextgtcn}.  
\begin{theorem'}
The mapping $\Lambda$ from the WP Alexandrov tangent cone $AC_p$ to $\mathbb R_{\ge 0}^{|\sigma|}\times T_p\mathcal T(\sigma)$ is an isometry of cones with restrictions of inner products.   A WP terminating geodesic $\gamma$ with a root geodesic-length function initial derivative $\frac{d\ell_{\alpha}^{1/2}(\gamma)}{dt}(0)$ vanishing is contained in the stratum $\{\ell_{\alpha}=0\},\ \mathcal T(\sigma)\subset\{\ell_{\alpha}=0\}$.  Geodesics $\gamma_0$ and $\gamma_1$ at zero angle have comparison angles $\angle(p,\gamma_0(t),\gamma_1(t))$ bounded as $O(t)$.
\end{theorem'}
We present applications of the geometry of the tangent cone.

\subsection{The model metric $4dr^2+r^6d\vartheta^2$ and Fenchel-Nielsen coordinates}

Understanding of WP geometry is forwarded by the comparison to the model metric $4dr^2+r^6d\vartheta^2$ for the plane $\mathbb R^2$.  Preliminary forms of a model metric were presented earlier, \cite{DW2,Wlthur,Yam2}.  We describe the geometry of the model metric and the comparison between a product of model metrics and WP.  The comparison accounts for the expansion for the WP connection, recent properties of WP sectional curvature, as well as the description of the Alexandrov tangent cone.  We find that the model metric also provides expansions for the WP metric in Fenchel-Nielsen coordinates.

Basic properties include the following.  The $2$-dimensional K\"{a}hler model metric is non complete on $\mathbb R^2-\{0\}$.  Coordinate vector fields are evaluated with $\langle\dr,\dr\rangle=4,\ \langle\dr,\dthet\rangle=0$ and 
$\langle\dthet,\dthet\rangle=r^6$.  The Riemannian connection $D$ is determined through the relations for vector fields $D_U\langle V,W\rangle=\langle D_UV,W\rangle+\langle V, D_UW\rangle$ and $D_UV-D_VU=[U,V]$.   The formulas for the model connection are
\[
D_{\dr}\tfrac{\partial}{\partial r}=0,\quad D_{\dthet}\tfrac{\partial}{\partial r}=D_{\dr}\tfrac{\partial}{\partial\vartheta}=\tfrac3r\tfrac{\partial}{\partial\vartheta}\quad\mbox{and}\quad D_{\dthet}\tfrac{\partial}{\partial\vartheta}=\tfrac{-3}{4}r^5\tfrac{\partial}{\partial r}.
\]
The formulas are combined with the derivation property to evaluate the curvature tensor $R(U,V)W=D_UD_VW-D_VD_UW-D_{[U,V]}W$.  The Lie bracket of coordinate vector fields $[\dr,\dthet]$ vanishes and from the above evaluations 
$R(\dr,\dthet)\dthet=\frac{-3}{2}r^4\dr$.  The $2$-plane $\dr\wedge\dthet$ has area $2r^3$ and the sectional curvature of the model metric is
\[
\frac{\langle R(\dr,\dthet)\dthet,\dr\rangle}{\|\dr\wedge\dthet\|^2}=\frac{-3}{2r^2}. 
\]

A comparison of metrics involves the deformation parameter for the angle $\vartheta$, in particular involves counterparts for the vector field $\dthet$ and the differential $d\vartheta$.  We explain in Section \ref{WPconn} that for an increment of the angle of $2\pi$ to correspond to a full rotation (a Dehn right twist in the mapping class group) the corresponding unit Fenchel-Nielsen infinitesimal angle variation is $T_{\alpha}=(2\pi)^{-1}\ell_{\alpha}^{3/2}J\lambda_{\alpha}=(2\pi)^{-1}\ell_{\alpha}t_{\alpha}$.  The variation $T_{\alpha}$ is the counterpart for 
$\dthet$.  The definition of $T_{\alpha}$ is intrinsic not involving auxiliary choices.  A counterpart candidate for the differential $d\vartheta$ is the differential of a Fenchel-Nielsen angle; see Section \ref{WPcompar}.  There are limitations to the selection of the Fenchel-Nielsen angle.  The first is that the definition requires the choice of free homotopy classes of two auxiliary simple closed curves; in fact arbitrary pairs of disjoint classes disjoint from the original curve occur.  The second is that the differential geometric properties of the Fenchel-Nielsen angle on Teichm\"{u}ller space are not available.  In Definition \ref{gauge} we use the Fenchel-Nielsen angle variation and WP pairing to introduce the Fenchel-Nielsen gauge $1$-form $\varrho_{\alpha}=2\pi(\ell_{\alpha}^{3/2}\langle\lambda_{\alpha},\lambda_{\alpha}\rangle)^{-1}\langle\ ,J\lambda_{\alpha}\rangle$ where $d\ell_{\alpha}^{1/2}\circ J=\langle\ ,J\lambda_{\alpha}\rangle$.  The Fenchel-Nielsen gauge satisfies the basic property $\varrho_{\alpha}(T_{\alpha})=1$ and the definition does not involve auxiliary choices.  The differential geometric properties of Fenchel-Nielsen gauges are presented in Lemma \ref{drho}, including an expansion relating gauges to angles.

We are ready to relate the WP and model geometry.  Following the lead of existing expansions we set for a geodesic-length $2\pi^2r_{\alpha}^2=\ell_{\alpha}$ for which $2^{1/2} \pi dr_{\alpha}=\langle\ ,\lambda_{\alpha}\rangle$ in terms of the root geodesic-length gradient.  The above comparison Theorem becomes the following.

\begin{theorem'}
The WP metric is comparable for a partition $\sigma$ as follows
\[
\langle\ ,\ \rangle \ \asymp\ \sum_{\alpha\in\sigma} 4dr_{\alpha}^2+r_{\alpha}^6\varrho_{\alpha}^2 
\]
for the Fenchel-Nielsen gauges where given $c_0$ positive there are positive constants $c_1$ and $c_2$ for the comparability on the Bers region $\mathcal B=\{\ell_{\alpha}\le c_0\mid\alpha\in\sigma\}$.
The WP metric has the expansion
\[
\langle\ ,\ \rangle \ = \pi^3 \sum_{\alpha\in\sigma} 4dr_{\alpha}^2+r_{\alpha}^6\varrho_{\alpha}^2 \ +\ O((\sum_{\alpha\in\sigma}\ell_{\alpha}^3)\,\langle\ ,\ \rangle )\\ 
\]
at the maximal frontier point of $\overline{\mathcal T}$ for the partition $\sigma$.
\end{theorem'}
Again the comparison is for a larger Bers region.  Constants agree with earlier expansions, \cite[Coro. 4]{Wlcomp}.  

We are ready to relate the product of model metrics $\sum_{\alpha\in\sigma} 4dr_{\alpha}^2+r_{\alpha}^6d\vartheta_{\alpha}^2$ to the WP pairing $\pi^{-3}\langle\ ,\ \rangle$.  The comparison is based on $\frac{\partial}{\partial r_{\alpha}}$ having as analog $2^{3/2}\pi^2\lambda_{\alpha}$ and $\frac{\partial}{\partial \vartheta_{\alpha}}$ having as analog the Fenchel-Nielsen angle variation $T_{\alpha}=(2\pi)^{-1}\ell_{\alpha}^{3/2}J\lambda_{\alpha}$.  The K\"{a}hler form for the model metric is $\sum_{\alpha\in\sigma}2r^3_{\alpha}dr_{\alpha}d\vartheta_{\alpha}=(4\pi^4)^{-1}\sum_{\alpha\in\sigma}\ell_{\alpha}d\ell_{\alpha}d\vartheta_{\alpha}$ which for the Fenchel-Nielsen angle $2\pi\tau_{\alpha}/\ell_{\alpha}$ corresponds to the known formula $(2\pi^3)^{-1}\sum_{\alpha\in\sigma}d\ell_{\alpha}d\tau_{\alpha}$, 
\cite{Wldtau, Wlcusps}.  The above Theorem provides equality of norms modulo the higher order $O$-term. The above expansion for the pairing 
$\langle\lambda_{\alpha},\lambda_{\alpha}\rangle=\frac{1}{2\pi}+O(\ell_{\alpha}^3)$,  Theorem on the WP connection and definitions for $2\pi^2r_{\alpha}^2=\ell_{\alpha}$ and $T_{\alpha}$ combine to provide the following covariant derivative formulas
\begin{gather*}
D_{\lambda_{\alpha}}\lambda_{\alpha}=O(\ell_{\alpha}^{3/2}),\quad D_{T_{\alpha}}(2^{3/2}\pi^2\lambda_{\alpha})=D_{2^{3/2}\pi^2\lambda_{\alpha}}T_{\alpha}=\frac{3}{r_{\alpha}}T_{\alpha}+O(\ell_{\alpha}^3)\\ \mbox{and}\quad
D_{T_{\alpha}}T_{\alpha}=\frac{-3}{4}r_{\alpha}^5(2^{3/2}\pi^2\lambda_{\alpha})+O(\ell_{\alpha}^{9/2})
\end{gather*}
in direct correspondence to the formulas for the model connection.  Relatedly the root geodesic-length along a WP geodesic with $\ell_{\alpha}^{1/2}$ small is modeled by the radius function along a geodesic of the model metric.  The principal terms of the derivative of $2^{3/2}\pi^2\lambda_{\alpha}$ and $T_{\alpha}$ correspond to the derivative of $\frac{\partial}{\partial r}$ and $\frac{\partial}{\partial \theta}$. 

We next describe comparisons for curvature.  WP sectional curvature is negative with an upper bound of $O(-\ell)$ and a lower bound of $O(-\ell^{-1})$ for $\ell$ the length of the shortest closed geodesic, \cite{Zh2}.  Huang describes tangent planes realizing the bounds.  In particular in agreement with the model metric the curvature is $O(-\ell_{\alpha}^{-1})$ for a plane $\frac{\partial}{\partial r_{\alpha}}\wedge\frac{\partial}{\partial \vartheta_{\alpha}}$.  Huang also shows from the {\em almost product} structure that planes spanned by an $\alpha$-factor tangent and an $\alpha'$-factor tangent have curvature $O(-\ell_{\alpha}-\ell_{\alpha'})$.  Additional planes of small curvature correspond to components of surfaces represented in the augmentation.  He further finds that the subset $\mathcal T_{\ell\ge c}\subset\mathcal T$ of surfaces with a fixed positive lower bound for the length of the shortest closed geodesic has a lower bound for sectional curvature independent of genus and number of punctures, \cite{Zh3}.  The uniform curvature bounds are consistent with a uniform bound for the magnitude of the $C^2$-norm of the difference of the WP and product model metrics. 

I would like to thank Sumio Yamada for conversations.

\section{Preliminaries}

\subsection{Collars and cusp regions}
A Riemann surface with hyperbolic metric can be considered as the union of a $thick$ region where the injectivity radius is bounded below by a positive constant and the complementary $thin$ region.  The totality of all $thick$ regions of Riemann surfaces of a given topological type forms a compact set of metric spaces in the Gromov-Hausdorff topology.  A $thin$ region is a disjoint union of collar and cusp regions.  We describe basic properties of collars and cusp regions including bounds for the injectivity radius and separation. 

We follow Buser's presentation \cite[Chap. 4]{Busbook}.   For a geodesic $\alpha$ of length $\ell_{\alpha}$ on a Riemann surface the collar about the geodesic is $\mathcal C(\alpha)=\{d(p,\alpha)\le w(\alpha)\}$ for the width $w(\alpha)$, $\sinh w(\alpha)\sinh \ell_{\alpha}/2=1$.  For $\mathbb H$ the upper half plane with hyperbolic distance $d(\ ,\ )$ a collar is covered by the region  
$\{d(z,i\mathbb R^+)\le w(\alpha)\}\subset \mathbb H$ with deck transformations generated by $z\rightarrow e^{\ell_{\alpha}}z$.  The quotient $\{d(z,i\mathbb R^+)\le w(\alpha)\}\slash \langle z\rightarrow e^{\ell_{\alpha}}z\bigr >$ embeds into the Riemann surface. For $z$ in $\mathbb H$ the prescribed region is approximately $\{\ell_{\alpha}/2\le \arg z\le \pi-\ell_{\alpha}/2\}$.   A cusp region $\mathcal C_{\infty}$ is covered by the region $\{\Im z\ge 1/2\}\subset\mathbb H$ with deck transformations generated by $z\rightarrow z+1$.  The quotient $\{\Im z\ge 1/2\}\slash \langle z\rightarrow z+1\rangle $ embeds into the Riemann surface.  The boundary of a collar $\mathcal C(\alpha)$ for $\ell_{\alpha}$ bounded and boundary of a cusp region $\mathcal C_{\infty}$ have length approximately $2$.

\begin{theorem}
\label{collars}
For a Riemann surface of genus $g$ with $n$ punctures given pairwise disjoint simple closed geodesics $\alpha_1,\dots,\alpha_m$ there exist simple closed geodesics $\alpha_{m+1},\dots,\alpha_{3g-3+n}$ such that $\alpha_1,\dots,\alpha_{3g-3+n}$ are pairwise disjoint.  The collars $\mathcal C(\alpha_j)$ about $\alpha_j$, $1\le j\le 3g-3+n$, and the cusp regions are mutually pairwise disjoint.
\end{theorem}

On $thin$ the injectivity radius is bounded below in terms of the distance into a collar or cusp region.  For a point $p$ of a collar or cusp region write $inj(p)$ for the injectivity radius of the Riemann surface and $\delta(p)$ for the distance to the boundary of the collar  or cusp region.   The injectivity radius is bounded as follows, \cite[II, Lemma 2.1]{Wlspeclim}.

\begin{lemma}
\label{enhcollar}
The product $inj(p)\,e^{\delta(p)}$ of injectivity radius and exponential distance to the boundary is bounded below by a positive constant.  
\end{lemma}

Simple infinite geodesics occur as leaves of laminations.  We are interested in the behavior of $inj(p)$ along a simple complete (infinite or closed) geodesic.  The standard consideration for cusp regions generalizes as follows.

\begin{lemma}
\label{separ}
A simple geodesic is either disjoint from the thin region or is the core of an included collar or crosses at least half an included collar or crosses a cusp region.
\end{lemma}
\prf  On $\mathbb H$ let $A$ be the deck transformation $z\rightarrow e^{\ell_{\alpha}}z$ for a collar and $z\rightarrow z+1$ for a cusp.   Consider a lift $\tilde\beta$ of the simple geodesic with endpoints $b_1<b_2$ in $\mathbb R\cup\{\infty\}$.  If for a collar either $b_1,b_2$ is infinity or zero then $\beta$ spirals to $\alpha$ and so crosses half the collar; if similarly the product $b_1b_2$ is negative then $\beta$ intersects $\alpha$ and crosses the collar.  If for a cusp either $b_1,b_2$ is infinity then $\beta$ crosses the cusp region.  If $\beta$ does not cross a region then we consider the interval $[b_1,b_2]$ and its translates by $A$ and $A^{-1}$.    If the interval and a translate overlap then $\tilde\beta$ intersects either $A(\tilde\beta)$ or $A^{-1}(\tilde\beta)$ interior to $\mathbb H$.  Since $\beta$ is simple it follows that the interval $[b_1,b_2]$ is contained in a fundamental domain for $A$.  The extremal position for an interval is given by the configuration with $A(b_1)=b_2$ or with $A(b_2)=b_1$.  
For a cusp region the height of $\tilde\beta$ is $1/2$ and $\beta$ only intersects the boundary of the region.  
For a collar (with $b_1>0$) we consider the quadrilateral with vertices $b_1,ib_1,ib_2,b_2$ and vertex angles $0,\frac{\pi}{2},\frac{\pi}{2},0$.  The hyperbolic trigonometry of the quadrilateral provides the equality $\sinh d(\alpha,\beta)\sinh \ell_{\alpha}/2=1$, \cite[Thrm. 2.3.1]{Busbook}. The geodesic $\beta$ only intersects the boundary of the collar.
The proof is complete.

\subsection{Mean value estimates}

A mean value estimate provides a simple method for bounding a function.  Our considerations involve bounds for holomorphic quadratic differentials and sums of translates of the exponential-distance on the upper half plane.  Holomorphic functions as well as eigenfunctions with positive eigenvalue of the hyperbolic Laplacian satisfy the area mean value property, \cite[Coro. 1.3]{Fay}.  A holomorphic function or $\lambda>0$ eigenfunction satisfies
\[
f(\rho)=c(\rho,\lambda)\int_{\mathbf B(p,\rho)}f\,dA
\]
on the upper half plane for $\mathbf B$ a hyperbolic metric ball and $dA$ the hyperbolic area element.  We consider the following specializations. 

A harmonic Beltrami differential is given as $\mu=\overline\varphi(ds^2)^{-1}$ for $\varphi$ a holomorphic quadratic differential and $ds^2$ the hyperbolic metric.  A harmonic Beltrami differential satisfies a mean value estimate 
\[
|\mu|(p)\le c(\rho)\int_{\mathbf B(p,\rho)}|\mu|\,dA
\]
as a consequence of the rotational symmetry of the metric and the mean value property for holomorphic functions for circles.   Eigenfunctions with positive eigenvalue satisfy the mean value property.  For positive eigenvalue $\lambda=s(s-1)$ the rotationally invariant potential about a point of $\mathbb H$ is a negative function of distance asymptotic to $-ce^{-sd(p,q)}$ for large distance, \cite[pg. 155]{Fay}.  It follows directly for $s>1$ that the comparable function $e^{-sd(p,q)}$ satisfies a mean value estimate.  Finally for a point in $\mathbb H$ given in polar form $re^{i\theta}$ the function $u(\theta)=1-\theta\cot \theta$ is a positive $2$-eigenfunction and satisfies the mean value property.  It follows directly that $\min\{u(\theta),u(\pi-\theta)\}$ and the comparable functions 
$\sin^2\theta= \operatorname{sech}^2\,d(p,i\mathbb R^+)$ and $e^{-2d(p,i\mathbb R^+)}$ satisfy mean value estimates. 

We will consider metric balls on the Riemann surface $R$ for the mean value estimate and use the injectivity radius to bound the covering number from the upper half plane.  From the description of collars and cusp regions the covering number for a fixed radius ball in $\mathbb H$ is bounded in terms of $inj^{-1}$ with constant determined by the radius of the ball.  We summarize the considerations for a Riemann surface $R=\mathbb H/\Gamma,\,\Gamma\subset PSL(2;\mathbb R)$ with the following.
\begin{lemma}
\label{enhmv}
Harmonic Beltrami differentials and the exponential-distance functions $e^{-2d(p,q)},\,\operatorname{sech}^2\,d(p,\tilde\alpha)$ and $e^{-2d(p,\tilde\alpha)}$ for a geodesic $\tilde\alpha$ satisfy a mean value estimate on $\mathbb H$ with constant determined by the radius of the ball.  A non negative function $f$ satisfying a mean value estimate 
on $\mathbb H$ and a subset of the discrete group $\mathcal G\subset\Gamma$ satisfy
\[
\sum_{A\in\mathcal G} f(A(p))\le c\,inj(p)^{-1}\int_{\cup_{A\in \mathcal G}A(\mathbf B(p,\rho))}f\,dA
\]
with constant determined by the mean value constant.
\end{lemma}
   
\subsection{Teichm\"{u}ller space and geodesic-length functions}

Let $\caT$ be the Teichm\"{u}ller space of genus $g$, $n$ punctured Riemann surfaces with hyperbolic metrics $ds^2$, \cite{Ahsome, Busbook, ImTan, Ngbook, Trmbook}.  From Kodaira-Spencer deformation theory the infinitesimal deformations of a surface $R$ are represented by the Beltrami differentials $\caH (R)$ harmonic with respect to the hyperbolic metric, \cite{Ahsome}. Also the cotangent space of $\caT$ at $R$ is $Q(R)$ the space of holomorphic quadratic differentials with at most simple poles at the punctures of $R$.  The holomorphic tangent-cotangent pairing is
\[
(\mu,\varphi)=\int_R \mu \varphi
\]
for $\mu\in\caH(R)$ and $\varphi\in Q(R)$.  Elements of $\caH (R)$ are symmetric tensors given as $\overline\varphi(ds^2)^{-1}$ for $\varphi\in Q(R)$ and $ds^2$ the hyperbolic metric.  The Weil-Petersson (WP) Hermitian metric and cometric are given as
\[
\langle\mu,\nu\rangle_{Herm}=\int_R\mu\overline\nu dA\quad\mbox{and}
\quad\langle\varphi,\psi\rangle_{Herm}=\int_R\varphi\overline\psi (ds^2)^{-1}
\]
for $\mu,\nu\in\caH (R)$ and $\varphi, \psi\in Q(R)$ and $dA$ the hyperbolic area element.  The WP Riemannian metric is $\langle\ ,\ \rangle =\Re\langle\ ,\ \rangle_{Herm}$.  The WP  metric is K\"{a}hler, non complete, with non pinched negative sectional curvature and determines a $CAT(0)$ geometry, see \cite{Ahsome, Zh2, Ngbook, Trmbook, Wlcomp, Wlpers} for references.  Ahlfors found for infinitesimal deformations defined by elements of $\caH(R)$ that the first derivatives of the WP metric tensor initially vanish, \cite{Ahsome}.  Equivalently a basis for 
$\caH(R)$ provides local coordinates for $\caT$ {\em normal} at the origin for the WP metric.  Equivalently to first order the WP Levi-Civita connection is initially Euclidean for deformations defined by elements of $\caH(R)$.

Basic invariants of a hyperbolic metric are the lengths of the unique geodesic representatives of the non peripheral free homotopy classes.  Points of the Teichm\"{u}ller space $\mathcal T$ are equivalence classes $\{(R,f)\}$ of marked Riemann surfaces with reference homeomorphisms $f:F\rightarrow R$ from a base surface $F$.  For the non peripheral free homotopy class $\alpha$ for $F$ the length of the geodesic representative for $f(\alpha)$ is the value of the geodesic-length $\ell_{\alpha}$ at the marked surface.  For  $R$ with uniformization representation $f_*:\pi_1(F)\rightarrow \Gamma,\,\Gamma\subset PSL(2;\mathbb R)$ and $\alpha$ corresponding to the conjugacy class of an element $A$ then $\cosh\ell_{\alpha}/2=\operatorname{tr} A/2$.   Traces of collections of elements provide real analytic coordinates for Teichm\"{u}ller space; for examples see \cite{Abbook, Busbook, ImTan, KnFr1,KnFr2,Thsurf}.     

We review the formulas for the differential of geodesic-length, the gradient of geodesic-length and the WP pairing.   Gardiner provided formulas for the differentials of geodesic-length and length of a measured geodesic lamination, \cite{Gardtheta, Gardmeas}.  In particular for a closed geodesic $\alpha$ conjugate the group $\Gamma$ for the geodesic to correspond to the deck transformation $A:z\rightarrow e^{\ell_{\alpha}}z$ and consider the series
\begin{equation}
\label{theta}
\Theta_{\alpha}=\sum_{B\in\langle A\rangle\backslash\Gamma}B^*(\frac{dz}{z})^2
\end{equation}
for $\langle A\rangle$ the cyclic group generated by $A$.  The differentials of geodesic-length are given as
\begin{equation}
\label{dgl}
d\ell_{\alpha}[\mu]=\frac{2}{\pi}\Re(\mu,\Theta_{\alpha})\quad\mbox{with}\quad \partial\ell_{\alpha}[\mu]=\frac{1}{\pi}(\mu,\Theta)
\end{equation}
for $\mu\in\mathcal H(\Gamma)$.  The WP Riemannian dual of the cotangent $\varphi\in Q(\Gamma)$ is the tangent 
$\overline{\varphi}(ds^2)^{-1}\in\mathcal H(\Gamma)$.   Accordingly the gradient of geodesic-length is $\grad\ell_{\alpha}=\frac{2}{\pi}\overline{\Theta_{\alpha}}(ds^2)^{-1}$.  The tangent vector $t_{\alpha}=\frac{i}{2}\grad\ell_{\alpha}$ has the geometric description as the infinitesimal Fenchel-Nielsen right twist deformation, \cite[display (1.1), Coro. 2.8, Thrm. 2.10]{WlFN}.  The Fenchel-Nielsen twist deformation and geodesic-length are related by duality in the WP K\"{a}hler form $\omega(\ ,t_{\alpha})=\frac12d\ell_{\alpha}$, 
\cite{WlFN, Wlcusps}.  The symplectic geometry of twists and lengths was considered in \cite[Thrm. 3.3]{Wlsymp} and \cite{Gdsymp} including the twist-length cosine formula
\[
4\omega(t_{\alpha},t_{\beta})= \langle \grad\ell_{\alpha},-i\grad\ell_{\beta}\rangle = 2\sum_{p\in\alpha\cap\beta}\cos \theta_p. 
\]
The sum is for transverse intersections of geodesics.  In particular the sum vanishes if the geodesics are disjoint or coincide.  In Section \ref{gradhess} we will consider Riera's length-length formula for the pairing $\langle \grad\ell_{\alpha},\grad\ell_{\beta}\rangle$ of geodesic-length gradients, \cite{Rier}

\section{Variations of geodesic-length}       

\subsection{Harmonic Beltrami differentials and Eichler integrals}    

The WP Riemannian metric is $\langle\ ,\ \rangle =\Re\langle\ ,\ \rangle_{Herm}$. To first order the WP Levi-Civita connection is initially Euclidean for deformations defined by elements of $\caH(R)$, \cite{Ahsome}.  The benefit of choosing $\caH(R)$ as the model for the tangent space is that expressions involving at most the initial first derivatives of the metric are WP intrinsic quantities.  For $t$ real, small and $\mu\in\caH(R)$, let $R^{t\mu}$ denote the deformation determined by the Beltrami differential $t\mu$.  For $z$ a local conformal coordinate for $R$ then  
$dw=dz+\mu\overline{dz}$ is the differential of a local conformal coordinate $w$  for $R^{t\mu}$.  At $t=0$ the family of Riemann surfaces $\{R^{t\mu}\}$ agrees to second order with a WP geodesic.  Relatedly for a smooth function defined in a neighborhood of $R$ in $\caT$ then $\frac{d^2}{dt^2}h(R^{tu})\mid_{t=0}$ is the WP Riemannian Hessian $\hess h$, \cite{BO'N, Wlnielsen}.  In general for vector fields $U,V$ the Hessian is given as $\hess h(U,V)=UVh-(D_UV)h$.  In particular for  $\mu,\nu\in\caH (R)$ and the family of Riemann surfaces $\{R^{\,t\mu+t'\nu}\}$  the WP covariant derivative $D_{\frac{\partial}{\partial t}}\frac{\partial}{\partial t'}$ vanishes at the origin for the tangent fields $\frac{\partial}{\partial t},\frac{\partial}{\partial t'}$.

A deformation of Riemann surfaces can be presented by a quasiconformal (qc) homeomorphism of the upper 
half plane $\mathH$.  For $R$ given as $\mathH/\Gamma$, $\Gamma\subset PSL(2;\mathR)$, a Beltrami differential $\mu\in\caH(R)$ with $\|\mu\|_{\infty}<1$ defines a deformation as follows.  There exists a qc self homeomorphism $w(z)$ of $\mathH$ satisfying $w_{\overline z}=\mu w_z$, \cite{AB}.  The map $w$ conjugates $\Gamma$ to a discrete group $w\circ\Gamma\circ w^{-1}\subset PSL(2;\mathR )$ uniformizing a homeomorphic Riemann surface $R^{\mu}$.  Homeomorphisms $w^{\mu}$ and $w^{\nu}$ are equivalent provided $w^{\mu}\circ (w^{\nu})^{-1}$ is homotopic to a conformal map.  Teichm\"{u}ller space is the space of equivalence classes of qc maps of $R$.  For $t$ real and small, we consider the family of qc homeomorphisms $w^{t\mu}$ and corresponding family of surfaces $\{R^{t\mu}\}$.  The initial deformation vector field  of holomorphic type $\dot w=\frac{d}{dt}w^{t\mu}\mid_{t=0}$ on $\mathH$ satisfies:
\begin{enumerate}
\item $\dot w_{\overline z}=\mu$,
\item\label{ES} the coboundary map $A\in \Gamma \rightarrow A_*(\dot w)-\dot w=\dot w(A(z))(A'(z))^{-1}-\dot w(z)$ is the class of the deformation in group cohomology with coefficients in real quadratic polynomials 
on $\mathH$, \cite{Ahsome, Eich, Weidis}. 
\end{enumerate} 
The standard approach presents $\dot w$ as the singular $\overline\partial$-potential integral of $\mu$, \cite{Ahsome, AB}.  

We  presented in \cite[Sec. 2]{Wlnielsen} for $\mu\in\caH (R)$ a formula for $\dot w$ in terms of a line integral.  The integral is a form of the Eichler-Shimura isomorphism of \ref{ES}. above.  
We introduced the line integral
\[
F(z)=\overline{\int^z_{z_0}(\overline z-\tau)^2\varphi(\tau)d\tau}
\]
for $\tau,z,z_0$ in $\mathH$, $z_0$ fixed, and $\varphi\in Q(\Gamma)$. The integral is path independent since the integrand is holomorphic.  The potential equation $F_{\overline z}=(z-\overline z)^2\overline{\varphi(z)}$ is an immediate property.  We established the following. 
\begin{theorem}
\label{lineint}
The deformation vector field for Beltrami differential $\mu=(z-\overline z)^2\overline{\varphi}$  is $$\dot w=\overline{\int^z_{z_0}(\overline z-\tau)^2\varphi(\tau)d\tau}\ +\ \int^z_{z_0}( z-\tau)^2\varphi(\tau)d\tau\ +\ q(z)$$
where $q$ is a quadratic polynomial with real coefficients. 
\end{theorem}
The polynomial corresponds to an infinitesimal $PSL(2;\mathR)$ conjugation.  We observe that for the quadratic differential $\varphi$ given by a suitably convergent series the integrals may be computed term-by-term.

We now specialize the formulas for the situation of $\Gamma$ conjugated to contain the transformation $z\rightarrow\lambda z,\ \lambda>1$ and the maps $w^{t\mu}$ fixing $0$ and infinity ($\dot w$ is $O(|z|)$ at the origin and $o(|z|^2)$ at infinity.)   The normalization provides that $\dot w(\lambda z)\lambda^{-1}-\dot w(z)$ is a multiple of $z$ (in fact $\dot\lambda z$ since $\lambda^t=w^{t\mu}(\lambda z)/w^{t\mu}(z)$ ) and consequently   
\begin{equation}
\label{ABC}
\dot w=z^22\Re\bA + z2\Re\bB+2\Re\bC + bz
\end{equation}
where
\begin{gather*}
\bA[\varphi]=\int^z_{z_0}\varphi\, d\tau\ +\ \frac{1}{1-\lambda}\int^z_{\lambda^{-1}z_0}\varphi\, d\tau \\
\bB[\varphi]=-2\int^z_{z_0}\tau\varphi\, d\tau\\
\bC[\varphi]=\int^z_{z_0}\tau^2\varphi\, d\tau\ +\ \frac{\lambda}{\lambda-1}\int^{z_0}_{\lambda^{-1}z_0}\tau^2\varphi\, d\tau,
\end{gather*}
\cite[see (4.1)]{Wlnielsen}.  The Eichler integral operators $\bA,\bB$ and $\bC$ define endomorphisms on the space $\operatorname{Holo}(\mathH)$  of functions holomorphic on $\mathH$ and the operators have formal properties of group cohomology.  Hejhal's variational formula $\dot\lambda= -2 \overline{\int^{z_0}_{\lambda^{-1}z_0}\tau \varphi d\tau}$  is a  consequence, \cite{Hejmono}.  The functions $\bA[\varphi]$ and $\bC[\varphi]$ are weighted equivariant for $z\rightarrow\lambda z$ acting on $\operatorname{Holo}(\mathH)$ by composition: $A(\lambda z)=\lambda^{-1}A(z)$ and $C(\lambda z)=\lambda C(z)$.  

\subsection{The deformation vector field on the horizontal strip}
\label{defvf}
It is natural to introduce a model for the hyperbolic plane presenting the transformation $z\rightarrow\lambda z$ as a translation.  The infinite horizontal strip $\mathS$ is a suitable model.  We present formulas for the deformation vector field on $\mathS$.  We then apply the formulas to study the first and second variation of geodesic-length.

The upper half plane $\mathH$ and infinite horizontal strip $\mathS=\{\zeta \mid 0< \Im \zeta<\pi\}$ are related by the map $\zeta=\log z$.  We consider now for $\mu\in\caH (\Gamma)$ the family $w^{t\mu}$ of qc homeomorphisms of $\mathH$ satisfying $w^{t\mu}_{\overline z}=t\mu w^{t\mu}_{z}$, $w^{t\mu}$ fixing $0,1$ and infinity, \cite{AB}. The homeomorphism $w^{t\mu}$ conjugates the transformation  $z\rightarrow\lambda z$ to a transformation fixing $0$ and infinity  with multiplier $\lambda^t$.   The family $w^{t\mu}$ is conjugated by $\zeta=\log z$ to a family of qc homeomorphisms $g^t$ of $\mathS$.  The homeomorphism $g^t$ conjugates the translation $\zeta\rightarrow\zeta+\ell$, $\ell=\log \lambda$, to a translation by $\ell^t=\log \lambda^t$.  We have from   the translation equivariance equation
$$
g^t(\zeta+\ell)=g^t(\zeta)\ + \ \ell^t,
$$
the expansion in $t$
$$
g^t(\zeta)=\zeta+t\dot g(\zeta)+\frac{t^2}{2}\Ddot g(\zeta)+O(t^3)
$$
and from the Beltrami equation the perturbation relation
$$
t\dot g_{\overline \zeta}+\frac{t^2}{2}\Ddot g_{\overline \zeta}=t\mu(1+t\dot g_{\zeta})
$$
provides that $\dot g_{\overline \zeta}=\mu$ and $\Ddot g_{\overline\zeta}=2\mu\dot g_{\zeta}$, \cite{AB}.  We choose $\caF=\{\zeta\in\mathS \mid 0\le\Re \zeta<\ell\}$ as a fundamental domain for the translation.  The following presents the first and second variation of $\ell^t$ as integrals of the deformation vector field $\dot g$.
\begin{lemma}
\label{ellvar}
For $\mu\in\caH (\Gamma)$ and $g^t$ the family of qc homeomorphisms of the strip $\mathS$ then
$$
\dot\ell=\frac{2}{\pi}\Re\int_{\caF}\mu\,dE\quad \mbox{and}\quad \Ddot\ell=
\frac{4}{\pi}\Re\int_{\caF}\mu\dot g_{\zeta}\,dE
$$
for $dE$ the Euclidean area element.  
\end{lemma}
\prf Ahlfors and Bers showed that the solution $g^t$ of the Beltrami equation varies real analytically in the Banach space with norm the sum of the H\"{o}lder norm and $L^p$ norm of the first distributional derivatives, \cite{AB}. For smooth Beltrami differentials it follows immediately from elliptic regularity that the $t$-derivatives of $g^t$ are smooth, \cite[Chap. 6]{GT}.  

We write $\zeta=x+iy$, $g=u+iv$ and begin with
$$
\ell^t=u^t(\ell,y)-u^t(0,y)=\int_0^{\ell}u^t_x\,dx
$$
and for the $n^{th}$ derivative at $t=0$ 
$$
\ell^{(n)}=\int_0^{\ell}u_x^{(n)}dx \quad\mbox{and}\quad \ell^{(n)}=\frac{1}{\pi}\int_{\caF}u_x^{(n)}dx.
$$
The map $g^t$ preserves the boundary of $\mathS$ and thus $\Im g^t$ is constant on the boundary while $\Im g^{(n)}$ vanishes on the boundary.  It follows that
$$
0=v^{(n)}(x,\pi)-v^{(n)}(x,0)=\int_0^{\pi}v_y^{(n)}\,dy.
$$
The general formula
$$
\ell^{(n)}=\frac{2}{\pi}\Re\int_{\caF}\frac{\partial}{\partial\overline\zeta}g^{(n)}\,dE
$$
is now a consequence and the argument is completed with the relation 
$g^{(2)}_{\overline\zeta}=2\mu g^{(1)}_{\zeta}$.  The proof is complete.   

Gardiner's formula (\ref{dgl}) for the differential of geodesic-length follows since for $\zeta=\log z$ then $(d\zeta)^2=(\frac{dz}{z})^2$  and $\caF$ maps to the half annulus $\mathcal A=\{1\le |z|<\lambda, \Im z> 0\}\subset\mathH$, 
\cite{Gardtheta}.  As always the pairing $\int_{\mathcal A} \mu (\frac{dz}{z})^2$ can be {\em folded} to $\int_{\mathcal F_{\Gamma}}\mu\sum_{A\in\Gamma_0\backslash\Gamma}A^*(\frac{dz}{z})^2$ where $\mathcal F_{\Gamma}$ is now a fundamental domain for the group and $\Gamma_0$ is the cyclic subgroup generated by $z\rightarrow\lambda z$.  

We are ready to calculate the deformation vector field quantity $\dot g_{\zeta}$ on $\mathS$.  Vector fields on the upper half plane and the strip are related for $\zeta=\log z$ by $\frac{\partial}{\partial \zeta}=z\frac{\partial}{\partial z}$ and in particular $\dot w \frac{\partial}{\partial z}=\dot w z^{-1}\frac{\partial}{\partial \zeta}=\dot g \frac{\partial}{\partial \zeta}$.  And we have for the coefficients $\frac{\partial}{\partial \zeta}\dot g=z\frac{\partial}{\partial z}(\dot wz^{-1})$. Lemma \ref{ellvar} above and the change of variables $\zeta=\log z$  combine to give Theorem 3.2 of \cite{Wlnielsen} $\Ddot \ell=\frac{4}{\pi}\Re\int_{\mathcal A}\mu z^{-1}\frac{\partial}{\partial z}(\dot w z^{-1})\,dE$.

We are ready to consider $\dot g_{\zeta}$ and the term-by-term contributions for a series expansion of the quadratic differential $\varphi$.  We begin with (\ref{ABC})
\[
\dot w=z^22\Re\bA+z2\Re\bB+2\Re\bC + bz
\]
and find
\begin{equation}
\label{zetaform}
\dot g_{\zeta}=z\frac{\partial}{\partial z} (\dot wz^{-1})=z2\Re\bA-z^{-1}2\Re\bC
\end{equation} 
since $\bA, \bB, \bC \in \operatorname{Holo}(\mathH)$ with $\bA'=\varphi,\ \bB'=-2z\varphi$ and $\bC'=z^2\varphi$.  We now find the contribution of an individual term $z^{\alpha}(\frac{dz}{z})^2$ for $\alpha=\epsilon n$, $\epsilon=\frac{2\pi i}{\ell}, \ell=\log \lambda$ and $n\in \mathZ$
\[
\bA[z^{\alpha}(\frac{dz}{z})^2]=\frac{1}{\alpha-1}z^{\alpha-1}\quad\mbox{and}\quad\bC[z^{\alpha}(\frac{dz}{z})^2]=\frac{1}{\alpha+1}z^{\alpha+1}
\]
\begin{equation}
\label{alpha}
z\bA[z^{\alpha}(\frac{dz}{z})^2]-z^{-1}\bC[z^{\alpha}(\frac{dz}{z})^2]=(\frac{1}{\alpha-1}-\frac{1}{\alpha+1})z^{\alpha}=\frac{-2}{1+|\alpha|^2}z^{\alpha}
\end{equation}
and for $\zeta=x+iy=\log z$
\begin{multline}
\label{alphabar}
z\overline{\bA[z^{\alpha}(\frac{dz}{z})^2]}-z^{-1}\overline{\bC[z^{\alpha}(\frac{dz}{z})^2]}\\=(\frac{e^{2iy}}{\alpha-1}-\frac{e^{-2iy}}{\alpha+1})\overline{z^{\alpha}}
=\frac{-2\cos {2y}-2i\alpha\sin{2y}}{1+|\alpha|^2}\overline{z^{\alpha}}.
\end{multline}
The right hand sides of (\ref{alpha}) and (\ref{alphabar}) will be paired below with $\varphi$.  

The quadratic differential $\varphi$ on $\mathH$ is invariant with respect to $z\rightarrow\lambda z$ and admits the expansion
\begin{equation}
\label{zA}
\varphi=\sum_na_nz^{\epsilon n}(\frac{dz}{z})^2\quad\mbox{with integral}\quad z\bA[\varphi]=\sum_n\frac{a_nz^{\epsilon n}}{\epsilon n -1}.
\end{equation}
We can consider $z\bA[\varphi]$ and $z^{-1}\bC[\varphi]$ as $0$-tensors since the quantities are coefficients for the $\zeta$-derivative of the vector field $\dot g$.  The $0$-tensors are expressed on $\mathS$ by a change of variable.
  We introduce a Hermitian pairing for suitable $\ell$-translation invariant functions on $\mathS$
\[
\caQ(\beta,\delta)=\int_{\caF}\beta\overline{\delta}\,\sin^2y\,dE
\]
with $\mathcal F$ a fundamental domain for $\ell$-translation.
We observe for $\zeta=\log z$, $z\in\mathH,\ \zeta\in\mathS$ the Euclidean area elements satisfy $\frac{dE(z)}{|z|^2}=dE(\zeta)$ with $\Im z=|z|\sin{\Im \zeta}$.  The pairing is also given as
\[
\caQ(\beta,\delta)=\int_{\mathcal A}\beta\overline{\delta}\,\frac{(\Im z)^2}{|z|^4}\,dE
\]
for $\mathcal A$ the image of $\caF$ by $z=e^{\zeta}$ and $\beta,\delta$ the compositions with $\zeta=\log z$.  We include formula (\ref{dgl}) in the following statement of variational formulas.  
\begin{theorem}
\label{quadform}
The first variation of geodesic-length is given as
\[
\dot \ell[\mu]=-4\ell\,\Re a_0
\]
and the second variation as
\[
\Ddot\ell[\mu,\mu]=\frac{32}{\pi}\,\Re\caQ(z\bA,z\bA)\ - \ \frac{16}{\pi}\,\Re\caQ(z\bA,\overline{z\bA})
\]
for Hermitian pairing $\caQ$ and a harmonic Beltrami differential $\mu=(z-\overline z)^2\overline\varphi$, $\varphi\in Q(\Gamma)$, with zeroth coefficient $a_0$, and integral $\bA=\bA[\varphi]$.  
\end{theorem}
\prf We combine Lemma \ref{ellvar} with formula (\ref{zetaform}) to find for $\mu=-4\sin^2y\,\overline\varphi$ and $\zeta=x+iy$ that
\[\Ddot\ell=\frac{-16}{\pi}\Re\int_{\caF}\overline\varphi((z\bA-z^{-1}\bC)+(z\overline\bA-z^{-1}\overline\bC))\,\sin^2y\,dE.
\]
The Beltrami differential $\mu$ is bounded and as noted Ahlfors and Bers established that $\dot g_{\zeta}$ is in $L^p$.  The series expansion for $\varphi,\,\bA$ and $\bC$ converge uniformly on compacta and consequently we can consider term-by-term contributions.  
	
The orthogonality of exponentials for $x$-integration for $x+iy=\zeta=\log z$ provides  simplifications.  The weight $n$ term $\overline{z^{\epsilon n}}$, $\epsilon=\frac{2\pi i}{\ell}$,  in the expansion for $\overline\varphi$ only has a nonzero $x$-integral when paired with the weight $n$ term of $(z\bA-z^{-1}\bC)$.  The contribution from (\ref{alpha}) is
\[
-2\int_{\caF}\overline{(\frac{z^{\epsilon n}}{\epsilon n-1})} (\frac{z^{\epsilon n}}{\epsilon n-1})\,\sin^2y\,dE.
\]
The weight $n$ term $\overline{z^{\epsilon n}}$ in the expansion for $\overline\varphi$ only has a nonzero $x$-integral when paired with the weight $-n$ term of $(z\overline\bA-z^{-1}\overline\bC)$.  The contribution from (\ref{alphabar}) is
\[
-2\int_{\caF}\overline{z^{\epsilon n}}\bigl(\frac{\cos{2y}+\frac{2\pi n}{\ell}\sin{2y}}{1+(\frac{2\pi n}{\ell})^2}
\bigr)\overline{z^{-\epsilon n}}\,\sin^2y\,dE
\]
and with trigonometric integrals in $y$
\[
=\int_{\caF}\overline{(\frac{z^{\epsilon n}}{\epsilon n-1})\,(\frac{z^{-\epsilon n}}{-\epsilon n-1})}\,\sin^2y\,dE.
\]
The desired conclusion now follows from the expansion (\ref{zA}).  The proof is complete. 

We first consider immediate consequences for the geometry of Teichm\"{u}ller space.  A form of the following based on \cite{Wlnielsen} was given by Yeung, \cite[Proposition 1]{Yun1}.  
\begin{corollary}
\label{first}
The complex Hessian of geodesic-length is given as 
\[
\partial\overline\partial\ell[\mu,\mu]=\frac{16}{\pi}\caQ(z\bA,z\bA)
\] 
for $\mu=(z-\overline z)^2\overline\varphi,\ \varphi\in Q(\Gamma)$ and integral $\bA=\bA[\varphi]$.  The complex and real Hessians are uniformly comparable 
\[
\partial\overline\partial\ell\le\hess\ell\le3\partial\overline\partial\ell.
\]  
The first and second derivatives of geodesic-length satisfy 
\[
2\ell\Ddot\ell[\mu,\mu]-\dot\ell^2[\mu]-3\dot\ell^2[i\mu]\ge0\quad\mbox{and}\quad\ell\partial\overline\partial\ell-2\partial\ell\overline\partial\ell\ge 0
\] 
with equality only for the elementary $\varphi=a(\frac{dz}{z})^2$.  
\end{corollary}
\prf The first two statements are consequences of the observation that $\caQ(z\bA,z\bA)$ is a Hermitian form, $\caQ(z\bA,\overline{z\bA})$ is a complex bilinear form and the inequality $|\caQ(z\bA,\overline{z\bA})|\le \caQ(z\bA,z\bA)$. 
From the given formula the contribution  of the weight zero term $a_0(\frac{dz}{z})^2$ to $\Ddot\ell$ is 
$8\ell(\Re a_0)^2+24\ell(\Im a_0)^2$ and to $\partial\overline\partial\ell$ is $8\ell|a_0|^2$.  
The inequality for derivatives follows from  formula (\ref{dgl}) and that $\Ddot\ell$ is a positive definite quadratic form and $\partial\overline\partial\ell$ a positive definite Hermitian form on the coefficients of $\varphi$.  The proof is complete.

\begin{example} Convexity for geodesic-length functions.
\end{example}
For closed curves $\gamma_1,\dots,\gamma_n$ the geodesic-length sum $\ell_{\gamma_1}+\cdots+\ell_{\gamma_n}$ satisfies:   $(\ell_{\gamma_1}+\cdots+\ell_{\gamma_n})^{1/2}$ is strictly convex along WP geodesics, $\log (\ell_{\gamma_1}+\cdots+\ell_{\gamma_n})$ is strictly plurisubharmonic, and $(\ell_{\gamma_1}+\cdots+\ell_{\gamma_n})^{-1}$ is strictly plurisuperharmonic.

We consider the Teichm\"{u}ller space of the annulus as an example that foreshadows the behavior of the WP metric for surfaces with short geodesics.
\begin{example}
The Teichm\"{u}ller space of the annulus.
\end{example}
The annulus with geodesic of length $\ell$ is given as the quotient $\mathH/\langle A\rangle$ for the cyclic group generated by $z\rightarrow e^{\ell} z$.  Geodesic-length provides a global coordinate for the real one-dimensional Teichm\"{u}ller space.  The space $Q(\langle A\rangle)$ of quadratic differentials is the span over $\mathbb R$ of $(\frac{dz}{z})^2$.  The norm of $\mu=(\Im z)^2\overline{(\frac{dz}{z})^2}$ is $\langle\mu,\mu\rangle_{Herm}=\int_{\mathcal A}\mu\overline\mu\,dA=\frac{\pi}{2}\ell$ for $\mathcal A=\{1<|z|<\lambda, \Im z>0\}$.  From formula (\ref{dgl}) we observe that $\grad\ell^{1/2}=\frac{1}{\pi\ell^{1/2}}(\Im z)^2\overline{(\frac{dz}{z})^2}$ with $\langle\grad(2\pi\ell)^{1/2},\grad(2\pi\ell)^{1/2}\rangle=1$.  The parameter $(2\pi\ell)^{1/2}$ is the WP unit-speed parameterization of the Teichm\"{u}ller space.  Correspondingly from Corollary \ref{first} and $\dot\ell[i\mu]=0$ the Hessian of $(2\pi\ell)^{1/2}$ vanishes and the vector field $\grad(2\pi\ell)^{1/2}$ is WP parallel, \cite{BO'N}.  

\subsection{Bounding the Hessian of geodesic-length}
\label{bdhess}
We bound the Hermitian form $\mathcal Q(z\bA,z\bA)$ in terms of the WP pairing and a geometrically defined weight function.  The weight function given in the universal cover is the sum of the inverse square of exponential-distance from a point to the components of the lift of the geodesic.  The description provides the approach for considering the weight function of a measured geodesic lamination. The weight is the corresponding integral on the space of leaves of the lamination for the distance from the point to the leaves of the lamination.   Gardiner established that the first variation formula for geodesic-length generalizes to a formula for the length of a measured geodesic lamination \cite{Gardmeas}.  Overall we find that the Hessian of length of a measured geodesic lamination is comparable to the product of the length and the WP pairing with constants depending on the reduced injectivity radius of the surface.  The length is always strictly convex along WP geodesics.

Considerations for bounding the Hessian begin with the horizontal strip model $\mathbb S$, the integral $\bA[\varphi]$ and a second Hermitian pairing for suitable $\ell$-translation invariant functions on $\mathbb S$
\[
\mathcal{QS}(\beta,\delta)=\int_{\mathcal F}\beta\overline{\delta}\sin ^4y\,dE
\]
for $\mathcal F$ a fundamental domain for $\ell$-translation.
\begin{lemma}
\label{QS}
The Hermitian forms $\mathcal Q$ and $\mathcal{QS}$ are bounded as
\[
\mathcal{QS}(\varphi,\varphi)\le3\mathcal{Q}(\bA[\varphi],\bA[\varphi])\le4\mathcal{QS}(\varphi,\varphi)
\]
for a holomorphic quadratic differential $\varphi\in Q(\Gamma)$ with integral $\bA[\varphi]$.
\end{lemma}
\prf The quantities $\varphi$ and $z\bA[\varphi]$ considered on $\mathbb S$ are given as sums (\ref{zA}) of exponentials $z^{\epsilon n}=e^{\epsilon n(x+iy)},\ x+iy=\log z,\ \epsilon=\frac{2\pi i}{\ell}$.  The orthogonality of exponentials for $x$-integration provides that $\mathcal Q$ and $\mathcal{QS}$ are diagonalized by the expansions.  It is sufficient to consider the contribution of individual terms.  The weight $n$ term 
$z^{\epsilon n}(\frac{dz}{z})^2$ contributes to $\mathcal Q$ the term
\[
\frac{1}{|\epsilon n|^2+1}\int_{\mathcal F}|z^{\epsilon n}|^2\sin ^2y\,dxdy\quad=\quad\frac{\ell}{|\epsilon n|^2+1}\int^{\pi}_0e^{-2|\epsilon|ny}\sin^2y\,dy
\]
where for $a=|\epsilon|n$ by \cite[pg. 12, formula (29)]{ErdHI}
\[
=\frac{\ell}{4}\frac{(1-e^{-2\pi a})}{a(a^2+1)^2}.
\]
The weight $n$ term $z^{\epsilon n}(\frac{dz}{z})^2$ contributes to $\mathcal{QS}$ the term
\[
\int_{\mathcal F}|z^{\epsilon n}|^2\sin ^4y\,dxdy\quad=\quad\ell\int^{\pi}_0e^{-2|\epsilon|ny}\sin^4y\,dy
\]
where  by \cite[pg. 12, formula (29)]{ErdHI}
\[
=\frac{3\ell}{4}\frac{(1-e^{-2\pi a})}{a(a^2+1)(a^2+4)}.
\]
The desired inequalities now follow from the observation that $(1-e^{-2\pi a})a^{-1}$ is positive for all real $a$ and the elementary inequality $(a^2+4)^{-1}\le(a^2+1)^{-1}\le4(a^2+4)^{-1}$.  The proof is complete.  

We next express the pairing $\mathcal{QS}$ as an integral on $\mathbb H$ with the weight of inverse square exponential-distance to the given geodesic.  A geodesic on the hyperbolic plane determines a projection of the plane to the geodesic.  We write $d(\ ,\ )$ for hyperbolic distance and $dA$ for the area element of the hyperbolic plane.
\begin{lemma}
\label{QE}
The Hessian of geodesic-length is bounded as
\[
\mathcal{QE}(\mu,\mu)\le3\pi\partial\overline{\partial}\ell_{\gamma}[\mu,\mu]\le 16\mathcal{QE}(\mu,\mu)
\]
for
\[
\mathcal{QE}(\mu,\mu)=\int_{\mathcal A}|\mu|^2e^{-2d(\ ,\tilde\gamma)}\,dA
\]
for $\gamma$ the geodesic corresponding to the deck transformation $z\rightarrow\lambda z$ with axis $\tilde\gamma$ with fundamental domain $\mathcal A=\{1\le|z|<\lambda,\Im z>0\}$ and $\mu\in\mathcal H(R)$.  
\end{lemma}
\prf We begin with the Hermitian form $\mathcal{QS}$ given for the strip $\mathbb S$ and the change of variables to $\mathbb H$.  As already noted for $\zeta=\log z$, we have that $ \frac{dE(z)}{|z|^2}=dE(\zeta)$ and $\Im z=|z|\sin\Im\zeta$.  We have that
\[
\mathcal {QS}=\int_{\mathcal A}\,\bigl|\sum_na_nz^{\epsilon n-2}\bigr|^2\frac{(\Im z)^4}{|z|^2}\,dE
\]
and recall that the Beltrami differential $\mu=-4(\Im z)^2\overline{\varphi}$ to obtain
\[
\label{QS1}
\mathcal{QS}=\int_{\mathcal A}|\varphi|^2(\frac{\Im z}{|z|})^2\,(ds^2)^{-1}=\frac{1}{16}\int_{\mathcal A}|\mu|^2(\frac{\Im z}{|z|})^2\,dA.
\]
The next matter is to relate $(\frac{\Im z}{|z|})^2$ to hyperbolic distance.  For the variable $z$ of $\mathbb H$ given in polar form $re^{i\theta}$ the hyperbolic distance to the imaginary axis, the axis of the transformation, is $\log|\csc\theta+|\cot\theta||$ and $e^{-2d(\ ,i\mathbb R^+)}\le\sin^2\theta=(\frac{\Im z}{|z|})^2\le 4e^{-2d(\ ,i\mathbb R^+)}$.  The bound for the Hessian now follows from Corollary \ref{first}, Lemma \ref{QS} and the inequality.  The proof is complete.  

We now proceed to obtain a general description of the Hessian in terms of a weight function $\mathbb P$ and the WP pairing.  The above description provides an approach for both closed geodesics and measured geodesic laminations.  Convexity will be bounded in terms of the maximum and minimum of the weight function $\mathbb P$ on the complement of horoball neighborhoods of any possible cusps.  We continue with considerations for hyperbolic geometry.

A geodesic $\tilde\gamma$ on the hyperbolic plane determines a projection and a distance function 
$d(\ ,\tilde\gamma)$.  Both the projection and distance considered in the compact-open topology on the plane depend continuously on the geodesic $\tilde\gamma$.   A Riemann surface $R$ with cusps contains disjoint area two standard horoball neighborhoods of the cusps, see Theorem \ref{collars} or \cite{Leut}.  The {\em reduced surface} $R_{red}$ is defined as the complement in $R$ of the unit area horoball neighborhoods of cusps.  The injectivity radius $inj_{red}$ and diameter $diam_{red}$ of $R_{red}$ are finite positive and are determined from the lengths of closed geodesics.  Geodesics on the upper half plane $\mathbb H$ are determined by their endpoints.  For the upper half plane with boundary $\check{\mathbb R}=\mathbb R\cup\{\infty\}$ the space of unoriented complete geodesics is given as the M\"{o}bius band $\mathcal G=\check{\mathbb R}\times\check{\mathbb R}\backslash\{diagonal\}/\{interchange\}$.  

A {\em geodesic lamination} for $R=\mathbb H/\Gamma$ is a closed subset which is foliated by complete geodesics (see \cite{Bonmgl} for background on measured geodesic laminations; the following is summarized from Bonahon's presentation.)  A {\em measured geodesic lamination} $\upsilon$ with compact support is a compact geodesic lamination $\lambda(\upsilon)$ together with a transverse invariant positive measure with support $\lambda(\upsilon)$.   A measured geodesic lamination can be considered as a locally finite Borel measure on the space of unoriented geodesics on $\mathbb H$ with the measure $\Gamma$ invariant.  The space $\mathcal{ML}_0$ of measured geodesic laminations with compact support is equipped with the weak$^*$-topology for the space of invariant measures. An arc is generic with respect to simple geodesics if it is transverse to every simple geodesic of $R$.  Every arc can be approximated by generic arcs.  For a generic arc $\alpha$ each measured geodesic lamination defines a measure on $\alpha$: for $h$ a continuous function on $\alpha$ and $\upsilon$ a measured geodesic lamination associate the integral $\upsilon(h)=\int_{\alpha}h\,d\upsilon$.  A finite system $\alpha_1,\dots,\alpha_k$ of generic arcs exists such that the map of $\mathcal{ML}_0$ to $\mathbb R^k$ defined by $\upsilon\rightarrow(\upsilon(\alpha_j))$ is a homeomorphism to its (piecewise linear) image.   The weighted simple closed curves form a dense subset of $\mathcal{ML}_0$.  From Lemma \ref{separ} measured geodesic laminations  are disjoint from the area two horoball neighborhoods of cusps.

The {\em length} $\ell_{\upsilon}$ of a measured geodesic lamination $\upsilon$ is defined in terms of a finite family of generic arcs $\alpha_1,\dots,\alpha_k$, \cite{Bonmgl}.  The components of $\lambda(\upsilon)-\cup_{j=1}^k\alpha_j$ are geodesic segments of finite length characterized be either of their endpoints in $\cup_{j=1}^k\alpha_j$.  The measure defined by $\upsilon$ on $\cup_{j=1}^k\alpha_j$ provides a measure on the space of components of $\lambda(\upsilon)-\cup_{j=1}^k\alpha_j$.  The length $\ell_{\upsilon}$ is defined as the integral of the segment lengths with respect to the measure defined by $\upsilon$.  Length is a continuous function on $\mathcal{ML}_0$.

We are ready to introduce the $\Gamma$-invariant weight function on $\mathbb H$ for a measured geodesic lamination.
\begin{definition}
\label{P}
The inverse square exponential-distance is defined for a measured geodesic lamination $\upsilon$ as
\[
\mathbb P_{\upsilon}(p)=\int_{\mathcal G}e^{-2d(p,\beta)}\,d\upsilon(\beta).
\]
\end{definition} 

For the case of $\upsilon$ given by a closed geodesic $\gamma$ on $R$ corresponding to a deck transformation $A$ with axis $\tilde\gamma$ the integral can be expressed as a sum.  The support of the lamination in $\mathbb H$ is the union $\cup_{B\in\langle A\rangle\backslash\Gamma}B^{-1}(\tilde\gamma)$.  Provided convergence the integral is the sum of $e^{-2d(p,B^{-1}(\tilde\gamma))}=e^{-2d(B(p),\tilde\gamma)}$ over the union.  In particular for $m(\gamma)$ the mass of a single intersection with $\gamma$ the integral over $\mathcal G$ is $m(\gamma)\mathbb P_{\gamma}$ for the series
\[
\label{Pseries}
\mathbb P_{\gamma}(p)=\sum_{B\in\langle A\rangle\backslash\Gamma}e^{-2d(B(p),\check\gamma)}.
\]
\begin{lemma}
\label{Pbound1}
The inverse square exponential-distance is bounded in terms of the length $\ell_{\upsilon}$ of the measured geodesic lamination and the injectivity radius of the reduced surface:   $c_1(inj_{red})\ell_{\upsilon}\le \mathbb P_{\upsilon}$ on $R_{red}$ and $\mathbb P_{\upsilon}\le c_2(inj_{red})\ell_{\upsilon}$ on $R$.  Furthermore the map $\upsilon\rightarrow\mathbb P_{\upsilon}$ from $\mathcal{ML}_0$ to $C(R)$ is continuous for the compact-open topology on $C(R)$.
\end{lemma}
\prf  First lift $\upsilon$ to the upper half plane and consider for $q\in\mathbb H$ the intersection $\Sigma(q)=\bB(q,2\,diam_{red})\cap\upsilon$ for $\bB$ the hyperbolic ball and consider the local contribution to length $L(q)=\int_{\Sigma(q)}\upsilon\times d\ell$. The function $L$ is $\Gamma$-invariant.  Basic properties of $\Sigma$ and $L$ are as follows.  For $q\in\widetilde{R_{red}}$, the lift of $R_{red}$, then $\Sigma(q)$ contains in $\overline{\bB(q,diam_{red})}$ a representative of each point of $R_{red}$ and for $q$ in a horoball neighborhood of area $e^{-2\,diam_{red}}$ of a cusp then $\Sigma(q)$ is empty.  It follows from the characterization of the half-radius ball $\mathcal B(q)=\overline{\bB(q,diam_{red})}$ that given a point on a leaf of the lamination $\upsilon\subset\mathbb H$ there exists a $B\in\Gamma$ such that $\Sigma(B(q))$ contains the point and a segment of the leaf of length at least $2\,diam_{red}$.  In particular $L(q)/2\,diam_{red}$ is an upper bound for the mass of the leaves intersecting $\mathcal B(q)$.

We compare the integral $\mathbb P$ over $\mathcal G$ to a sum over the $\Gamma$ orbit of the point $p$ in $\mathbb H$.  The comparison is in terms of groupings of projected images of $p$ to the leaves of $\upsilon$.  For a reference point $q_0\in R_{red}$ the projected images are contained in the union $\cup_{B\in\Gamma}\mathcal B(B(q_0))$ with the projected points of $\mathcal B(B(q_0))$ at distance at least $d(q_0,B(q_0))-2\,diam_{red}$ from $p$.  As noted the leaves of $\upsilon$ intersecting $\mathcal B(B(q_0))$ have mass at most $L(q_0)/2\,diam_{red}$.  It now follows that the projected points contained in $\mathcal B(B(q_0))$ contribute at most $e^{-2(d(q_0,B(q_0))-2\,diam_{red})}L(q_0)/2\,diam_{red}$ to the integral for $\mathbb P_{\upsilon}$.  The next step is to bound $L(q_0)$.  The covering number for the map from $\bB(q_0,2\,diam_{red})\rightarrow R_{red}$ is bounded in terms of $diam_{red}$ and $inj_{red}$.  Since the diameter is bounded in terms of the injectivity radius, $L(q_0)$ is bounded above as $c(inj_{red})\ell_{\upsilon}$.  Next from Lemma \ref{enhmv} we have that $e^{-2d(q_0,z)}$ satisfies a mean value estimate and consequently
\[
\sum_{B\in\Gamma}e^{-2d(q_0,B(q_0))}\le c(inj_{red})\int_{\bB(B(q_0),inj_{red})}e^{-2d(p,z)}\,dA.
\]  
Gathering the bounds we have that $\mathbb P_{\upsilon}$ is bounded by $c(inj_{red})\ell_{\upsilon}\int_{\mathbb H}e^{-2d(p,z)}\,dA$.   For the unit disc model  for the hyperbolic plane $ds^2=(\frac{2|dz|}{(1-|z|^2)})^2$ and $d(0,z)=\log\frac{1+|z|}{1-|z|}$.  The comparison integral converges and the upper bound for $\mathbb P_{\upsilon}$ is established.  

The lower bound is obtained by considering the projected points in a single ball.  Consider 
$p\in\widetilde{R_{red}}$.  The points of $\bB(p,2\,diam_{red})$ are at distance at most $2\,diam_{red}$ to $p$ and thus $e^{-2d(p,\beta)}\ge e^{-2\,diam_{red}}$.  Leaves of the lamination in $\bB(p,2\,diam_{red})$ have length at most $4\,diam_{red}$ and thus the mass of the leaves intersecting the ball is at least $L(p)/4\,diam_{red}$.  The half-radius ball $\mathcal B$ surjects to $R_{red}$ which contains the lamination and thus $L(p)>\ell_{\upsilon}$.  The contribution of projected points in $\bB(p,2\,diam_{red})$ to $\mathbb P_{\upsilon}$ is bounded below by $c'(inj_{red})\ell_{\upsilon}$.   The lower bound is established.

The last matter is to consider the continuity of $\upsilon\rightarrow\mathbb P_{\upsilon}$.  It is equivalent to show that the map is continuous to $C(\mathbb H)$ with the compact-open topology.   Given a compact set $K$ in $\mathbb H$ the contribution to $\mathbb P_{\upsilon}$ from projected points in the exterior of $\bB(K,\rho)$, $\rho$ large, is bounded from the above by the tail of a convergent integral and $\ell_{\upsilon}$.  For $\rho$ sufficiently large the exterior contribution is sufficiently small.  The remaining interior integral is supported on a compact set and is straightforward to bound.  The proof is complete.

Length functions of measured geodesic laminations are real analytic and Gardiner has provided a formula for the differential $d\ell_{\upsilon}$, \cite{Gardmeas}.  The formula provides that the map $\upsilon\rightarrow d\ell_{\upsilon}$ from $\mathcal{ML}_0$ to $Q(R)$ is continuous.

\begin{theorem}
\label{hessmgl}
The Hessian of length of a measured geodesic lamination $\upsilon$ is bounded in terms of the WP pairing and the weight $\mathbb P_{\upsilon}$ as follows
\[
\langle\mu,\mu\mathbb{P}_{\upsilon}\rangle \le 3\pi\partial\overline{\partial}\ell_{\upsilon}[\mu,\mu]\le 16\langle\mu,\mu\mathbb{P}_{\upsilon}\rangle 
\]
and
\[
\langle\mu,\mu\mathbb{P}_{\upsilon}\rangle \le 3\pi\Ddot\ell_{\upsilon}[\mu,\mu]\le 48\langle\mu,\mu\mathbb{P}_{\upsilon}\rangle 
\]
for $\mu$ a harmonic Beltrami differential.  There are positive functions $c_1$ and $c_2$ such that
\[
c_1(inj_{red})\ell_{\upsilon}\langle\ ,\ \rangle 
\le\partial\overline{\partial}\ell_{\upsilon},\ \Ddot\ell_{\upsilon}\le c_2(inj_{red})\ell_{\upsilon}
\langle\ ,\ \rangle 
\]
for the reduced injectivity radius with $c_1(\rho)$ an increasing function vanishing at the origin and $c_2(\rho)$ a decreasing function tending to infinity at the origin.  
\end{theorem}
A bound for $\mathbb P_{\upsilon}$ independent of injectivity radius is provided in Lemma \ref{Pbound2} below.

{\bf Proof of Theorem.}  The matter is to bound the Hessian in terms of the WP pairing with weight $\mathbb P_{\upsilon}$.  From Corollary \ref{first} we can consider the real or complex Hessian.  From Lemma \ref{QE} for a closed geodesic $\gamma$ with mass $m(\gamma)$ the Hessian is bounded by the integral $\int_{\mathcal A}|\mu|^2e^{-2d(\ ,\tilde\gamma)}\,dA$ for the geodesic $\gamma$ corresponding to the deck transformation $A$ with fundamental domain $\mathcal A$ and axis $\tilde\gamma$.  The integral is {\em folded} by first expressing $\mathcal A=\cup_{B\in\langle A\rangle\backslash\Gamma}B(\mathcal F)$ for $\mathcal F$ a $\Gamma$ fundamental domain and then changing variables in the integrand to obtain the integral of $|\mu|^2\sum_{B\in\langle A\rangle\backslash\Gamma}e^{-2d(B(p),\tilde\gamma)}$ over $\mathcal F$.  As noted following Definition \ref{P} the product of the mass $m(\gamma)$ and the sum is the integral $\mathbb P_{\upsilon}$ for the measured geodesic lamination $\upsilon=(\gamma,m(\gamma))$.  The bound $\langle\mu,\mu\mathbb{P}_{\upsilon}\rangle \le 3\pi\Ddot\ell[\mu,\mu]\le48\langle\mu,\mu\mathbb P_{\upsilon}\rangle $ is now established for simple closed curve laminations.  

Since simple closed curves are dense in $\mathcal{ML}_0$ it suffices to establish the bound for limits in $\mathcal{ML}_0$.  The bound for limits is a property of functions of one variable as follows.  Consider a WP geodesic $\eta$ with arc length parameter $t$ and a sequence of weighted simple closed curves $(\gamma_j,m_j)$ converging to the measured geodesic lamination $\upsilon$.  The functions $m_j\ell_{\gamma_j}$ restricted to $\eta$ converge $C^1$ to $\ell_{\upsilon}$, \cite{Bonmgl, Gardmeas}.  As above at the initial point of $\eta$ the increasing functions 
$f_j(t)=m_j\dot\ell_j(t)$  have first derivatives bounded in terms of $\langle\mu_j,\mu_j\mathbb P_{(\gamma_j,m_j)}\rangle $.  We select a convergent subsequence from the sequence of bounded positive derivatives $\{f_j''(t)\}$ and find that the initial value of $\Ddot\ell_{\upsilon}(t)$ is bounded in terms of the limit of initial values $\lim_jf_j'(t)$.  The desired bound now follows from the continuity of the map $\upsilon\rightarrow\mathbb P_{\upsilon}$ of $\mathcal{ML}_0$ to $C(R)$.   

The final matter is to show that the lower bound for $\mathbb P_{\upsilon}$ on $R_{red}$ provides a lower bound for the pairing $\int_R|\mu|^2\,\mathbb P_{\upsilon}\,dA$.  It suffices to show that the mass of $|\mu|^2$ on $R-R_{red}$ is bounded in terms of the mass on $R_{red}$.  It is enough to consider the mass on a horoball for a single cusp.   For $z$ the standard coordinate for a cusp the hyperbolic metric is $ds^2=(\frac{|dz|}{|z|\log|z|})^2$ with $|z|\le\rho$ the cusp region for $\pi=-\log \rho$.  For the harmonic Beltrami differential $\mu=\overline{\varphi}(ds^2)^{-1}$ the holomorphic quadratic differential $\varphi$ has at most a simple pole at the origin.  The product $\varphi z^2$ is holomorphic on the coordinate chart and vanishes at the origin.  From the Schwarz Lemma the product satisfies  $|\varphi z^2|\le c|z|\max_{|z|=\rho}|\varphi z^2|$ and on multiplying by $\log ^2|z|$ we find the 
bound $|\mu|\le c'|z|\log ^2|z|\max_{|z|=\rho}|\mu|$.  For a suitable value $|z|=\rho$ is contained in $R_{red}$ and the maximum is bounded by Lemma \ref{enhmv}.  The bound on mass in the cusps is established and from Lemma \ref{Pbound1} we have the overall bound $c_1(inj_{red})\ell_{\upsilon}\le \langle\mu,\mu\mathbb P_{\upsilon}\rangle \le c_2(inj_{red})\ell_{\upsilon}$.   The proof is complete.

\subsection{The gradient and Hessian for small geodesic-length}
\label{gradhess}
We present expansions and bounds for the WP pairing of gradients and Hessian of geodesic-length.  The first consideration is Riera's formula for gradients, \cite{Rier}, followed by consideration of the integral pairings bounding the Hessian.  Bounds are developed by comparing the {\em sum} of inverse square exponential-distance to the {\em integral} of inverse square exponential-distance.   

We begin with Riera's length-length formula for the WP pairing of geodesic-length gradients, \cite{Rier}. A less explicit formula was presented in \cite{Wlthur}.   For closed geodesics $\alpha,\beta$ with corresponding deck transformations $A, B$ with corresponding axes $\tilde\alpha,\tilde\beta$ then
\[
\langle\grad\ell_{\alpha},\grad\ell_{\beta}\rangle =\frac{2}{\pi}\bigl(\ell_{\alpha}\delta_{\alpha\beta}\ +\ 
\sum_{\langle A\rangle\backslash\Gamma/\langle B\rangle}(u\log\frac{u+1}{u-1}-2\bigr)\bigr)
\]
for the Kronecker delta $\delta_*$, where for $C\in\langle A\rangle\backslash\Gamma/\langle B\rangle$ then $u=u(\tilde\alpha,C(\tilde\beta))$ is the cosine of the intersection angle if $\tilde\alpha$ and $C(\tilde\beta)$ intersect and is otherwise $\cosh d(\tilde\alpha,C(\tilde\beta))$; for $\alpha=\beta$ the double-coset of the identity element is omitted from the sum, \cite[Theorem 2]{Rier}.  We are interested in the expansion for $\ell_{\alpha},\ell_{\beta}$ small in which case  axes are disjoint and at large distance.   For $u>1$ the positive summand has the expansion for $u$ large
\[
u\log\frac{u+1}{u-1}\,-\,2=\frac23 u^{-2}(1+O(u^{-1}))
\]
and thus for $\ell_{\alpha},\ell_{\beta}$ bounded the double-coset sum is  bounded in terms of $\sum_{\langle A\rangle\backslash\Gamma/\langle B\rangle}e^{-2d}$ the sum of inverse square exponential-distances between axes.  The above formula enables an improvement of the prior bound for pairing gradients, 
\cite[II, Lemmas 2.3 and 2.4]{Wlspeclim}.      
\begin{lemma}
\label{gradest}
The WP pairing of geodesic-length gradients of disjoint geodesics $\alpha,\beta$ satisfies 
\[
0<\langle\grad\ell_{\alpha},\grad\ell_{\beta}\rangle -\frac{2}{\pi}\ell_{\alpha}\delta_{\alpha\beta}
\quad\mbox{is}\quad O(\ell_{\alpha}^2\ell_{\beta}^2)
\]
where for $c_0$ positive the remainder term constant is uniform for $\ell_{\alpha},\ell_{\beta}\le c_0$.  The shortest nontrivial segment on $\mathbb H/\Gamma$  connecting $\alpha$ and $\beta$ contributes a term of order $\ell_{\alpha}^2\ell_{\beta}^2$ to the pairing.
\end{lemma}
\prf  The pairing formula summand is positive and the first inequality is immediate.  The approach follows the considerations for Lemma \ref{Pbound1}.  The sum $\sum e^{-2d}$ is compared to the integral $\int e^{-2d}\,dA$ by introducing comparison balls about projected nearest points.  The projected points will be determined in $\mathbb H$ approximately on the $\beta$ collar boundaries.  To begin, the axis $\tilde\alpha$ has a unique (projected) nearest point on each axis $C(\tilde\beta)$ for $C\in\Gamma/\langle B\rangle$; for $\alpha=\beta$ the 
coset of the identity element is not included.   Double-coset representatives are chosen so that the nearest points for $C(\tilde\beta)$ for $C\in\langle A\rangle\backslash\Gamma/\langle B\rangle$ lie in a fundamental domain $\mathcal A$ for the action of $\langle A\rangle$ on $\mathbb H$; for $\alpha=\beta$ the double-coset of the identity is not included.  From the description of collars for $\ell_{\alpha},\ell_{\beta}\le c_0$ we can set a positive constant $c_1$ such that points in the collars at distance 
$c_1$ from a boundary have injectivity radius at least $c_1$.  We arrange that $c_1$ is at most the $\alpha,\beta$ collar half-width.  For an axes pair we now introduce a secondary projected point along the minimal connecting geodesic segment.  In particular for $\tilde\alpha,C(\tilde\beta)$ the secondary point is inside the collar about $C(\tilde\beta)$ on the $\tilde\alpha$ to $C(\tilde\beta)$ minimal connecting segment at distance $c_1$ from the $C(\tilde\beta)$  boundary.  We introduce the $c_1$-radius ball $\mathcal B(C(\tilde\beta))$ about the secondary point.  The collars about the distinct axes $C(\tilde\beta)$ are disjoint and consequently the balls $\mathcal B$ for distinct axes are disjoint.  Now from the mean value property the double-coset sum of inverse square exponential-distance of $\tilde\alpha$ to the secondary points is bounded by the integral of $e^{-2d(\tilde\alpha,\ )}\,dA$ over the union of the balls.  We next describe a region containing the union $\mathcal U=\cup_{\langle A\rangle\backslash\Gamma/\langle B\rangle}\mathcal B(C(\tilde\beta))$ of balls. 

We specialize the situation for $A$ the transformation $z\rightarrow e^{\ell_{\alpha}}z$ with axis $\tilde\alpha$ the imaginary axis,  fundamental domain $\mathcal A=\{1\le |z|<e^{\ell_{\alpha}}, \Im z>0\}$ for $\langle A\rangle$  and collar $\mathcal C(\tilde\alpha)=\{1\le |z|< e^{\ell_{\alpha}},\,\sinh d(z,\tilde\alpha) \sinh \ell_{\alpha}/2 \le 1 \}$.  The secondary points are determined in $\mathcal A$; since $\alpha$ and $\beta$ collars are disjoint the secondary points are determined in $\mathcal A-\mathcal C(\tilde\alpha)$ at least distance $c_1$ from $\mathcal C(\tilde\alpha)$.  By specification at a secondary point the region $\mathcal R=\{e^{-\ell_{\alpha}}\le |z| \le e^{2\ell_{\alpha}},\, \sinh d(z,\tilde\alpha) \sinh \ell_{\alpha}/2 \ge 1 \}$ contains the $c_1$-ball about the  point.  The region 
$\mathcal R$ contains the union $\mathcal U$.  And as already noted the inverse square exponential-distance to 
$\tilde\alpha$  is bounded by $\sin ^2\theta$.   The comparison integral is $\int_{\mathcal R}\sin ^2\theta\, dA$ is $O(\ell_{\alpha}^2)$.  As the last step we observe that a minimal connecting segment from $\tilde\alpha$ to $C(\tilde\beta)$ consists of a segment from $\tilde\alpha$ to the secondary point and a segment from the secondary point to $C(\tilde\beta)$.  The second segment has length $w_{\beta}-c_1$, for $w_{\beta}$  the half-width of the $\beta$ collar.  In particular $e^{-2w_{\beta}+2c_1}$ is bounded in terms of $\ell_{\beta}^2$.  The sum for distances between axes is now bounded in terms of $\ell_{\alpha}^2\ell_{\beta}^2$, the desired bound.  

The shortest segment on the surface connecting $\alpha$ and $\beta$ consists of a half-width segment crossing the 
$\alpha$ collar, a segment in the complement of the $\alpha$ and $\beta$ collars, and a half-width segment crossing 
the $\beta$ collar.  For the corresponding component lengths $w_{\alpha}, w_0$ and $w_{\beta}$ the segment contributes 
$\ell_{\alpha}^2\ell_{\beta}^2e^{-2w_0}$ to the inverse square exponential-distance sum.  The proof is complete.

We begin considerations for the Hessian with the closed geodesic $\alpha$ with corresponding deck transformation $A$ given as $z\rightarrow e^{\ell_{\alpha}}z$ for $z\in\mathbb H$, with fundamental domain $\mathcal A=\{1\le |z|<e^{\ell_{\alpha}}, \Im z>0\}$ for $\langle A\rangle$  and with a collar $\mathcal C(\tilde\alpha)=\{1\le |z|< e^{\ell_{\alpha}},\,\ell_{\alpha}\le \arg z \le \pi -\ell_{\alpha}\}$ about the axis   of $A$.  Consider further a holomorphic quadratic differential $\varphi\in Q(\Gamma)$ with expansion (\ref{zA}) and zeroth coefficient $a_0(\varphi)$.  Introduce the alteration
\[
\varphi^\flat=\varphi-a_0(\varphi)(\frac{dz}{z})^2.
\]
The linear functional $\varphi\rightarrow\Re a_0(\varphi)$ by Theorem \ref{quadform} is represented by $\eta_{\alpha}=-(4\ell_{\alpha})^{-1}\grad\ell_{\alpha}$ considered as an element of the dual $Q(\Gamma)^*$.  By Lemma \ref{gradest} the WP norm of $\eta_{\alpha}$ has magnitude $\ell_{\alpha}^{-1/2}$ and thus by the Cauchy-Schwarz inequality $|\Re a_0(\varphi)|=|\langle\varphi,\eta_{\alpha}\rangle |$ is bounded in terms of $\ell_{\alpha}^{-1/2}\|\varphi\|_{WP}$ for $\ell_{\alpha}$ bounded.  

We introduce supplementary Hermitian pairings for $\ell_{\alpha}$ sufficiently small, for the subcollar $\mathcal{SC}(\tilde\alpha)=\{1\le |z|< e^{\ell_{\alpha}},\,2\ell_{\alpha}\le \arg z \le \pi -2\ell_{\alpha}\}$ and $z$ in $\mathbb H$ given in polar form $re^{i\theta}$
\[
\mathcal{QS}_{collar}(\psi)=\int_{\mathcal{SC}(\tilde\alpha)}|\psi|^2\sin^2\theta\, (ds^2)^{-1}
\]
and 
\[
\mathcal{QS}_{comp}(\psi)=\int_{\mathcal A-\mathcal{SC}(\tilde\alpha)}|\psi|^2\sin^2\theta\, (ds^2)^{-1}.
\]
In Theorem \ref{qshess} below we bound $\hess \ell_{\alpha}$ in terms of the pairings $\mathcal{QS}_{collar}(\varphi^\flat)$ and $\mathcal{QS}_{comp}(\varphi^\flat)$.   As the first step in Lemma \ref{qscomp} we bound $\mathcal{QS}_{comp}(\varphi^\flat)$ by first replacing $\varphi^\flat$ with $\varphi$ and then considering the $L^1$ and mean value estimate for a series similar to $\mathbb P_{\alpha}$.  As the second step in Lemma \ref{qscollar} we bound $\mathcal{QS}_{collar}(\varphi^\flat)$ by introducing the quotient $\mathbb H/\langle A\rangle$,  applying the Cauchy Integral Formula to represent $\varphi^\flat$ and  applying the Schwarz Lemma (using that the zeroth coefficient vanishes.)  We find that $\varphi^\flat$ is small in the core of the collar.

\begin{lemma}
\label{qscomp}
The Hermitian pairing  $\mathcal{QS}_{comp}(\varphi^\flat)$ is uniformly $O(\ell^2_{\alpha}\|\varphi\|^2_{WP})$  for  $\ell_{\alpha}$ sufficiently small.
\end{lemma}
\prf We first consider replacing $\varphi^\flat$ by $\varphi$.  The pairing
\[
\mathcal{QS}_{comp}((\frac{dz}{z})^2)=2\int_0^{2\ell_{\alpha}}\int_1^{e^{\ell_{\alpha}}}\sin^4\theta\,d\theta d\log r
\]
has magnitude $\ell_{\alpha}^6$ and as noted above $|a_0(\varphi)|^2$ has magnitude $\ell^{-1}_{\alpha}\|\varphi\|^2_{WP}$.  It follows that $\mathcal{QS}_{comp}(\varphi)=\mathcal{QS}_{comp}(\varphi^\flat)\,+\,O(\ell_{\alpha}^5\|\varphi\|^2_{WP})$ and that it is sufficient to bound $\mathcal{QS}_{comp}(\varphi)$.

We next introduce the characteristic function $\chi_a(\theta)$ of the subinterval $[a,\pi-a]$ in the interval $(0,\pi)$ and consider the series 
\[
S(z)=\sum_{B\in\langle A\rangle\backslash\Gamma}s(B(z))\quad\mbox{for}\quad s(w)=(1-\chi_a(\arg w))\sin^2\arg w.
\]
For $f$ a bounded $\Gamma$-invariant function on $\mathbb H$ and $\mathcal F$ a fundamental domain for $\Gamma$ the pairing with $S$ can be unfolded as follows
\[
\int_{\mathcal F}Sf\,dA=\int_{\mathcal F}\sum_{B\in\langle A\rangle\backslash\Gamma}S(B(z))f\,dA=\sum_{B\in\langle A\rangle\backslash\Gamma}\int_{B(\mathcal F)}S(B(z))f\,dA=\int_{\mathcal A}sf\,dA.
\]
For $f$ identically unity we have that the $L^1$-norm of $S$ is the integral 
$\int_{\mathcal A}s\,dA=2a\ell_{\alpha}$.  

The next matter is to bound the location of points in $\mathcal A-\mathcal{SC}(\tilde\alpha)$ with a given distance to $\widetilde{thick}$, the lift of the $thick$ region of $\mathbb H/\Gamma$.  The geodesic segment in $\mathcal A$ from $z$ to $\tilde\alpha$ necessarily enters $\widetilde{thick}$.  It follows that the distance $\delta(z)$ of a point $z$ to $\widetilde{thick}$ is at most $d(z,\tilde\alpha)-w_{\alpha}$ for $w_{\alpha}$ the half-width of the $\alpha$ collar.  From the formula for the distance to $\tilde\alpha$ and the description of the collar it follows that $e^{\delta(z)}\le  \ell_{\alpha}\csc\arg z$ or equivalently that $\sin \arg z \le  \ell_{\alpha}e^{-\delta(z)}$.  It further follows that points of $\mathcal A-\mathcal{SC}(\tilde\alpha)$ at a given distance $\delta$ to $\widetilde{thick}$ are contained in $\mathcal R_{\delta}=\{z\in\mathcal A\mid \arg z\le 2\ell_{\alpha}e^{-\delta }\mbox{ or }\pi-2\ell_{\alpha}e^{-\delta }\le \arg z\}$.

We now proceed for $\ell_{\alpha}$ sufficiently small to bound the value of $S$ for a point $z$ at distance $\delta$ to $\widetilde{thick}$.  For $\rho=\min\{1,inj(z)\}$ the larger region $\mathcal R_{\delta-1}$ contains disjoint balls of radius $\rho$ about the $\langle A\rangle\backslash\Gamma$-orbit of $z$ in $\mathcal A-\mathcal{SC}(\tilde\alpha)$.  It  follows that $\int_{\mathbf B(z,\rho)}S\,dA$ is bounded by $4e\ell_{\alpha}^2e^{-\delta}$ by choosing $a=2e\ell_{\alpha}e^{-\delta}$ for the above unfolding.  The function $s$ satisfies a mean value estimate and thus from Lemma \ref{enhmv} the value of $S$ at $z$ is now bounded in terms of $\ell_{\alpha}^2e^{-\delta}inj(z)^{-1}$.  Finally from Lemma \ref{enhcollar} the $L^{\infty}$-norm of $S$ is bounded in terms of $\ell_{\alpha}^2$, the desired overall bound for $S$.  And for the selection above of $f=|\varphi(ds^2)^{-1}|^2$ and $a=2\ell_{\alpha}$ we have that 
\[
\mathcal{QS}_{comp}(\varphi)=\int_{\mathcal F}S|\varphi(ds^2)^{-1}|^2\,dA\le\|S\|_{\infty}\|\varphi\|^2_{WP}.
\]
The proof is complete.

\begin{lemma}
\label{qscollar}
The Hermitian pairing  $\mathcal{QS}_{collar}(\varphi^\flat)$ is uniformly $O(\ell^2_{\alpha}\|\varphi\|^2_{WP})$ for  $\ell_{\alpha}$ sufficiently small.
\end{lemma}
\prf We first bound $\varphi$ near the $\alpha$ collar boundaries and then consider $\varphi$ in the annulus model for $\mathbb H/\langle A\rangle$.  Consider the collar $\mathcal C(\tilde\alpha)$ for the axis of $A$.   For $\ell_{\alpha}$ sufficiently small the injectivity radius of $\mathbb H/\Gamma$ is bounded away from zero on the rays $\arg z=2\ell_{\alpha},\pi-2\ell_{\alpha}$.  It follows from the mean value estimate Lemma \ref{enhmv} that on the two rays the quantity $|\varphi|(ds^2)^{-1}$ is uniformly bounded in terms of $\|\varphi\|_{WP}$; for $\ell_{\alpha}$ sufficiently small $|\varphi^\flat|(ds^2)^{-1}$ is similarly bounded.  The quotient annulus $\mathbb H/\langle A\rangle$ is described by setting
\[
w=e^{\frac{2\pi i}{\ell_{\alpha}}\log z} \quad \mbox{for } z \mbox{ in } \mathbb H,\ \tau<|w|<1 \mbox{ and } \tau=e^{-\frac{2\pi^2}{\ell_{\alpha}}}
\]
and has hyperbolic metric
\[
ds^2=\bigl(\frac{\ell_{\alpha}}{2\pi}\csc \frac{\ell_{\alpha}}{2\pi}\log |w| \bigr)^2\bigl|\frac{dw}{w}\bigr|^2.
\]
We set
\[
\varphi^\flat=f(w)\bigl(\frac{dw}{w}\bigr)^2
\]
and note that the rays $\arg z=2\ell_{\alpha},\pi-2\ell_{\alpha}$ correspond to the circles 
$|w|=e^{-4\pi},\,e^{4\pi}\tau$ where on $\bigl|\bigl(\frac{dw}{w}\bigr)^2\bigr|(ds^2)^{-1}$ is approximately $(4\pi)^2$ for $\ell_{\alpha}$ sufficiently small.  It follows that on the circles $|w|=e^{-4\pi},\,e^{4\pi}\tau$ the function $f(w)$ is bounded in terms of $\|\varphi\|_{WP}$, the desired bound.

We next apply the Cauchy Integral Formula and write
\[
f(w)=F(w)\ +\ G(w)
\]
for
\[
F(w)=\frac{1}{2\pi i}\int_{|u|=e^{-4\pi}}\frac{f(u)}{u-w}\,du \quad \mbox{and}
\quad G(w)=\frac{-1}{2\pi i}\int_{|u|=e^{4\pi}\tau}\frac{f(u)}{u-w}\,du. 
\]
The function $F(w)$ is holomorphic in $|w|<e^{-4\pi}$ and vanishes at the origin since $F(0)=\frac{1}{2\pi i}\int\frac{f(u)}{u}\,du$ is the zeroth coefficient of $\varphi^\flat$.  From the above paragraph on $|w|\le e^{-8\pi}$ the function $F$ is bounded in terms of $\|\varphi\|_{WP}$. It follows from the Schwarz Lemma that $|F(w)|\le c|w|\|\varphi\|_{WP}$ on $|w|\le e^{-8\pi}$, the desired bound.  The integral for $G(w)$ is first transformed by the change of variable $uv=\tau$ 
\[
G(w)=\frac{-1}{2\pi i}\int_{|v|=e^{-4\pi}}\frac{f(\frac{\tau}{v})}{\frac{\tau}{v}-w}\tau\,\frac{dv}{v^2}
=\frac{1}{2\pi i}\int_{|v|=e^{-4\pi}}\frac{\frac{\tau}{w}f(\frac{\tau}{v})}{v-\frac{\tau}{w}}\,\frac{dv}{v}.
\]
The function $G(w)$ is holomorphic in the quantity $\frac{\tau}{w}$ for $|\frac{\tau}{w}|\le e^{-8\pi}$, vanishes at the origin and from the Schwarz Lemma satisfies $|G(w)|\le c|\frac{\tau}{w}|\|\varphi\|_{WP}$ on $|\frac{\tau}{w}|\le e^{-8\pi}$, the desired bound.  

We apply the bounds for $F$ and $G$ to bound the pairing $\mathcal{QS}_{collar}$.  The estimate is demonstrated by considering
\begin{multline*}
\mathcal{QS}_{collar}(F\bigl(\frac{dw}{w}\bigr)^2)=\int_{\mathcal{SC}(\tilde\alpha)}(\sin^2\frac{\ell_{\alpha}}{2\pi}\log |w|)\bigl|F\bigl(\frac{dw}{w}\bigr)^2\bigr|^2\,(ds^2)^{-1}\\
\le c\|\varphi\|^2_{WP}\int_{\mathcal{SC}(\tilde\alpha)}(\frac{2\pi}{\ell_{\alpha}})^2 (\sin^4\frac{\ell_{\alpha}}{2\pi}\log |w|)\,|dw|^2\,
\end{multline*}
and substituting $|\sin \mu|\le |\mu|$
\[
\le c\|\varphi\|^2_{WP}\int_{|w|\le e^{-8\pi}}(\frac{\ell_{\alpha}}{2\pi})^2 \log ^4 |w|\,|dw|^2\le c\ell^2_{\alpha}\|\varphi\|^2_{WP},
\]
the desired bound.  The estimate for $G$ is obtained in the same manner.  The proof is complete.

We are now ready to present the overall bound for the Hessian.  
\begin{theorem}
\label{qshess}
The variation of geodesic-length $\ell_{\alpha}$ satisfies
\[
2\ell\Ddot\ell[\mu,\mu]-\dot\ell^2[\mu]-3\dot\ell^2[i\mu]\quad\mbox{is}\quad O(\ell_{\alpha}^3\|\mu\|^2_{WP})
\]
and
\[
\ell\partial\overline\partial\ell[\mu,\mu]-2\partial\ell[\mu]\overline\partial\ell[\mu]\quad\mbox{is}\quad O(\ell_{\alpha}^3\|\mu\|^2_{WP})
\]
for a harmonic Beltrami differential $\mu\in\mathcal H(R)$  where for $c_0$ positive the remainder term constant is uniform for $\ell_{\alpha}\le c_0$.
\end{theorem}
\prf The specified variation of geodesic-length is  bounded by a constant multiple of $\ell_{\alpha}\mathcal{QS}(\varphi^\flat,\varphi^\flat)$ for $\mu=(z-\bar z)^2\overline\varphi$,  $\varphi\in Q(R)$, from Corollary \ref{first} and Lemma \ref{QE}.   By Lemmas \ref{qscomp} and \ref{qscollar} the pairing $\mathcal{QS}$ is bounded by the pairing sum $\mathcal{QS}_{collar}+\mathcal{QS}_{comp}$ with the sum uniformly $O(\ell_{\alpha}^2 \|\varphi\|^2_{WP})$ for $\ell_{\alpha}$ sufficiently small.  It only remains to bound the pairing $\mathcal{QS}$ for $\ell_{\alpha}\ge c_1$, $c_1$ a universal constant.  Given $c_0>c_1$ it suffices from the considerations of Theorem \ref{hessmgl} to bound the series $\mathbb P_{\alpha}$ for $\ell_{\alpha}\le c_0$.  An appropriate bound is provided in the following Lemma.  The proof is complete.

A general bound for the series $\mathbb P_{\alpha}$  provides a general bound for the gradient $\grad \ell_{\alpha}$, as well as for the Hessian $\hess \ell_{\alpha}$.  

\begin{lemma}
\label{Pbound2}
The series $\mathbb P_{\alpha}$ and gradient of geodesic-length for a simple geodesic $\alpha$ are bounded as follows
\[
\mathbb P_{\alpha}\le c(1+\ell_{\alpha}e^{\ell_{\alpha}/2})\quad\mbox{and}
\quad \langle\grad\ell_{\alpha},\grad\ell_{\alpha}\rangle  \le c(\ell_{\alpha}\,+\,\ell_{\alpha}^2 e^{\ell_{\alpha}/2})
\]
for a universal positive constant $c$.
\end{lemma}
\prf  The approach follows the considerations for the series $S$ in the proof of Lemma \ref{qscomp}.  We consider the geodesic $\alpha$ corresponding to the deck transformation $A$ with axis $\tilde\alpha$ the imaginary axis, and  fundamental domain $\mathcal A$.  As noted in the proof of Lemma \ref{QE}, 
$e^{-2d(z,\tilde\alpha)}\le \sin^2\theta \le 4e^{-2d(z,\tilde\alpha)}$ for $z$ in $\mathbb H$ in polar form $z=re^{i\theta}$.  The function $e^{-2d(z,\tilde\alpha)}$ has $L^1$-norm bounded in terms of $\int_{\mathcal A}\sin^2\theta\, dA=\pi\ell_{\alpha}$ and satisfies a mean value estimate.  Furthermore the injectivity radius is bounded below on $thick$ by a positive constant and by Lemma \ref{separ} is bounded below in terms of $\min\{\ell_{\alpha},e^{-\ell_{\alpha}/2}\}$ on collars that $\alpha$ intersects.  It follows from Lemma \ref{enhmv} that on $thick$ and  collars that $\alpha$ intersects the series $\mathbb P_{\alpha}$ is bounded in terms of $1+\ell_{\alpha}e^{\ell_{\alpha}/2}$, a desired bound. 

It remains to bound the series on the components of $thin$ disjoint from $\alpha$.   We proceed and consider in the universal cover the location of points $z$ of $\mathcal A$ with distance $\delta(z)$ to $\widetilde{thick}$, the lift of the $thick$ region of $\mathbb H/\Gamma$.  Since $\alpha$ is disjoint from the remaining $thin$  components, the geodesic segment from $z$ to $\tilde\alpha$ enters $\widetilde{thick}$ and consequently 
$\delta(z)\le d(z,\tilde\alpha)$ and $\sin\theta \le 2 e^{-\delta(z)}$.   
We have that points of $\mathcal A$ at a given positive distance $\delta$ to $\widetilde{thick}$ are  contained in the region $\mathcal R'_{\delta}=\{z\in\mathcal A \mid \arg z\le \rho' \mbox{ or }\pi-\rho'\le \arg z \}$ for 
$\rho'=\min\{2e^{-\delta}, \frac{\pi}{2}\}$.
The larger region $\mathcal R'_{\delta -1}$ contains disjoint balls of radius $\min\{1,inj(z)\}$ about the $\langle A\rangle\backslash\Gamma$-orbit of $z$ in $\mathcal A$.   
It now follows from the formula for $\int_{\mathcal A}s\,dA$ of Lemma \ref{qscomp} that $\int_{\mathbf B(z,\rho)}\mathbb P_{\alpha}\,dA$ is bounded in terms of $\ell_{\alpha}e^{-\delta}$.  It then follows from Lemma \ref{enhmv} that $\mathbb P_{\alpha}(z)$ is bounded in terms of $\ell_{\alpha}e^{-\delta}inj(z)^{-1}$.   The desired  bound now follows from Lemma \ref{enhcollar}.  In combination with the bound of the above paragraph we have that $\mathbb P_{\alpha}$ is  bounded by $1+\ell_{\alpha}e^{\ell_{\alpha}/2}$, the desired bound.

We next bound the gradient of geodesic-length.  From the inequality 
$\bigl|(\frac{dz}{z})^2\bigr|(ds^2)^{-1}=\sin^2\theta \le4 e^{-2d(z,\tilde\alpha)}$ and remarks following Definition \ref{P} we have that 
$\bigl|\Theta_{\alpha}\bigr|(ds^2)^{-1} \le 4\,\mathbb P_{\alpha}$.
From formula (\ref{dgl}) the gradient $\grad \ell_{\alpha}$ is represented by $\frac{2}{\pi}\Theta_{\alpha}$ and the above $L^1$ and $L^{\infty}$ bounds for $\mathbb P_{\alpha}$ provide the desired bound for the inner product.    The proof is complete.   

\section{WP geometry}

\subsection{Comparisons for the WP metric}
\label{WPcompar}
Let $\mathcal T$ be the Teichm\"{u}ller space of genus $g$, $n$ punctured Riemann surfaces with $F$ the reference topological surface for the markings, \cite{ImTan, Ngbook}.  The points of $\mathcal T$ are equivalence classes 
$\{(R,f)\}$ of surfaces with reference homeomorphisms $f:F\rightarrow R$.  
The {\em complex of curves} $C(F)$ is defined as follows.  The vertices of
$C(F)$ are the free homotopy classes of homotopically nontrivial, non peripheral,
simple closed curves on $F$.  An edge of the complex consists of a pair of
homotopy classes of disjoint simple closed curves.  A $k$-simplex consists of
$k+1$ homotopy classes of mutually disjoint simple closed curves.  A maximal
set of mutually disjoint simple closed curves, a {\em partition},
has $3g-3+n$ elements.  The mapping class group $Mod$ acts on the complex $C(F)$.

The Fenchel-Nielsen coordinates for $\mathcal T$ are given in terms of geodesic-lengths and
lengths of auxiliary geodesic segments, \cite{Abbook,Busbook,ImTan,MasFN, WlFN}.  A partition
 $\mathcal{P}=\{\alpha_1,\dots,\alpha_{3g-3+n}\}$ decomposes the
reference surface $F$ into $2g-2+n$ components, each homeomorphic to
a sphere with a combination of three discs or points removed.  A marked
Riemann surface $(R,f)$ is likewise decomposed into pants by the geodesics
representing the elements of $\mathcal P$.  Each component pants, relative to its hyperbolic
metric, has a combination of three geodesic boundaries and cusps.  For each
component pants the shortest geodesic segments connecting boundaries determine
 designated points  on each boundary.  For each geodesic $\alpha$ in the pants
decomposition of $R$ a parameter $\tau_{\alpha}$ is defined as the displacement along
the geodesic between designated points, one  for each side of the geodesic.
For marked Riemann surfaces close to an initial reference marked Riemann surface, the
displacement $\tau_{\alpha}$ is the distance between the designated points; in
general the displacement is the analytic continuation (the lifting) of the
distance measurement.  For $\alpha$ in $\mathcal P$ define the {\em
Fenchel-Nielsen  angle} by $\vartheta_{\alpha}=2\pi\tau_{\alpha}/\ell_{\alpha}$.
The Fenchel-Nielsen coordinates for Teichm\"{u}ller space for the
decomposition $\mathcal P$ are
$(\ell_{\alpha_1},\vartheta_{\alpha_1},\dots,\ell_{\alpha_{3g-3+n}},\vartheta_{\alpha_{3g-3+n}})$.
The coordinates provide a real analytic equivalence of $\mathcal T$ to
$(\mathbb{R}_+\times \mathbb{R})^{3g-3+n}$, \cite{Abbook,Busbook,ImTan,MasFN,WlFN}.

A bordification (a partial compactification) of Teichm\"{u}ller space is introduced by extending the range
of the Fenchel-Nielsen parameters.  The added points correspond to unions of
hyperbolic surfaces with formal pairings of cusps.  The interpretation of {\em length
vanishing} is the key ingredient.  For an $\ell_{\alpha}$ equal to zero, the
angle $\vartheta_{\alpha}$ is not defined and in place of the geodesic for
$\alpha$ there appears a pair of cusps; the reference map $f$ is now a homeomorphism of
$F-\alpha$ to a union  of hyperbolic surfaces (curves parallel to $\alpha$
map to loops encircling the cusps).  The parameter space for a pair
$(\ell_{\alpha},\vartheta_{\alpha})$ will be the identification space $\mathbb{R}_{\ge
0}\times\mathbb{R}/\{(0,y)\sim(0,y')\}$.  More generally for the pants decomposition
$\mathcal P$ a frontier set $\mathcal{F}_{\mathcal P}$ is added to the
Teichm\"{u}ller space by extending the Fenchel-Nielsen parameter ranges: for
each $\alpha\in\mathcal{P}$, extend the range of $\ell_{\alpha}$ to include
the value $0$, with $\vartheta_{\alpha}$ not defined for $\ell_{\alpha}=0$.  The
points of $\mathcal{F}_{\mathcal P}$ parameterize unions of Riemann
surfaces with each $\ell_{\alpha}=0,\alpha\in\mathcal{P},$ specifying a pair
of cusps. The points of $\mathcal F_{\mathcal P}$ are Riemann surfaces with nodes in the sense of Bers, \cite{Bersdeg}.   For a simplex $\sigma\subset\mathcal{P}$, define the
$\sigma$-null stratum, a subset of $ \mathcal F_{\mathcal P}$, as $\mathcal{T}(\sigma)=\{R\mid \ell_{\alpha}(R)=0
\mbox{ iff }\alpha\in\sigma\}$. Null strata are given as products of lower dimensional
Teichm\"{u}ller spaces.   The frontier set $\mathcal{F}_{\mathcal P}$
is the union of the $\sigma$-null strata for the sub simplices of
$\mathcal{P}$.  Neighborhood bases for points of $\mathcal{F}_{\mathcal P}
\subset\mathcal{T}\cup\mathcal{F}_{\mathcal P}$  are specified by the
condition that for each simplex $\sigma\subset\mathcal P$  the projection
$((\ell_{\beta},\vartheta_{\beta}),\ell_{\alpha}):
\mathcal{T}\cup\mathcal{T}(\sigma)\rightarrow\prod_{\beta\notin\sigma}(\mathbb{R}_+\times\mathbb{R})\times\prod_{\alpha\in\sigma}(\mathbb{R}_{\ge
0})$  is continuous.    For a simplex $\sigma$ contained in partitions $\mathcal P$ and $\mathcal P'$ the specified neighborhood systems for $\mathcal T \cup\mathcal{T}(\sigma)$ are equivalent.  The {\em
augmented Teichm\"{u}ller space} $\Tbar=\mathcal{T}\cup_{\sigma\in
C(F)}\mathcal{T}(\sigma)$ is the resulting stratified topological space,
\cite{Abdegn, Bersdeg}.  $\Tbar$ is not locally compact since points of the
frontier do not have relatively compact neighborhoods; the  neighborhood bases are
unrestricted in the $\vartheta_{\alpha}$ parameters for $\alpha$ a $\sigma$-null.
The action of $Mod$ on $\mathcal T$ extends to an action by homeomorphisms on
$\Tbar$ (the action on $\Tbar$ is not properly discontinuous) and the quotient
$\Tbar/Mod$ is topologically the compactified moduli space of stable curves, \cite[see Math. Rev.56 \#679]{Abdegn}.

We present an alternate description of the frontier points in terms of representations of groups and the Chabauty topology.  A Riemann surface with punctures and hyperbolic metric is uniformized by a cofinite subgroup $\Gamma\subset PSL(2;\mathbb R)$.  A puncture corresponds to the $\Gamma$-conjugacy class of a maximal parabolic subgroup.  In general a Riemann surface with punctures corresponds to the $PSL(2;\mathbb R)$ conjugacy class of a tuple $(\Gamma,\langle\Gamma_{01}\rangle ,\dots,\langle\Gamma_{0n}\rangle )$ where $\langle\Gamma_{0j}\rangle $ are the maximal parabolic classes and a labeling for punctures is a labeling for conjugacy classes.  A {\em Riemann surface with nodes} $R'$ is a finite collection of $PSL(2;\mathbb R)$ conjugacy classes of tuples 
$(\Gamma^\ast,\langle\Gamma_{01}^\ast\rangle ,\dots,\langle\Gamma_{0n^\ast}^\ast\rangle )$ with a formal pairing of certain maximal parabolic classes.  The conjugacy class of a tuple is called a {\em part} of $R'$.  The unpaired maximal parabolic classes are the punctures of $R'$ and the genus of $R'$ is defined by the relation $Total\ area=2\pi(2g-2+n)$.  A cofinite $PSL(2;\mathbb R)$ injective representation of the fundamental group of a surface is topologically allowable provided peripheral elements correspond to peripheral elements.  A point of the Teichm\"{u}ller space $\mathcal T$ is given by the $PSL(2;\mathbb R)$ conjugacy class of a topologically allowable injective cofinite representation of the fundamental group $\pi_1(F)\rightarrow\Gamma\subset PSL(2;\mathbb R)$.  For a simplex $\sigma$ a point of $\mathcal T(\sigma)$ is given by a collection $\{(\Gamma^\ast,\langle\Gamma_{01}^\ast\rangle ,\dots,\langle\Gamma_{0n^\ast}^\ast\rangle )\}$ of tuples with: a bijection between $\sigma$ and the paired maximal parabolic classes; a bijection between the components $\{F_j\}$ of $F-\sigma$ and the conjugacy classes of parts $(\Gamma^j,\langle\Gamma_{01}^j\rangle ,\dots,\langle\Gamma_{0n^j}^j\rangle )$ and the $PSL(2;\mathbb R)$ conjugacy classes of topologically allowable isomorphisms $\pi_1(F_j)\rightarrow\Gamma^j$, \cite{Abdegn, Bersdeg}. We are interested in geodesic-lengths for a sequence of points of $\mathcal T$ converging to a point of $\mathcal T(\sigma)$.  The convergence of hyperbolic metrics provides that for closed curves of $F$ disjoint from $\sigma$ geodesic-lengths converge, while closed curves with essential $\sigma$ intersections have geodesic-lengths tending to infinity, \cite{Bersdeg, Wlhyp}.  

We refer to the Chabauty topology to describe the convergence for the $PSL(2;\mathbb R)$ representations.  
Chabauty introduced a topology for the space of discrete subgroups of a locally compact group, \cite{Chb}.   A neighborhood of $\Gamma\subset PSL(2;\mathbb R)$ is specified by a neighborhood $U$ of the identity in $PSL(2;\mathbb R)$ and a compact subset $K\subset PSL(2;\mathbb R)$.  A discrete group $\Gamma'$ is in the neighborhood $\mathcal N(\Gamma,U,K)$ provided $\Gamma'\cap K\subseteq\Gamma U$ and $\Gamma\cap K\subseteq\Gamma'U$.  The sets $\mathcal N(\Gamma,U,K)$ provide a neighborhood basis for the topology. We consider a sequence of points of $\mathcal T$ converging to a point of $\mathcal T(\sigma)$ corresponding to 
$\{(\Gamma^\ast,\langle\Gamma_{01}^\ast\rangle ,\dots,\langle\Gamma_{0n^\ast}^\ast\rangle )\}$.  Important for the present considerations is the following property.  Given a sequence of points of $\mathcal T$ converging to a point of $\mathcal T(\sigma)$ and a component $F_j$ of $F-\sigma$ there exist $PSL(2;\mathbb R)$ conjugations such that restricted to $\pi_1(F_j)$ the corresponding representations converge elementwise to $\pi_1(F_j)\rightarrow\Gamma^j$, \cite[Thrm. 2]{HrCh}.

We continue for a simplex $\sigma$ to consider geodesic-length functions in a neighborhood of an augmentation point $p$ of 
$\mathcal T(\sigma)\subset\overline{\mathcal T}$.  For the following considerations 
we refer to the elements of $\sigma$ as the short geodesics. 
We are interested in collections of geodesic-length functions giving coordinates in a neighborhood of $p$.
\begin{definition}
\label{relbas}
A relative length basis for a point $p$ of $\mathcal T(\sigma)$ is a collection $\tau$ of vertices of $C(F)$   disjoint from the elements of $\sigma$ such that at $p$ the gradients $\{\grad \ell_{\beta}\}_{\beta\in\tau}$ provide the germ of a frame over $\mathbb R$ for the tangent space $T\mathcal T(\sigma)$.
\end{definition}
Examples of relative length bases are given as the union of a partition and a dual partition 
for $R-\sigma$, see \cite[Chap. 3 Secs. 3 \& 4]{Busbook}, \cite[{\em markings} in Sec. 2.5]{MaMiII} and \cite[Thrm. 3.4]{WlFN}.   The lengths of the elements of a relative length basis are necessarily bounded on a neighborhood in $\overline{\mathcal T}$.  We consider below for a geodesic $\gamma$ the root geodesic-length $\ell^{1/2}_{\gamma}$ and gradient $\lambda_{\gamma}=\grad \ell^{1/2}_{\gamma}$.  We also introduce the convention that the pairing $\langle\grad\ell_{\alpha},\grad\ell_{\beta}\rangle $ vanishes for geodesics 
$\alpha,\beta$ on distinct components (parts) of a Riemann surface with nodes.  
\begin{lemma}
\label{paircont}
The WP pairing of geodesic-length gradients $\langle\grad \ell_{\alpha},\grad \ell_{\beta}\rangle $ is continuous in a neighborhood of a point $p$ of $\mathcal T(\sigma)\subset\overline{\mathcal T}$ for $\alpha$ and $\beta$ disjoint from the simplex $\sigma$.  The matrix of WP pairings for a combined short and relative length basis  $\{\lambda_{\alpha},i\lambda_{\alpha},\grad\ell_{\beta}\}_{\alpha\in\sigma,\,\beta\in\tau}$ determines a germ at $p$ of a continuous map from $\Tbar$ into a real general linear group $GL(\mathbb R)$.  The matrix of WP Hermitian pairings for the basis  
$\{\lambda_{\alpha}\}_{\alpha\in\sigma'}$ for $\sigma'$ a partition containing $\sigma$ determines a germ at $p$ of a continuous map into a complex general linear group $GL(\mathbb C)$.
\end{lemma}
\prf  We consider a sequence in $\mathcal T$ converging to the point $p$. Conjugate the representations to arrange for each that $\alpha$ has the imaginary axis $\tilde\alpha$ as a lift. Use the representations to label the elements of the discrete groups.  We suppress subscripts and write $\Gamma$ for a group in the sequence.  For $\langle A\rangle\subset\Gamma$ the stabilizer of $\tilde\alpha$, the group $\langle A\rangle$ acts on the lifts of $\beta$.  Select a lift $\tilde\beta$ of $\beta$ at minimal distance from $\tilde\alpha$ and with nearest point projection to $\tilde\alpha$ in $\mathcal I=\{iy\mid 1\le y<e^{\ell_{\alpha}}\}\subset\mathbb H$.  Either the selected lifts $\tilde\beta$ have distance to $\tilde\alpha$ tending to infinity or for a subsequence the lifts converge. In $\Gamma$ denote the stabilizer of $\tilde\beta$ by  $\langle B\rangle $ and select double-coset representatives with the nearest point for $C(\tilde\beta)$, 
$C\in\langle A\rangle\backslash\Gamma/\langle B\rangle$ also in $\mathcal I$.   

We first consider the fixed number of double-cosets for which $C(\tilde\beta)$ intersects $\tilde\alpha$ 
in $\mathcal I$.  The corresponding transformations $CBC^{-1}$ have translation length $\ell_{\beta}$ and axes intersecting a relatively compact set.  The Chabauty convergence of groups provides that the contributions of the double-cosets converges.  We next consider the double-cosets for the non intersecting axes $\tilde\alpha$ and $C(\tilde\beta)$.  As in the proof of Lemma \ref{gradest} introduce balls about the secondary points on the minimal geodesic segments connecting $C(\tilde\beta)$ to $\tilde\alpha$.  The balls are contained in the collars about the distinct axes $C(\tilde\beta)$ and consequently the balls for distinct axes are disjoint.  Now from the mean value property the $\langle A\rangle\backslash\Gamma/\langle B\rangle$ double-coset sum of $f(u)=u\log\frac{u+1}{u-1}-2$ for $u=d(\tilde\alpha,C(\tilde\beta))$ is bounded by the corresponding integral.  Furthermore if the summands are arranged in order of increasing distance then a uniform majorant for the sum is given in terms of the integrals $\int_{\mathbf B(q,r)}f(u)\,dA$.  Accordingly the limit of the double-coset sum as representations tend to a Riemann surface with nodes is the sum of the term-by-term limits.  We consider the term limits. 

Let $\Gamma'$ be a part of $p$ containing the limit of the stabilizer  $\langle A\rangle$. For Chabauty convergence  deck transformations corresponding to $\beta$ either limit to elements of  $\Gamma'$ or have axes with distance from $\tilde\alpha$ tending to infinity.  In the latter case since $f(u)\approx e^{-2d}$ the limiting contribution to the  sum is zero.  For the former case each deck transformation of $\Gamma'$ corresponding to $\beta$ is a limit of deck transformations with contributions converging.  The pairing $\langle\grad\ell_{\alpha},\grad\ell_{\beta}\rangle$ is continuous relative to the Chabauty topology.

The considerations for the point $p$ involve the pairing $\langle\grad \ell_{\alpha},i\,\grad \ell_{\beta}\rangle $ only for the circumstance of $\alpha$ disjoint from $\beta$. The twist-length formula of \cite[Thrm. 3.3]{Wlsymp} provides that the pairing vanishes in this situation.  The non zero pairings are all of the form $\langle\grad \ell_{\alpha}, \grad \ell_{\beta}\rangle $.    

The matrix $P$ of WP pairings for a combined short and relative length basis is a Chabauty continuous matrix-valued function. The desired conclusion follows from Lemma \ref{gradest} provided the determinant is non zero at $p$.  Let $\mathcal V$ 
 be the span for the short geodesics basis $\{\lambda_{\alpha},i\lambda_{\alpha}\}_{\alpha\in\sigma}$ and $\mathcal{V}'$
 be the span for the relative length basis $\{\grad\ell_{\beta}\}_{\beta\in\tau}$.  By hypothesis the pairing matrix $P'$ for $\mathcal{V}'$ has non zero determinant at $p$.  From Lemma \ref{gradest} and the twist-length formula the limit of the pairing matrix for $\mathcal V$ at each point of $\mathcal T(\sigma)$ is $\frac{1}{2\pi}I$ and for $\mathcal V$ with $\mathcal V'$ is the null matrix.  We have at $p$ that $\det(P)=(2\pi)^{-\dim\mathcal V}\det(P')$, the desired property.

The matrix $P$ of WP Hermitian pairings for $\{\lambda_{\alpha}\}_{\alpha\in\sigma'}$ is a Chabauty continuous matrix-valued function. The conclusion follows provided the determinant is non zero.  Let 
$\mathcal V$ be the span for the short geodesics basis $\{\lambda_{\alpha} \}_{\alpha\in\sigma}$ and 
$\mathcal V'$ be the span for the basis $\{\lambda_{\alpha} \}_{\alpha\in\sigma'-\sigma}$ for completing a partition.
  The complex differentials $\{\partial\ell_{\alpha}^{1/2}\}_{\alpha\in\sigma'-\sigma}$ provide a global frame over $\mathbb C$ for the tangent space of $\mathcal T(\sigma)$, \cite[Thrm. 3.7]{WlFN}.  The pairing matrix $P'$ for $\mathcal{V}'$ has non zero determinant at $p$.  As above the limit of the pairing matrix for $\mathcal V$  is $\frac{1}{2\pi}I$ and for $\mathcal V$ with $\mathcal V'$ is the null matrix.  We have  at $p$ that $\det(P)=(2\pi)^{-\dim\mathcal V}\det(P')$.  The proof is complete.       

We are ready to present expressions comparable to the WP metric.  Bers showed that a rough fundamental domain for $Mod$ acting on $\mathcal T$ is provided in terms of bounded length partitions.   For suitable positive constants $c_{g,n}$ a rough fundamental domain is the union of regions $\{\ell_{\alpha}\le c_{g,n},\,0\le\vartheta_{\alpha}\le 2\pi \mid \alpha\in\sigma\}$ for $\sigma$ ranging over the finite number of homotopically distinct partitions, \cite{Bersdeg, Busbook}.  Theorem \ref{hessmgl} already provides that $\hess\ell_{\alpha}$ is comparable to the weighted WP pairing $\langle\ ,\ \mathbb P_{\alpha}\rangle $.  We write $f\asymp g$ for tensors provided there exists positive constants $c',c''$ such that $c'g\le f\le c''g$ for all evaluations.  We have the following comparability for bounded length partitions.   
\begin{theorem}
\label{wpcomp}
The WP metric is comparable to a sum of first and second derivatives of  geodesic-length functions for a partition $\sigma$ and $J$ the almost complex structure as follows
\[
\langle\ ,\ \rangle \ \asymp\ \sum_{\alpha\in\sigma}(d\ell_{\alpha}^{1/2})^2+(d\ell_{\alpha}^{1/2}\circ J)^2\ \asymp\ \sum_{\alpha\in\sigma}\hess\ell_{\alpha} 
\]
where given $c_0$ positive there are positive constants $c_1$ and $c_2$ for the comparability on the Bers region $\mathcal B=\{\ell_{\alpha}\le c_0\mid\alpha\in\sigma\}$.
\end{theorem}
\prf There are a finite number of distinct partitions modulo the action of $Mod$.  It is enough to establish the existence of $c_1$ and $c_2$ depending on $c_0$ and the choice of $\sigma$.  The first comparability is equivalent to the boundedness in $GL(\mathbb C)$ of the Hermitian pairing matrix for $\{\partial\ell_{\alpha}^{1/2}\}_{\alpha\in\sigma}$.  The partition $\sigma$ and pairing matrix $\bigl(\langle\partial\ell_{\alpha}^{1/2},\partial\ell_{\alpha'}^{1/2}\rangle _{Herm}\bigr)$ are invariant under the group $Twist(\sigma)$ of Dehn twists about the elements of $\sigma$.  It is enough to consider the pairing matrix on the product of punctured discs 
$\mathcal B/Twist(\sigma)$.  WP metric balls provide a system of relatively compact neighborhoods of  $(\overline{\mathcal B}-\mathcal B)/Twist(\sigma)\subset\overline{\mathcal T}/Twist(\sigma)$.   Lemma \ref{paircont} provides for the boundedness on sufficiently small closed metric balls.  The region $\overline{\mathcal B}/Twist(\sigma)$ is compact and a finite collection of suitable metric balls is selected to cover $(\overline{\mathcal B}-\mathcal B)/Twist(\sigma)$.    The complex differentials $\{\partial\ell_{\alpha}^{1/2}\}_{\alpha\in\sigma}$ also provide a global frame over $\mathbb C$ for the tangent bundle of $\mathcal T$, \cite[Thrm. 3.7]{WlFN}.   The pairing is bounded in $GL(\mathbb C)$ on the complement of the finite collection.  The first comparability is established.  

The second comparability is a consequence of the bounds for the Hessian.  From Corollary \ref{first} the real and complex Hessians are uniformly comparable with 
$8\,\partial\ell ^{1/2}\overline{\partial}\ell ^{1/2}\le\partial\overline{\partial}\ell$.  From Theorem \ref{qshess} the Hessian $\partial\overline{\partial}\ell$ is uniformly bounded in terms of $\sum_{\alpha}|\partial\ell_{\alpha}^{1/2}|^2$ and the WP metric for bounded geodesic-lengths.  The proof is complete.  

From Lemma \ref{gradest}, the twist-length vanishing and Theorem \ref{qshess} we find the following.
\begin{corollary}
\label{wpexpan}
The WP metric has the geodesic-length differential and Hessian expansions
\begin{align*}
\langle\ ,\ \rangle \ =&\ 2\pi\sum_{\alpha\in\sigma}(d\ell_{\alpha}^{1/2})^2+(d\ell_{\alpha}^{1/2}\circ J)^2 \ +\ O((\sum_{\alpha\in\sigma}\ell_{\alpha}^3)\,\langle\ ,\ \rangle )\\ =&\ \frac{\pi}{6}\sum_{\alpha\in\sigma}\frac{\hess\ell_{\alpha}^2}{\ell_{\alpha}} \ +\ O((\sum_{\alpha\in\sigma}\ell_{\alpha}^2)\,\langle\ ,\ \rangle )
\end{align*}
at the maximal frontier point of $\overline{\mathcal T}$ for the partition $\sigma$.
\end{corollary}

We present a local form of the comparability for a relative length basis.  
\begin{corollary}
\label{wpcomp2}
The WP metric is comparable to a sum of differentials of geodesic-length functions for a simplex $\sigma$ of short geodesics and corresponding relative length basis $\tau$ as follows
\[
\langle\ ,\ \rangle \ \asymp\ \sum_{\alpha\in\sigma}(d\ell_{\alpha}^{1/2})^2+(d\ell_{\alpha}^{1/2}\circ J)^2\ + \sum_{\beta\in\tau}(d\ell_{\beta})^2
\]
where given  $p\in\mathcal T(\sigma)$ and $\sigma,\,\tau$ there is a neighborhood $U$ 
of $p$ in $\mathcal T\subset\overline{\mathcal T}$ and positive constants $c_1$ and $c_2$ for the comparability.
\end{corollary}

We sketch an application for the compactified moduli space of stable curves $\overline{\mathcal M}=\Tbar\slash Mod$.  The compactified moduli space is a complex $V$-manifold and projective variety, \cite{DM}.  The quotient $\mathcal M=\mathcal T\slash Mod$ is the classical moduli space of curves and the compactification locus $\mathcal D=\overline{\mathcal M}-\mathcal M$ is the divisor of noded stable curves.  The divisor $\mathcal D$ is stratified corresponding to the $Mod$-orbits of the augmentation locus of $\Tbar$.  A point of $\mathcal D$ has a neighborhood with a local manifold cover given as the product of an open set in a stratum $\mathcal T(\sigma)$ and parameters for the transverse to $\mathcal D$. A point of $\mathcal T(\sigma)$ has nodes labeled by the elements of $\sigma$.  For each node of a stable curve a local plumbing family $\{zw=t\}\rightarrow\{|t|<1\}$ contributes a parameter $t$ for the transverse.
  Let $\mathcal D(\sigma)$ be a stratum of $\mathcal D$ represented  by a stratum $\mathcal T(\sigma)\subset\Tbar$.  The complex tangent bundle of $\mathcal D(\sigma)\subset\overline{\mathcal M}$ has a local extension as follows. Let $\mathcal N(\sigma)$ be the subbundle of the tangent bundle of $\overline{\mathcal M}$ with frame $\{\lambda_{\alpha},J\lambda_{\alpha}\}_{\alpha\in\sigma}$.  As noted above for a partition $\sigma'\supset\sigma$ the complex differentials $\{\partial\ell_{\alpha}^{1/2}\}_{\alpha\in\sigma'}$ provide a local frame for the tangent bundle of $\overline{\mathcal M}$.  The complex tangent bundle of $\mathcal D(\sigma)$ is characterized as the restriction of the orthogonal complement $\mathcal N(\sigma)^{\perp}$ to the locus $\mathcal D(\sigma)$.  We write $\langle\ ,\ \rangle_{\mathcal N(\sigma)^{\perp}}$ for the restriction of the WP metric to $\mathcal N(\sigma)^{\perp}$.  The restriction of the WP metric to the tangent space of $\mathcal D(\sigma)$ coincides with the lower dimensional WP metric of $\mathcal T(\sigma)$, \cite{Msext, Wlcomp}.  The WP metric for $\mathcal N(\sigma)^{\perp}$ is an extension of the lower dimensional WP metric.    The expansion $\ell_{\alpha}=2\pi^2(-\log |t_{\alpha}|)^{-1}+O((\log |t_{\alpha}|)^{-2})$ for the geodesic-length in terms of the plumbing parameter \cite[Example 4.3]{Wlhyp} and the expansion $d_{\mathcal T(\sigma)}=(2\pi\sum_{\alpha\in\sigma}\ell_{\alpha})^{1/2}+O(\sum_{\alpha\in\sigma}\ell_{\alpha}^{5/2})$ for the distance to $\mathcal T(\sigma)$ (see Corollary \ref{termbhv} below) suggest the following. 

\begin{corollary}
The geodesic-length functions for the simplex $\sigma$ provide an expansion for the WP K\"{a}hler potential for the transverse to $\mathcal D(\sigma)\subset\overline{\mathcal M}$.  The WP metric has the expansion
\[
\langle\ ,\ \rangle\ =\ \langle\ ,\ \rangle_{\mathcal N(\sigma)^{\perp}}\ +\ \pi\partial\overline{\partial}\sum_{\alpha\in\sigma}\ell_{\alpha}\ +\ O((\sum_{\alpha\in\sigma}\ell_{\alpha}^2)\langle\ ,\ \rangle)
\]
in a local manifold cover for a neighborhood of a point of $\mathcal D(\sigma)$ where for $c_0$ positive the remainder term constant is uniform for $\ell_{\alpha}\le c_0,\,\alpha\in\sigma$.  
\end{corollary}
\prf We consider the orthogonal decomposition $\mathcal N(\sigma)\oplus N(\sigma)^{\perp}$.  The considerations for Corollary \ref{wpexpan} provide for the expansion $\ \langle\ ,\ \rangle_{\mathcal N(\sigma)}=8\pi\sum_{\alpha\in\sigma}|\partial\ell_{\alpha}^{1/2}|^2+O(\sum_{\alpha}\ell_{\alpha}^3\langle\ ,\ \rangle)$.  Theorem \ref{qshess} provides for the expansion $8|\partial\ell_{\alpha}^{1/2}|^2=\partial\overline{\partial}\ell_{\alpha}+O(\ell_{\alpha}^2\langle\ ,\ \rangle)$.  The proof is complete.

The expansion for $\ell_{\alpha}$ also suggests consideration of $\partial\overline{\partial}(-\ell_{\alpha}^{-1})$.  The complex Hessian is positive definite by Corollary \ref{first} and by Theorem \ref{qshess} is bounded by the WP metric and a constant uniform for $\ell_{\alpha}\le c_0$.        

\subsection{A formula for the WP connection}
\label{WPconn}
We present an expansion for the WP covariant derivative of gradients of geodesic-length functions.  The formula enables further study of gradients of geodesic-length functions, calculation of their Lie brackets and an examination of geodesics ending at the augmentation $\Tbar-\mathcal T$. To show the formal nature of the connection for geodesic-length functions we continue and write $J$ for the almost complex structure of $\mathcal T$ and for a geodesic $\gamma$ for the hyperbolic surface we write $\lambda_{\gamma}=\grad\ell_{\gamma}^{1/2}$ for the root geodesic-length gradient.  The Riemannian Hessian is directly related to covariant differentiation $\hess h(U,V)=\langle D_U\grad h,V\rangle $ for vector fields $U,V$, \cite{BO'N}.  
\begin{theorem}
\label{D}
The WP covariant derivative $D$ of the root geodesic-length gradient $\lambda_{\alpha}$ satisfies
\[
D_U\lambda_{\alpha}\ =\ 3\ell_{\alpha}^{-1/2}\langle J\lambda_{\alpha},U\rangle J\lambda_{\alpha}\ +\ O(\ell_{\alpha}^{3/2}\|U\|)
\]
in terms of the WP norm $\|\ \|$ where for $c_0$ positive the remainder term constant is uniform for $\ell_{\alpha}\le c_0$.
\end{theorem}
\prf We have from Theorem \ref{qshess} that
\begin{multline*}
\hess \ell_{\alpha}^{1/2}(U,V)=\frac34 \ell_{\alpha}^{-3/2}(\dot\ell_{\alpha}\circ J)^2\ +
\ O(\ell_{\alpha}^{3/2}\|U\|\|V\|)\\
=3\ell_{\alpha}^{-1/2}\langle J\lambda_{\alpha},U\rangle \langle J\lambda_{\alpha},V\rangle \ +\ O(\ell_{\alpha}^{3/2}).
\end{multline*}
The conclusion now follows from the characterization of $D_U\grad h$ as the dual in the inner product of the linear functional $\hess H(U,V)$ of $V$.  The proof is complete.

We use the formula to study the covariant derivative and behavior of WP geodesics in a neighborhood of a point $p$ of $\mathcal T(\sigma)\subset\Tbar$ for the simplex $\sigma$.  Basic properties of a K\"{a}hler metric include that $J$ is orthogonal, that $J$ is parallel $D_UJV=JD_UV$, and the commutator relation $[U,V]=D_UV-D_VU$.  We continue following Definition \ref{relbas} and consider $\mathcal S=\{\lambda_{\alpha}\}_{\alpha\in\sigma}$ for the short geodesic-lengths and $\mathcal R=\{\grad\ell_{\beta}\}_{\beta\in\tau}$ for a relative length basis. The collection 
$\{\lambda_{\alpha},J\lambda_{\alpha},\grad\ell_{\beta}\}_{\alpha\in\sigma,\,\beta\in\tau}$  provides a local tangent frame for $\mathcal T$ in a neighborhood of $p$.  Lemma \ref{gradest} describes the inner products of elements of the frame.  We further consider the behavior of the frame.   

The elements of $J\mathcal S$ have the geometric description as the infinitesimal Fenchel-Nielsen right twist deformations, \cite[display (1.1), Coro. 2.8, Thrm. 2.10]{WlFN}.  In particular $2\ell_{\alpha}^{1/2}J\lambda_{\alpha}=J\grad\ell_{\alpha}=2t_{\alpha}$ and $J\lambda_{\alpha}=\ell_{\alpha}^{-1/2}t_{\alpha}$ for the infinitesimal right twist normalized for the hyperbolic distance between the designated reference points having unit derivative.  
The Lie bracket for the rescaled twist deformations $\widetilde{T_{\gamma}}=4\ell_{\gamma}^{1/2}\sinh\ell_{\gamma}J\lambda_{\gamma}$ has the structure of a Lie algebra over the integers, \cite[Thrm. 4.8]{Wlsymp}.
For Fenchel-Nielsen coordinates the angle is defined as $\vartheta_{\alpha}=2\pi\tau_{\alpha}/\ell_{\alpha}$ with an increment of $2\pi$ corresponding to a full rotation (a Dehn twist.)  The corresponding Fenchel-Nielsen infinitesimal angle variation is given as $T_{\alpha}=(2\pi)^{-1}\ell_{\alpha}^{3/2}J\lambda_{\alpha}=(2\pi)^{-1}\ell_{\alpha}t_{\alpha}$ for which  $\vartheta_{\alpha}$ has unit derivative.    From Lemma \ref{gradest} the norm $\|T_{\alpha}\|$ is bounded as $O(\ell_{\alpha}^{3/2})$.  We abbreviate notation for the simplex $\sigma$ and write $\mathcal A=\{T_{\alpha}\}_{\alpha\in\sigma}$ for the Fenchel-Nielsen angle variations, $\mathcal S=\{\lambda_{\alpha}\}_{\alpha\in\sigma}$ for the short geodesic root length gradients and $\mathcal R=\{\grad\ell_{\beta}\}$ for the relative length basis gradients.  The elements of $\sigma$ are disjoint and the angle variations $T_{\alpha},\alpha\in\sigma$ commute.  

\begin{corollary}
\label{comm}
The Lie brackets of the elements of $\mathcal S$ and $\mathcal R$ on $\mathcal T$ in a neighborhood of $p\in\mathcal T(\sigma)\subset\Tbar$ are bounded as follows: $[T_{\alpha},\lambda_{\alpha'}]$ is $O(\ell_{\alpha}^3+\ell_{\alpha}^{3/2}\ell_{\alpha'}^{3/2})$, $[T_{\alpha},\grad\ell_{\beta}]$ is $O(\ell_{\alpha}^{3/2})$, $[\lambda_{\alpha},\lambda_{\alpha'}]$ is $O(\ell_{\alpha}^{3/2}+\ell_{\alpha'}^{3/2})$, $[\lambda_{\alpha},\lambda_{\beta}]$ is $O(\ell_{\alpha}^{3/2}+\ell_{\beta}^{3/2})$ and $[\lambda_{\beta},\lambda_{\beta'}]$ is $O(\ell_{\beta}^{-1/2}+\ell_{\beta'}^{-1/2})$, for $\alpha,\alpha'\in\sigma$ and $\beta,\beta'\in\tau$.  On $\mathcal T$ the elements of $\mathcal S$ have covariant derivatives with projections onto the span of $\mathcal S$ bounded as $O(\sum_{\alpha\in\sigma}\ell_{\alpha}^{3/2})$.  For $c_0$ positive all constants are uniform for $\ell_{\alpha},\ell_{\beta}\le c_0$, $\alpha\in\sigma$ and $\beta\in\tau$.
\end{corollary}
\prf  We begin with $[T_{\alpha},\lambda_{\alpha'}]=D_{T_{\alpha}}\lambda_{\alpha'}-D_{\lambda_{\alpha'}}T_{\alpha}$ where from Theorem \ref{D} and the definition of $T_{\alpha}$
\begin{align*}
D_{T_{\alpha}}\lambda_{\alpha'}=&3\ell_{\alpha'}^{-1/2}\langle J\lambda_{\alpha'},T_{\alpha}\rangle J\lambda_{\alpha'}\ +\ O(\ell_{\alpha}^{3/2}\ell_{\alpha'}^{3/2})\\
=&3(2\pi)^{-1}\ell_{\alpha'}^{-1/2}\ell_{\alpha}^{3/2}\langle\lambda_{\alpha'},\lambda_{\alpha}\rangle J\lambda_{\alpha'}\ +\ O(\ell_{\alpha}^{3/2}\ell_{\alpha'}^{3/2})
\end{align*}
and
\begin{align*}
D_{\lambda_{\alpha'}}T_{\alpha}=&(2\pi)^{-1}(\lambda_{\alpha'}\ell_{\alpha}^{3/2})J\lambda_{\alpha}\ +\ (2\pi)^{-1}\ell_{\alpha}^{3/2}JD_{\lambda_{\alpha'}}\lambda_{\alpha}\\
=&3(2\pi)^{-1} \ell_{\alpha}\langle\lambda_{\alpha'},\lambda_{\alpha}\rangle J\lambda_{\alpha}\ +\ (2\pi)^{-1}\ell_{\alpha}^{3/2}JD_{\lambda_{\alpha'}}\lambda_{\alpha}.
\end{align*}
The desired bound for $[T_{\alpha},\lambda_{\alpha'}]$ now follows from Lemma \ref{gradest} and the twist-length vanishing, \cite[Thrm. 3.3]{Wlsymp}.  The bound for $[T_{\alpha},\grad\ell_{\beta}]$ follows similarly from Theorem \ref{D} and Lemma \ref{gradest}.  The bounds for $[\lambda_*,\lambda_{**}]$ follow from Theorem \ref{D}, 
Lemma \ref{gradest} and the twist-length vanishing.  For the covariant derivatives of elements of $\mathcal S$ we consider
$D_U\lambda_{\alpha}=3\ell_{\alpha}^{-1/2}\langle J\lambda_{\alpha},U\rangle \operatorname{proj}_{\mathcal S}J\lambda_{\alpha}\ +\ O(\ell_{\alpha}^{3/2}\|U\|)$
and note that the projection $\operatorname{proj}_{\mathcal S}J\lambda_{\alpha}$ onto the span of $\mathcal S$ vanishes by the twist-length vanishing.  The proof is complete.

We are ready to consider geodesics that end at the augmentation point $p$ in $\mathcal T(\sigma)\subset\Tbar$.  Yamada discovered that the augmented Teichm\"{u}ller space is a complete $CAT(0)$ metric space with strata $\mathcal T(\sigma)$ forming convex subsets, \cite{DW2, MW, Wlcomp, Yam2}.  A closed convex set in a $CAT(0)$ metric space is the base of a distance non increasing projection, \cite[Chap. II.2, Prop. 2.4]{BH}.  Accordingly there is a projection $\operatorname{proj}_{\mathcal T(\sigma)}$ of $\Tbar$ onto the closure of 
$\mathcal T(\sigma)$.  The unique geodesic from $q\in\Tbar$ to $\operatorname{proj}_{\mathcal T(\sigma)}(q)$ realizes distance from its points to the closure of $\mathcal T(\sigma)$.  
\begin{definition}
\label{termproj}
A terminating WP geodesic has endpoint on a stratum $\mathcal T(\sigma)$.  A projecting WP geodesic realizes distance between its points and a stratum $\mathcal T(\sigma)$.  
\end{definition} 

We consider geodesic-length functions and tangent fields for unit-speed geodesics $\kappa$ of $\Tbar$ ending at the augmentation point $p$.  Such WP geodesics $\kappa$ are parameterized in $t$ for $0\le t\le t_0$ with $\kappa(0)=p$.  For the topology of $\Tbar$ and geodesics $\gamma$ on the hyperbolic surface not transversely intersecting the short geodesics $\alpha\,,\alpha\in\sigma$, the geodesic-length functions $\ell_{\gamma}$ are continuous on $\kappa$ and real analytic for $0<t\le t_0$.  Furthermore by Corollary \ref{first} the root geodesic-length functions $\ell_{\gamma}^{1/2}$ are convex along $\kappa$.  It follows that on $\kappa$ the functions $\ell_{\gamma}^{1/2}$ have one-sided derivatives at $t=0$.  And as noted above the collection $\mathcal S\cup J\mathcal S\cup\mathcal R$ provides a local tangent frame for $\mathcal T$ in a neighborhood of $p\in\mathcal T(\sigma)\subset\Tbar$.  We have the following description of the tangent field of a terminating geodesic.   
 
\begin{corollary}
\label{termbhv}
A geodesic $\kappa$ that terminates at $p$ has tangent field $\dt$ almost in the span of $\mathcal S\cup\mathcal R$ as follows: $\|\operatorname{proj}_{J\mathcal S}\dt\|$ is $O(t\sum_{\alpha\in\sigma}\ell_{\alpha}^{3/2})$ for small $t$ and $\operatorname{proj}_{J\mathcal S}$ the projection onto the span of $J\mathcal S$.  On a terminating geodesic $\kappa$ a geodesic-length function $\ell_{\alpha},\,\alpha\in\sigma$, either is identically zero or satisfies 
$\frac{d\ell_{\alpha}^{1/2}}{dt}(t)=a_{\alpha}\,+\,O(t\ell_{\alpha}^{3/2})$ for small $t$ with $a_{\alpha}=\frac{d\ell_{\alpha}^{1/2}}{dt}(0^+)$ nonzero.

A unit-speed projecting geodesic $\varsigma$ has tangent field almost constant in the span of $\mathcal S$ as follows:
\[
\dt= (2\pi)^{-1}\sum_{\alpha\in\sigma}a_{\alpha}\lambda_{\alpha}\ +\ O(t^2 \sum_{\alpha\in\sigma}\ell_{\alpha})
\]
for small $t$ with $(2\pi)^{1/2}\|(a_{\alpha})\|_{Euclid}=1$ in terms of the Euclidean norm.   The distance $d_{\mathcal T(\sigma)}$ on $\Tbar$ to the closure of $\mathcal T(\sigma)$ satisfies
\[
d_{\mathcal T(\sigma)}\le(2\pi\sum_{\alpha\in\sigma}\ell_{\alpha})^{1/2}\quad\mbox{and}\quad d_{\mathcal T(\sigma)}=(2\pi\sum_{\alpha\in\sigma}\ell_{\alpha})^{1/2}\, +\, O(\sum_{\alpha\in\sigma}\ell_{\alpha}^{5/2}).  
\]
For $c_0$ positive all constants are uniform for $\ell_{\alpha}\le c_0,\alpha\in\sigma$.  
\end{corollary}
\prf The first matter is the estimate for the pairing $\langle\lambda_{\alpha},J\dt\rangle $ for a terminating geodesic.  We introduce the function $f(t)=\langle\grad\ell_{\alpha}^2,J\dt\rangle $ for which 
$f'(t)=\langle D_{\dt}\grad\ell_{\alpha}^2,J\dt\rangle =\hess \ell_{\alpha}^2(\dt,J\dt)$. From Lemma \ref{gradest} we note that $f(0)=0$.  From Theorem \ref{qshess} we observe that $\hess \ell_{\alpha}^2$ has principal term $12\ell_{\alpha}(\langle\lambda_{\alpha},\dt\rangle \langle\lambda_{\alpha},J\dt\rangle \, +\, \langle J\lambda_{\alpha},\dt\rangle \langle J\lambda_{\alpha},J\dt\rangle )$ vanishing and consequently $\hess \ell_{\alpha}^2$ is bounded as $O(\ell_{\alpha}^3)$.  Since $\ell_{\alpha}$ is increasing in $t$ it follows on integration in $t$ that $f(t)$ is bounded as $O(t\ell_{\alpha}^3)$ and since $\grad\ell_{\alpha}^2=2\ell_{\alpha}^{3/2}\lambda_{\alpha}$, it further follows that 
$\langle\lambda_{\alpha},J\dt\rangle $ is bounded as $O(t\ell_{\alpha}^{3/2})$, the first conclusion.  It then follows from Theorem \ref{D} that $\hess\ell_{\alpha}^{1/2}(\dt,\dt)$ is bounded as $O(\ell_{\alpha}^{3/2})$ and on integration that $\frac{d\ell_{\alpha}^{1/2}}{dt}(t)=a_{\alpha}\,+\, O(t\ell_{\alpha}^{3/2})$ for $a_{\alpha}$ the initial one-sided derivative of $\ell_{\alpha}^{1/2}$. We pose that the derivative is nonzero provided $\ell_{\alpha}$ is not identically zero.  Proceeding by contradiction if $a_{\alpha}=0$ then $\ell_{\alpha}^{1/2}$ satisfies the description of $g(t)$ increasing with $g(0)=0$ and $\mid\frac{dg}{dt}\mid\le Ctg(t)^3$.  Integration provides $g(t_0)^{-2}-g(t_1)^{-2}\le C(t_1^2-t_0^2)$ which implies $g(t)^{-2}$ is bounded as $t$ tends to zero in contradiction 
to $g(0)=0$, as desired.  A further integration of $\frac{d\ell_{\alpha}^{1/2}}{dt}$ provides that
$\ell_{\alpha}^{1/2}=a_{\alpha}t\ +\ O(t^2\ell_{\alpha}^{3/2})$ and
\[
(2\pi\sum_{\alpha\in\sigma}\ell_{\alpha})^{1/2}=(2\pi)^{1/2}\|(a_{\alpha})\|_{Euclid}\,t\ +\ O(t^2\sum_{\alpha\in\sigma}\ell_{\alpha}^{3/2})
\]
for the Euclidean norm of the vector $(a_{\alpha})_{\alpha\in\sigma}$.

We further consider the projecting geodesic $\varsigma$.  We recall the prior expansion $d_{\mathcal T(\sigma)}=(2\pi\sum_{\alpha\in\sigma}\ell_{\alpha})^{1/2}\,+ \, O(\sum_{\alpha\in\sigma}\ell_{\alpha}^2)$ for the distance to the stratum $\mathcal T(\sigma)$, \cite[Coro. 21]{Wlcomp}.  Specialize by restricting $d_{\mathcal T(\sigma)}$ to the unit-speed projecting geodesic $\varsigma$ for which $d_{\mathcal T(\sigma)}(\varsigma(t))=t$ and substitute the expansion from the above paragraph to derive the unity equation $1=(2\pi)^{1/2}\|(a_{\alpha})\|_{Euclid}$. The expansion for $d_{\mathcal T(\sigma)}$ is established.  The inequality for distance and geodesic-length follows since on $\varsigma$ the functions  $d_{\mathcal T(\sigma)}$ and $(2\pi\sum_{\alpha\in\sigma}\ell_{\alpha})^{1/2}$ have the same initial derivative with the former linear and the latter convex.  Finally we consider the vector field $V=\dt-(2\pi)^{-1}\sum_{\alpha\in\sigma}a_{\alpha}\lambda_{\alpha}$ along $\varsigma$.   Noting that $\dt$ is parallel along the geodesic and applying Theorem \ref{D} with the bound that $\langle J\lambda_{\alpha},\dt\rangle $ is $O(t\ell_{\alpha}^{3/2})$ we observe that $D_{\dt}V$ is bounded as $O(t\sum_{\alpha\in\sigma}\ell_{\alpha})$ and on integration that $V$ is bounded as $O(t^2\sum_{\alpha\in\sigma}\ell_{\alpha})$, the final bound.  The proof is complete. 

We now apply the considerations to give an improvement for Lemma \ref{paircont}.  We adopt the convention that at an augmentation point the WP pairings with $\lambda_{\alpha}$ and $J\lambda_{\alpha}$  are determined by continuity.

\begin{corollary}
\label{WPcont}
The WP pairing matrix $P$ for the frame $\{\lambda_{\alpha},J\lambda_{\alpha},\grad \ell_{\beta}\}$ for $\alpha\in\sigma,\,\beta\in\tau$ satisfies $P(q)=P(p)+O(d(p,q))$ in a neighborhood of the point $p$.
\end{corollary}
\prf  First, the twist-length vanishing provides that the vectors $\{J\lambda_{\alpha}\}_{\alpha\in\sigma}$ are orthogonal to $\{\lambda_{\alpha},\grad \ell_{\beta}\}_{\alpha\in\sigma,\,\beta\in\tau}$.  Second, Lemma \ref{gradest} and Corollary \ref{termbhv} provide that the pairings $\langle\lambda_{\alpha},\lambda_{\alpha'}\rangle $ and $\langle\lambda_{\alpha},\grad \ell_{\beta}\rangle $ have an expansion with an $O(d(p,q)^3)$ remainder.  Third, to analyze $\langle\grad \ell_{\beta},\grad \ell_{\beta'}\rangle $ we consider its derivative along the geodesic connecting $p$ to $q$ and apply Theorem \ref{D} to estimate the covariant derivatives.  The geodesic-lengths $\{\ell_{\beta}\}_{\beta\in\tau}$ are bounded in a neighborhood and consequently the derivative of the pairing is bounded.  The proof is complete. 

We now use properties of geodesics to describe families of surfaces with $\langle\grad\ell_{\alpha},\grad\ell_{\alpha}\rangle\ge e^{\ell_{\alpha}/2}$ in comparison to Lemma \ref{Pbound2}.  Begin with the family of hyperbolic pants with two boundaries with a common length and combine the boundaries to form a one-holed torus with a reflection symmetry fixing the combined geodesic $\beta$.  Let $\mathbf T$ be the marked family of one-holed tori with reflection.  Let $\alpha$ be the geodesic which crosses the $\beta$ collar and intersects $\beta$ orthogonally.  In general consider a sub locus of Teichm\"{u}ller space of surfaces with reflection containing the family of tori $\mathbf T$ as marked subsurfaces.  Further consider a terminating geodesic in the sub locus with unit-speed parameter $t$ and $\ell_{\beta}$ tending to zero.  From the description of collars and the Gromov-Hausdorff compactness of the {\em thick} regions it follows that along the geodesic $\ell_{\alpha}=2\log 1/\ell_{\beta}+O(1)$.  From the above Corollary it follows for the geodesic that $\ell_{\beta}=ct^2+O(t^5)$.  From convexity and the mean value theorem it follows that the derivative $\frac{d\ell_{\alpha}}{dt}$ is comparable to $-1/t$ or equivalently comparable to $-\ell_{\beta}^{-1/2}$.  We have from the relation between $\ell_{\alpha}$ and $\ell_{\beta}$ that $(\frac{d\ell_{\alpha}}{dt})^2$ is comparable to $e^{\ell_{\alpha}/2}$ as proposed.   
  
\subsection{The Alexandrov tangent cone at the augmentation}
\label{Alextgtcn}
We consider the infinitesimal structure in a neighborhood of the augmentation.  The augmentation locus has codimension at least $2$ with neighborhoods not locally compact.   We find that the Alexandrov tangent cone is suited to describe the behavior of geodesics terminating at the augmentation.   In Theorem \ref{Alextgt} we present an isometry from the Alexandrov tangent cone to the product of a Euclidean orthant and the tangent space of the stratum with WP metric. The Alexandrov tangent cone has the structure of a cone in an inner product space.   The isometry is in terms of root geodesic-length functions.  The analysis combines the Euclidean geometry of Fenchel-Nielsen angle variations about disjoint geodesics, the expansion of the WP pairing for root geodesic-lengths, the expansion for the WP connection and basic properties for a $CAT(0)$ geometry.   Throughout the following discussion we consider the short geodesic root lengths $\{\ell_{\alpha}^{1/2}\}_{\alpha\in\sigma}$ and following Definition \ref{relbas} a relative length basis $\{\ell_{\beta}\}_{\beta\in\tau}$.  

\begin{definition} A length number for an augmentation point $p$ of $\mathcal T(\sigma)$ is a tuple $\mathcal L=(\ell_{\alpha}^{1/2},\ell_{\beta}^{1/2})_{\alpha\in\sigma,\,\beta\in\tau}$ of root geodesic-length functions for the short geodesics $\sigma$ and a relative length basis $\tau$.  A length number variation for the point $p$ is the tuple of gradients $(\lambda_{\alpha},\lambda_{\beta})_{\alpha\in\sigma,\,\beta\in\tau}$.
\end{definition} 
In the following without distinguishing the index sets $\sigma$ and $\tau$ we may write $(\ell_{\alpha}^{1/2},\ell_{\beta}^{1/2})$ to denote length numbers and $(\lambda_{\alpha},\lambda_{\beta})$ to denote length numbers variation;  we may also write $(\ell_{\alpha}^{1/2})$ and $(\lambda_j)$ for the tuples.  Beginning properties are as follows.  Fenchel-Nielsen angles are defined by introducing a partition containing $\sigma$.  Length numbers $(\ell_{\alpha}^{1/2},\ell_{\beta}^{1/2})$ in combination with a partial collection $(\vartheta_{\alpha})_{\alpha\in\sigma}$ of Fenchel-Nielsen angles provide local coordinates at the point $p$. Equivalently the length numbers $(\ell_{\alpha}^{1/2},\ell_{\beta}^{1/2})$ determine the corresponding Riemann surface modulo the Fenchel-Nielsen angles for $\alpha\in\sigma$, \cite{Abbook, Busbook, ImTan, MasFN, WlFN}.  Let $|\sigma|$ denote the number of short geodesics and $Flow(\sigma)$ the Lie group isomorphic to $\mathbb R^{|\sigma|}$ of Fenchel-Nielsen angle variations.  The $Flow(\sigma)$ invariant neighborhoods of the point $p$ provide a neighborhood basis.  For a $Flow(\sigma)$ invariant neighborhood $\mathcal N$ length numbers $(\ell_{\sigma}^{1/2},\ell_{\beta}^{1/2})$ provide coordinates for the differentiable manifold with boundary $\mathcal N/Flow(\sigma)$.  As already noted $\{\ell_{\alpha}=0\}$ is locally a product of lower dimensional Teichm\"{u}ller spaces. For $\ell_p=(\ell_{\beta}^{1/2})_{\beta\in\tau}\in\mathbb R^{\dim \mathcal T-2|\sigma|}$ then the length number $(\ell_{\alpha}^{1/2},\ell_{\beta}^{1/2})$ image of $\mathcal N$ is a neighborhood of $(0,\ell_p)$ in $\mathbb R_{\ge 0}^{|\sigma|}\times\mathbb R^{\dim \mathcal T-2|\sigma|}$.      

Given WP geodesics $\gamma_0,\gamma_1$ terminating at $p$ we will introduce $2$-parameter families to interpolate between the pair.  The interpolation is based on the observation from Corollary \ref{termbhv} that on a terminating geodesic the short geodesic-lengths are almost linear in WP arc length.  The first family will interpolate $\gamma_0$ to a curve $\gamma_*$ with length numbers matching $\gamma_1$. The second family will interpolate the Fenchel-Nielsen angles of $\gamma_*$ and $\gamma_1$.   Begin with the convex combination of length numbers $F(s,t)=(1-s)\mathcal L(\gamma_0(t))+s\mathcal L(\gamma_1(t))$, $0\le s\le 1$, contained in $\mathbb R_{\ge 0}^{|\sigma|}\times\mathbb R^{\dim \mathcal T-2|\sigma|}$.  
From the above  for $t$ sufficiently small, $0\le s\le 1$, the family $F(s,t)$ lies in the range of $\mathcal L$ on $\mathcal N/Flow(\sigma)$.
The inverse image $\mathcal L^{-1}(F(s,t))$ is a differentiable $2$-parameter family in 
$\mathcal N/Flow(\sigma)$ or equivalently $\mathcal F=\mathcal L^{-1}(F(s,t))\subset\Tbar$ is a $2+|\sigma|$-parameter differentiable family parameterized by $0\le s\le 1,\,0< t\le t'$ and  $\mathbb R^{|\sigma|}$.   Points $q$ in the family $\mathcal F$ are characterized by $\mathcal L(q)=F(s,t)$ for some $0\le s\le1$ and 
$0\le t\le t'$.  We parameterize the map $\mathcal F$ for the parameter vector field $\frac{\partial}{\partial s}$ to push forward to the orthogonal complement of $\{J\lambda_{\alpha}\}_{\alpha\in\sigma}$.  Below we introduce a second vector field on the map $\mathcal F$.  We now assume that $\gamma_0(t)$ and $\gamma_1(t)$ satisfy the generic condition that each $\ell_{\alpha}^{1/2},\,\alpha\in\sigma,$ restricted to the geodesics is not identically zero and thus by Corollary \ref{termbhv} has non zero initial derivative.  The condition ensures that for $t$ positive the WP geodesic segments and image of $\mathcal F$ lie interior to $\mathcal T$.
\begin{definition}
\label{pQ}
The pseudo geodesic variation field $\mathcal Q=\sum_{\sigma\cup\tau}a_k\lambda_k$ satisfies $\mathcal Q\ell_j^{1/2}=(\mathcal L(\gamma_1(t))-\mathcal L(\gamma_0(t)))^{(j)}$, $\mathcal Q\perp J\lambda_{\alpha}$, $\alpha\in\sigma$, $t>0$ for the short geodesics $\sigma$, the relative length basis $\tau$ and $j^{th}$ tuple component.
\end{definition}
The generic condition ensures that for $t$ positive the image of $\mathcal F$  lies interior to $\mathcal T$ and that there is a flow of $\mathcal Q$ beginning at $\gamma_0$.  On an integral curve of $\mathcal Q$ the length numbers have a constant derivative $\mathcal L(\gamma_1(t_*))-\mathcal L(\gamma_0(t_*))$ for some value $t_*$.  The flow of $\gamma_0(t)$ is a $2$-parameter smooth family $\kappa(s,t)$ satisfying
\begin{multline}
\label{kappa}
\kappa(s,0)=p,\ \kappa(0,t)=\gamma_0(t),\ \frac{\partial}{\partial s}\kappa(s,t)\perp J\lambda_{\alpha},\ \alpha\in\sigma, \\ \mbox{ and }\mathcal L(\kappa(s,t))=(1-s)\mathcal L(\gamma_0(t))+s\mathcal L(\gamma_1(t)) 
\end{multline}
for $0\le s\le 1$ and sufficiently small positive $t$.  The family $\kappa(s,t)$ interpolates $\gamma_0(t)$ to a differentiable curve $\gamma_*(t)$ with length numbers coinciding with $\gamma_1(t)$, in particular $\mathcal L(\gamma_*(t))=\mathcal L(\gamma_1(t))$.
\begin{lemma}
\label{Qbound}
The pseudo geodesic variation field $\mathcal Q$ on the family $\kappa(s,t),\,0\le s\le 1,\,0<t<t'$ is bounded as follows: 
$\mathcal Q$, $[\mathcal Q,\lambda_j],\,j\in\sigma\cup\tau$ and $[\mathcal Q,J\lambda_{\alpha}],\,\alpha\in\sigma$ are  $O(t')$ for the short geodesics $\sigma$ and relative length basis $\tau$.
\end{lemma}
\prf  The derivatives of the root geodesic-lengths by $\mathcal Q$ are given as  $\mathcal Q\ell_j^{1/2}=\sum_ka_k\langle\lambda_k,\lambda_j\rangle$.  From Definition \ref{pQ} on the family $\kappa(s,t)$ the coefficient vector is given as $(a_k)=P^{-1}(\mathcal L(\gamma_1(t))-\mathcal L(\gamma_0(t)))$ for the pairing matrix $P=(\langle\lambda_j,\lambda_k\rangle )$.  From Lemma \ref{paircont} the inverse $P^{-1}$ is continuous at the point $p$.  From Corollary \ref{first} geodesic-lengths are convex along geodesics, 
$\mathcal L$ is initially one-sided differentiable and consequently $(a_k)$ and $\mathcal Q$ are suitably bounded.  Next from Corollary \ref{comm}  the Lie brackets $[\lambda_k,\lambda_j]$ and $[\lambda_k,J\lambda_{\alpha}]$ are bounded.  To show that the Lie brackets with $\mathcal Q$ are suitably bounded it is enough since $\mathcal Q=\sum a_k\lambda_k$ to show that $\lambda_ja_k$ and $J\lambda_{\alpha} a_k$ are bounded as $O(t)$.  Since $(a_k)=P^{-1}(\mathcal L(\gamma_1(t))-\mathcal L(\gamma_0(t)))$ the desired bound will follow provided the derivatives of $P^{-1}$ are bounded and the derivatives of $(\mathcal L(\gamma_1(t))-\mathcal L(\gamma_0(t)))$ are bounded as $O(t)$.  

The pairing matrix $P$ provides a continuous map into $GL(\mathbb R)$ with a bound for the derivative of $P$ providing a bound for the derivative of the inverse.  To consider a derivative $V\langle\lambda_m,\lambda_k\rangle =\langle D_V\lambda_m,\lambda_k\rangle \,+\,\langle\lambda_m,D_V\lambda_k\rangle $ for $V=\lambda_j$ or $J\lambda_{\alpha}$ apply Theorem \ref{D} and consider for example $D_V\lambda_m=3\ell_m^{-1/2}\langle J\lambda_m,V\rangle J\lambda_m+O(\ell_m^{3/2}\|V\|)$.  The principal term has direction $J\lambda_m$.  Either $m\in\sigma$ and by twist-length vanishing the pairing $\langle J\lambda_m,\lambda_k\big>$ vanishes or $m\in\tau$ and $\ell_m^{-1/2}$ is bounded in a neighborhood of the point $p$.  The derivative $\langle D_V\lambda_m,\lambda_k\rangle $ and derivatives of $P$ are bounded.

We consider the derivatives of $\mathcal L$.  As already noted the derivatives 
$\{J\lambda_{\alpha}\}_{\alpha\in\sigma}$ are the infinitesimal Fenchel-Nielsen twist deformations which stabilize $\mathcal L$.  In particular the evaluation $d\mathcal L(J\lambda_{\alpha})$ is zero.  It remains to consider the $\lambda_j$ derivative of the difference of $\mathcal L$ values in particular $\langle\lambda_j,\lambda_k\rangle (\gamma_1(t))-\langle\lambda_j,\lambda_k\rangle (\gamma_0(t))$.   The difference is $O(t)$ from Corollary \ref{WPcont}.  The proof is complete.  

The definition of geodesic-length is intrinsic.  A Fenchel-Nielsen angle is given in terms of the hyperbolic distance between designated points which are determined only after a partition is fixed, 
\cite{Abbook,Busbook,ImTan,MasFN,WlFN}.  The differential of Fenchel-Nielsen angles is characterized by $d\vartheta_{\alpha}(T_{\alpha'})=\delta_{\alpha\alpha'}$ for the infinitesimal angle variations and 
$\frac{\partial\vartheta_{\alpha}}{\partial \ell_{\alpha'}}=0$ for Fenchel-Nielsen coordinates $(\ell_{\alpha'},\vartheta_{\alpha'})$  with the definition of $\frac{\partial\ }{\partial \ell_{\alpha'}}$ requiring the data of the selected partition and designated points.  

We seek to bound the differential of angle and variation of angle along terminating WP geodesics.  For this purpose we introduce a quantity closely related to the Fenchel-Nielsen angle and consider its differential geometric properties.  The intrinsically defined Fenchel-Nielsen angle variation $T_{\alpha}=(2\pi)^{-1}\ell_{\alpha}^{3/2}J\lambda_{\alpha}$ satisfies $d\vartheta_{\alpha}(T_{\alpha})=1$ for the geodesic $\alpha$.
\begin{definition}
\label{gauge}
The Fenchel-Nielsen gauge is the differential $1$-form $\varrho_{\alpha}=2\pi(\ell_{\alpha}^{3/2}\langle\lambda_{\alpha},\lambda_{\alpha}\rangle )^{-1}\langle\ ,J\lambda_{\alpha}\rangle $.
\end{definition}
The gauge satisfies $\varrho_{\alpha}(T_{\alpha})=1$.  We expect that the gauge is not a closed $1$-form and consider its exterior derivative $d\varrho_{\alpha}(U,V)= U\varrho_{\alpha}(V)-V\varrho_{\alpha}(U)-\varrho_{\alpha}([U,V])$ for vector fields $U,V$.  The exterior derivative is a tensor and is determined by evaluation on a basis of tangent vectors.
\begin{lemma}
\label{drho}
The Fenchel-Nielsen gauges satisfy $\sum_{\alpha\in\sigma}\varrho_{\alpha}^2=\sum_{\alpha\in\sigma}(d\vartheta_{\alpha})^2\circ \operatorname{proj}_{\operatorname{Ker}\, d\mathcal L}(1\,+\,O(\sum_{\alpha\in\sigma}\ell_{\alpha}^3))$ for the orthogonal projection $\operatorname{proj}_{\operatorname{Ker}\, d\mathcal L}$ onto the kernel of the differential of length number $\mathcal L$.  The Fenchel-Nielsen gauge exterior derivatives are bounded as: $|d\varrho_{\alpha}(T_{\alpha'},\lambda_j)|$ is $O(\ell_{\alpha'}^{3/2}\ell_{\alpha}^{-3/2})$ and $|d\varrho_{\alpha}(\lambda_j,\lambda_k)|$ is $O(\ell_{\alpha}^{-3/2})$ for $\alpha,\alpha'\in\sigma,$ $j,k\in\sigma\cup\tau$ for the short geodesics $\sigma$ and relative length basis $\tau$.
\end{lemma}
\prf  We consider first the relation between gauges and angles.  The length number $\mathcal L$ and angles $(\vartheta_{\alpha})_{\alpha\in\sigma}$ combine to provide local coordinates at $p$ where from twist-length vanishing $\{T_{\alpha}\}_{\alpha\in\sigma}$ annihilate the differential $d\mathcal L$.  It follows that the kernel of $d\mathcal L$ has dimension $|\sigma|$ and is spanned by the vectors $\{T_{\alpha}\}_{\alpha\in\sigma}$.  It follows from the defining properties of Fenchel-Nielsen angles, angle variations and gauges that $d\vartheta_{\alpha}(T_{\alpha'})=\delta_{\alpha\alpha'}$ and $\varrho_{\alpha}(T_{\alpha})=1$.  It follows from Lemma \ref{gradest}  for $\alpha\ne \alpha'$ that $\varrho_{\alpha}(T_{\alpha'})$ is bounded as $O(\ell_{\alpha}^{3/2}\ell_{\alpha'}^{3/2})$.  The expansion follows.

We consider the exterior derivative of a gauge.  From twist-length vanishing it follows that $\varrho_{\alpha}(\lambda_j)=0,$ $j\in\sigma\cup\tau$ and consequently that $d\varrho_{\alpha}(T_{\alpha'},\lambda_j)=-\lambda_j\varrho_{\alpha}(T_{\alpha'})-\varrho_{\alpha}([T_{\alpha'},\lambda_j])$, $d\varrho_{\alpha}(\lambda_j,\lambda_k)=-\varrho_{\alpha}([\lambda_j,\lambda_k])$ for $\alpha,\alpha'\in\sigma$ and $j,k\in\sigma\cup\tau$.  It follows directly from the definition of gauge and Corollary \ref{comm} that $\varrho_{\alpha}([T_{\alpha'},\lambda_j])$ is bounded as $O(\ell_{\alpha'}^{3/2}\ell_{\alpha}^{-3/2})$ and $\varrho_{\alpha}([\lambda_j,\lambda_k])$ is bounded as $O(\ell_{\alpha}^{-3/2})$.  It only remains to bound the $\lambda_j$ derivative of $\varrho_{\alpha}(T_{\alpha'})=\ell_{\alpha'}^{3/2}\ell_{\alpha}^{-3/2}\langle\lambda_{\alpha},\lambda_{\alpha}\rangle ^{-1}\langle\lambda_{\alpha},\lambda_{\alpha'}\rangle $.  The derivative is zero for $\alpha=\alpha'$.  We continue and consider the term $\lambda_j(\ell_{\alpha'}^{3/2}\ell_{\alpha}^{-3/2})=3\ell_{\alpha'}\ell_{\alpha}^{-3/2}\langle\lambda_j,\lambda_{\alpha'}\rangle \ - \ 3\ell_{\alpha'}^{3/2}\ell_{\alpha}^{-2}\langle\lambda_j,\lambda_{\alpha}\rangle $ and the term $\langle\lambda_{\alpha},\lambda_{\alpha}\rangle ^{-1}\langle\lambda_{\alpha},\lambda_{\alpha'}\rangle $.  By Lemma \ref{gradest} for $\alpha\ne \alpha'$ the product of terms is smaller than $O(\ell_{\alpha'}^{3/2}\ell_{\alpha}^{-3/2})$, the desired bound.  

We next consider the derivative
\begin{multline*}
\lambda_j(\langle\lambda_{\alpha},\lambda_{\alpha}\rangle ^{-1}\langle\lambda_{\alpha},\lambda_{\alpha'}\rangle )=-2\langle\lambda_{\alpha},\lambda_{\alpha}\rangle ^{-2}\langle D_{\lambda_j}\lambda_{\alpha},\lambda_{\alpha}\rangle \langle\lambda_{\alpha},\lambda_{\alpha'}\rangle \\ +\ \langle\lambda_{\alpha},\lambda_{\alpha}\rangle ^{-1}(\langle D_{\lambda_j}\lambda_{\alpha},\lambda_{\alpha'}\rangle \ +\ \langle\lambda_{\alpha},D_{\lambda_j}\lambda_{\alpha'}\rangle ).
\end{multline*}
Again by twist-length vanishing the principal term of $D_{\lambda_j}\lambda_{\alpha}$ vanishes and the covariant derivative is bounded as $O(\ell_{\alpha}^{3/2})$.  The derivative $\lambda_j(\langle\lambda_{\alpha},\lambda_{\alpha}\rangle ^{-1}\langle\lambda_{\alpha},\lambda_{\alpha'}\rangle )$ is bounded as $O(\ell_{\alpha}^{3/2}+\ell_{\alpha'}^{3/2})$.  The bound for $\lambda_j\varrho_{\alpha}(T_{\alpha'})$ now follows.   The proof is complete.  

We review the notion of Alexandrov angle for $CAT(0)$ metric spaces; see \cite[Chaps. I.1, I.2, II.3, esp. Chap. II.3, Prop. 3.1]{BH} for properties of angles in metric spaces.  A triple of points $(p,q,r)$ has Euclidean comparison triangle with angle $\angle (p,q,r)$ at $p$ valued in the interval $[0,\pi]$ determined by the Law of Cosines $2d(p,q)\,d(p,r)\cos\angle (p,q,r)=d(p,q)^2+d(p,r)^2-d(q,r)^2$.  For constant speed geodesics $\gamma_0(t),\gamma_1(t)$ with common initial point from the $CAT(0)$ inequality the comparison angle for  $(\gamma_0(0),\gamma_0(t),\gamma_1(t'))$ is a non decreasing function of $t$ and $t'$; see \cite[Chap. II.3, Prop. 3.1]{BH}.  The Alexandrov angle between the constant speed geodesics $\gamma_0$ and $\gamma_1$ with initial point $p$ is defined by the limit
\[
\cos \angle(\gamma_0,\gamma_1)=\lim_{t\rightarrow 0} \frac{d(p,\gamma_0(t))^2\,+\, d(p,\gamma_1(t))^2-d(\gamma_0(t),\gamma_1(t))^2}{2d(p,\gamma_0(t))\,d(p,\gamma_1(t))}.
\]
The Alexandrov angle for a triple of points $(p,q,r)$ is defined in terms of the geodesics connecting $p$ to $q$ 
and $r$.  The Alexandrov angle $(p,q,r)\rightarrow\angle (p,q,r)$ is upper semi continuous.  Geodesics at zero angle are said to define the {\em same direction}.  At zero angle provides an equivalence relation on the geodesics beginning at a point $p$ with the Alexandrov angle providing a metric on the space of directions.  The Alexandrov tangent cone $AC_p$ is the set of constant speed geodesics beginning at $p$ modulo the equivalence relation same speed and at zero angle.  

We will compare Fenchel-Nielsen angles for terminating geodesics in terms of the geometry of level sets of the length number variation $d\mathcal L$ and the length number $\mathcal L$.  A level set of $d\mathcal L$ over a level set of $\mathcal L$ (with each $\ell_{\alpha}$ positive) has the formal structure of a $|\sigma|$-fold product of bundles $\mathbb R\rightarrow\mathbb R$ with $\mathbb R$ acting independently by translation on the fiber and the base.  The structure is explicit in terms of the local coordinates $(\vartheta_{\alpha},\ell_{\alpha}^{1/2},\ell_{\beta}^{1/2})_{\alpha\in\sigma,\,\beta\in\tau}$ for a neighborhood of the augmentation point $p$ with $(\vartheta_{\alpha})_{\alpha\in\sigma}$ a partial collection of Fenchel-Nielsen angles, $\sigma$ the short geodesics and $\tau$ a relative length basis.  In particular for a tangent vector $V$ at a point $q$ the level set of the pair $d\mathcal L(V)$ over $\mathcal L(q)$ is mapped to the pair $(d\vartheta_{\alpha})(V)\in\mathbb R^{|\sigma|}$ over $(\vartheta_{\alpha})(q)\in\mathbb R^{|\sigma|}$ with the Lie group $Flow(\sigma)$ acting by translation on the level set of $\mathcal L$ and the Lie algebra of $Flow(\sigma)$ acting by $(d\vartheta_{\alpha})(V+T_{\alpha'})=(d\vartheta_{\alpha})(V)+(\delta_{\alpha\alpha'})$.   Points $q$ and $r$ in a common level set of $\mathcal L$ are related by the element of $Flow(\sigma)$ given as $(\vartheta_{\alpha}(q)-\vartheta_{\alpha}(r))_{\alpha\in\sigma}$.  Curves $\gamma_*(t)$ and $\gamma_1(t)$ with common length numbers $\mathcal L(\gamma_*(t))=\mathcal L(\gamma_1(t))$ are related by a family of elements of $Flow(\sigma)$ given as $(\vartheta_{\alpha}(\gamma_*(t))-\vartheta_{\alpha}(\gamma_1(t)))_{\alpha\in\sigma}$ with $t$-derivative 
$(d\vartheta_{\alpha})(\gamma_*'(t))-(d\vartheta_{\alpha})(\gamma_1'(t))$.   We are ready to apply the above considerations to the Fenchel-Nielsen angle difference.

The $2$-parameter family $\kappa(s,t),\,0\le s\le 1,\,0\le t\le t'$ interpolates between $\gamma_0(t)=\kappa(0,t)$ and the terminating differentiable curve $\gamma_*(t)=\kappa(1,t)$ satisfying $\mathcal L(\gamma_*(t))=\mathcal L(\gamma_1(t))$.  We continue with the generic condition that each $\ell_{\alpha},\alpha\in\sigma,$ is non trivial on $\gamma_0(t)$ and on $\gamma_1(t)$.  

\begin{lemma}
\label{nospiral}
The geodesic $\gamma_1$ and terminating curve $\gamma_*$ have comparison angle $\angle(p,\gamma_1(t),\gamma_*(t))$ bounded as $O(t)$.  The angle is bounded as $O(t^2\log 1/t)$ for $\gamma_0$ and $\gamma_1$ with coinciding initial derivatives for $\mathcal L$.
\end{lemma}
\prf   We first use the family $\kappa(s,t)$ to bound the Fenchel-Nielsen gauge on $\gamma_*(t)$.  To this purpose consider Stoke's Theorem for the differential $1$-form $\varrho_{\alpha}$ and the rectangle: $\kappa(s,t),\,0\le s\le1,\,t''\le t\le t'$.  We bound the integral of $d\varrho_{\alpha}$ over the rectangle and the integral of $\varrho_{\alpha}$ on the three boundaries: $\kappa(0,t)=\gamma_0(t),\,t''\le t\le t'$; $\kappa(s,t'')$ and $\kappa(s,t'),\,0\le s\le1$ to obtain a bound for $\varrho_{\alpha}$ on $\gamma_*(t),\,t''\le t\le t'$.  
The gauge $\varrho_{\alpha}$ on $\gamma_0(t)$ is bounded as $O(t)$ by Corollary \ref{termbhv} and by orthogonality of $\mathcal Q$ to $\{J\lambda_{\alpha}\}$ the gauge $\varrho_{\alpha}$  vanishes on $\kappa(s,t'')$ and $\kappa(s,t')$, $0\le s\le1$.  The sum of the three boundary integrals is bounded as $O(t'(t'-t''))$.  
The $2$-dimensional integrand is $d\varrho_{\alpha}(\frac{\partial}{\partial s}\kappa,\frac{\partial}{\partial t}\kappa)$.  The derivative $\frac{\partial}{\partial s}\kappa$ is the pseudo geodesic variation field $\mathcal Q$ which by Lemma \ref{Qbound} is bounded in norm as $O(t)$ and as $O(t^2)$ for $\gamma_0$ and $\gamma_1$ with coinciding initial derivatives for $\mathcal L$.  The derivative $\frac{\partial}{\partial t}\kappa$ is the $\mathcal Q$-flow of the 
derivative $\gamma_0'(t)$ and from Lemma \ref{Qbound} the norm of $\|\frac{\partial}{\partial t}\kappa\|$ is bounded. 
Now from Lemma \ref{drho} the evaluation $d\varrho_{\alpha}(\frac{\partial}{\partial s}\kappa,\frac{\partial}{\partial t}\kappa)$ is bounded as $O(t\ell_{\alpha}^{-3/2})$ and as $O(t^2\ell_{\alpha}^{-3/2})$ for coinciding initial derivatives of $\mathcal L$. Upon integrating in $s$ it follows that the integral of $d\varrho_{\alpha}$ over the rectangle is bounded by the integral of $t\ell_{\alpha}^{-3/2}$ on $t''\le t\le t'$ and by the integral of $t^2\ell_{\alpha}^{-3/2}$ on $t''\le t\le t'$ for coinciding initial derivatives of $\mathcal L$.  The generic condition for $\gamma_0$ and $\gamma_1$, Corollary \ref{termbhv} and (\ref{kappa}) provide that $\ell_{\alpha}^{1/2}$ is bounded below by a positive multiple of $t$.  In summary it follows that the integral of the gauge on $\gamma_*(t)$ for $t''\le t\le t'$ is bounded by the integral of $t^{-2}$ on $t''\le t\le t'$ and by the integral of $t^{-1}$ on $t''\le t\le t'$ for coinciding initial derivatives of $\mathcal L$.  By the mean value theorem $\varrho_{\alpha}$ is bounded as 
$O(t^{-2})$ and as $O(t^{-1})$ for coinciding initial derivatives of $\mathcal L$.

We next introduce a second family to bound the angle displacement from $\gamma_1$ to $\gamma_*$.  As described the family in the Lie group  $Flow(\sigma)$ is given as $(\vartheta_{\alpha}(\gamma_*(t))-\vartheta_{\alpha}(\gamma_1(t)))_{\alpha\in\sigma}$  in the parameter $t$.  The differential $1$-form $d\vartheta_{\alpha}$ is closed and consequently the integrals of $d\vartheta_{\alpha}$ on $\gamma_*(t)$ and $\gamma_1(t)$ for $t''\le t\le t'$ bound the difference of the $t''$ and $t'$ displacement $(\vartheta_{\alpha}(\gamma_*(t''))-\vartheta_{\alpha}(\gamma_1(t''))) \,-\,(\vartheta_{\alpha}(\gamma_*(t'))-\vartheta_{\alpha}(\gamma_1(t')))$.   In fact the displacement is contained in $\operatorname{Ker}\,d\mathcal L$ in which case from Lemma \ref{drho} we have the comparability $\sum|\varrho_{\alpha}|\asymp(\sum\varrho_{\alpha}^2)^{1/2}\asymp(\sum d\vartheta_{\alpha}^2)^{1/2}\asymp\sum|d\vartheta_{\alpha}|$.  We now combine the Corollary \ref{termbhv} bound for the gauge on $\gamma_1(t)$ with the above bound for the gauge on $\gamma_*(t)$.  In total the difference of the $t''$ and $t'$ displacement in $Flow(\sigma)$ from $\gamma_1$ to $\gamma_*$ is bounded as $O(t''^{-1})$ in general and as $O(\log 1/t'')$ for coinciding initial derivatives of $\mathcal L$.  

We are ready to bound the WP length of the displacement from $\gamma_1(t'')$ to $\gamma_*(t'')$.  Again the sums $\sum|d\vartheta_{\alpha}|$ and $\sum |\varrho_{\alpha}|$ are comparable; the WP metric on a level set of $\mathcal L$  is comparable to $\sum\ell_{\alpha}^{3/2}|\varrho_{\alpha}|$ from Definition \ref{gauge} with $\ell_{\alpha}^{1/2}$ bounded as $O(t'')$.  Combining estimates the WP length of the displacement is bounded as $O(t''^2)$ and $O(t''^3\log 1/t'')$ for coinciding initial derivatives of $\mathcal L$.  The final item is to bound $d(p,\gamma_*(t))$ from below by a positive multiple of $t$.  A bound follows from the Corollary \ref{termbhv} expansion for the distance $d_{\mathcal T(\sigma)}$ and the lower bound for $\ell_{\alpha}^{1/2}$ in terms of $t$.  The bounds for the comparison angle follow.   The proof is complete.  

We are now ready to describe for the augmentation point $p$ an isometry between the Alexandrov tangent cone $AC_p$ and the product $\mathbb R_{\ge 0}^{|\sigma|}\times T_p\mathcal T(\sigma)$ with the first factor the Euclidean orthant and the second factor the stratum tangent space with WP metric.   The mapping for a geodesic $\gamma(t)$ terminating at $p$ is given by associating for the length numbers $\mathcal L(\gamma(t))=(\ell_{\alpha}^{1/2},\ell_{\beta}^{1/2})_{\alpha\in\sigma,\,\beta\in\tau}(\gamma(t))$  the initial one-sided derivative  
\[
\Lambda:\gamma\rightarrow(2\pi)^{1/2}\frac{d\mathcal L(\gamma)}{dt}(0).
\]
By the selection of the relative length basis the tuple $(\ell_{\beta}^{1/2})_{\beta\in\tau}$ provides local coordinates at $p$ for the stratum $\mathcal T(\sigma)$ and thus $(2\pi)^{1/2}\bigl(\frac{d\ell_{\beta}^{1/2}(\gamma)}{dt}(0)\bigr)_{\beta\in\tau}$ defines a vector in the tangent space $T_p\mathcal T(\sigma)$ with WP inner product.  The positive orthant $\mathbb R_{\ge 0}^{|\sigma|}\subset\mathbb R^{|\sigma|}$ is considered with the Euclidean inner product.  The Alexandrov tangent cone is given the structure of a cone in an inner product space through the formal relation $\langle\gamma_0,\gamma_1\rangle =\|\gamma_0'\|\|\gamma_1'\|\cos\angle (\gamma_0,\gamma_1)$.  We  present the main result.
\begin{theorem}
\label{Alextgt}
The mapping $\Lambda$ from the WP Alexandrov tangent cone $AC_p$ to $\mathbb R_{\ge 0}^{|\sigma|}\times T_p\mathcal T(\sigma)$ is an isometry of cones with restrictions of inner products.   A WP terminating geodesic $\gamma$ with a root geodesic-length function initial derivative $\frac{d\ell_{\alpha}^{1/2}(\gamma)}{dt}(0)$ vanishing is contained in the stratum $\{\ell_{\alpha}=0\},\ \mathcal T(\sigma)\subset\{\ell_{\alpha}=0\}$.  Geodesics $\gamma_0$ and $\gamma_1$ at zero angle have comparison angles $\angle(p,\gamma_0(t),\gamma_1(t))$ bounded as $O(t)$.
\end{theorem}
\prf  The primary matter is the expansions for the arc length of segments $\gamma_0(t),\gamma_*(t)$ for $0\le t\le t'$ and $\kappa(s,t')$ for $0\le s\le1$.  The goal is to express quantities in an appropriate form $ct\,+\,O(t^2)$.  By Corollary \ref{termbhv} the projection of the tangent field $\gamma_0'(t)$ onto the span of $\{J\lambda_{\alpha}\}_{\alpha\in\sigma}$ is bounded as $O(t^4)$.  From the proof of Lemma \ref{nospiral} the evaluation $\varrho_{\alpha}(\gamma_*'(t))$ is bounded as $O(t^{-2})$ and from Definition \ref{gauge} the projection of the tangent field $\gamma_*'(t)$ onto $\{J\lambda_{\alpha}\}_{\alpha\in\sigma}$ is bounded as $O(t)$ with resulting contribution $O(t'^2)$ to the length of $\gamma_*$.  
The bounds for $\langle\gamma_0',J\lambda_j\rangle $  and $\langle\gamma_1',J\lambda_j\rangle $ combined with Theorem \ref{D} also provide that the $t$-derivative of $\langle\gamma_0',\lambda_j\rangle $  and $\langle\gamma_*',\lambda_j\rangle =\langle\gamma_1',\lambda_j\rangle $ are bounded. 
The tangent fields of the curves $\gamma_0$ and $\gamma_*$ satisfy the desired form $\langle\gamma'(t),\lambda_j\rangle  =\langle\gamma'(0),\lambda_j\rangle \,+\,O(t)$.   The  curves $\gamma_0,\gamma_*$ and $\kappa$ are contained in a radius $t$ neighborhood of the point $p$.  From Corollary \ref{WPcont} the matrix of pairings $P=(\langle\lambda_j,\lambda_k\rangle )$ has an expansion of the form $ct\,+\,O(t^2)$.  Similarly the length numbers $\mathcal L(\gamma_0(t)),\,\mathcal L(\gamma_1(t))$ and $\mathcal L(\gamma_*(t))$, as well as the pseudo geodesic variation field $\mathcal Q=\sum a_k\lambda_k$ for $(a_k)=P^{-1}(\mathcal L(\gamma_1(t))-\mathcal L(\gamma_0(t)))$ all have expansions of 
the form $ct\,+\,O(t^2)$.  It now follows that the arc lengths of $\gamma_0,\gamma_*$ and $\kappa$ are given modulo terms of order $O(t^2)$ by evaluating initial tangents in the WP pairing matrix at $p$.  It further follows that the mapping $\Lambda$ is length preserving and that the coefficient for the principal term for the arc length of $\kappa$ is given by the length of $\Lambda(\gamma_0)-\Lambda(\gamma_*)$ in the tangent space of the cone. It further follows that the principal terms for the distance $d(\gamma_0(t'),\gamma_*(t'))$ and the arc length of $\kappa(s,t')$ coincide. In conclusion from Lemma \ref{nospiral} the mapping $\Lambda$ is an isometry of length and angle on the set of geodesics with all $\{\ell_{\alpha}\}_{\alpha\in\sigma}$ non trivial.  

The considerations provide that $\Lambda$ is length preserving for all elements of the Alexandrov tangent cone.  It now follows from the triangle inequality for Alexandrov angle that $\Lambda$ is an isometry on its domain, \cite[Chap. I.1, Prop. 1.14]{BH}.  To show that $\Lambda$ is a surjection consider the Taylor expansion for the length number.  It follows from Theorem \ref{D} and Corollary \ref{termbhv} that on a geodesic $\gamma$ then $\mathcal L(t)=\mathcal L(0)\,+\,\mathcal L'(0)t\,+\,O(t^2\|\gamma\|)$ with a uniform constant for the remainder.  We consider $(\mathcal L(t)-\mathcal L(0))t^{-1}= \mathcal L'(0) \,+\,O(t\|\gamma\|)$ and observe that since $\mathcal L$ is an open map the left hand side limits to a surjective map to $\mathbb R_{\ge 0}^{|\sigma|}\times T_p\mathcal T(\sigma)$. And finally Corollary \ref{termbhv} provides for the strata and initial geodesic-length derivative behavior for terminating geodesics.   The proof is complete.    

The structure of $AC_p$ provides simple invariants for the point $p$.   The dimension $\dim AC_p$ of the cone and the dimension $\dim \mathcal T(\sigma)$  of the maximal vector subspace are intrinsic.   Additionally the representation of an orthant as a product of half lines is unique modulo permutation of factors.   We now observe that a pair of $\Tbar$ closed strata are either disjoint or one contained in the other or intersect orthogonally. In particular for a component of a proper intersection of closed strata $\mathcal S_1$ and $\mathcal S_2$ there are simple free homotopy classes $\alpha_1$ and $\alpha_2$ with $\alpha_1$ a $\sigma$-null for $\mathcal S_1$ but not for $\mathcal S_2$ and $\alpha_2$ a $\sigma$-null for $\mathcal S_2$ but not for $\mathcal S_1$.  At a point $p$ of the intersection the half lines for $\alpha_1$ and $\alpha_2$ are tangents to $\mathcal S_1$ and $\mathcal S_2$ in $AC_p$ and are orthogonal.

An important property for non positively curved Riemannian manifolds is that the exponential map is distance non decreasing.  An inverse exponential map $exp^{-1}_p:\Tbar\rightarrow AC_p$ is defined by associating to $q\in\Tbar$ the unique geodesic connecting $p$ to $q$ with speed $d(p,q)$.  The map is not an injection since geodesics at zero angle with common speed are mapped to a common element of $AC_p$.  The map is distance non increasing as follows.  From the $CAT(0)$ inequality and definition of the Alexandrov angle points on geodesics beginning at $p$ have distance satisfying
\begin{multline*}
d(\gamma_0(t),\gamma_1(t))^2 \ge d(p,\gamma_0(t))^2\,+d(p,\gamma_1(t))^2 \\ -2d(p,\gamma_0(t))\,d(p,\gamma_1(t))\cos\angle (\gamma_0,\gamma_1).
\end{multline*}
For equality for a single value of $t$ the Flat Triangle Lemma \cite[Chap. II.2 Prop. 2.9]{BH} provides that the geodesics are contained in a flat subspace of $\Tbar$.   The flat subspaces were characterized in \cite[Prop. 16]{Wlcomp}.  

A further feature of the Alexandrov angle is a first variation formula for the distance from a point to a geodesic.  For a unit-speed geodesic $\gamma(t)$ beginning at $p$ the distance $d(\gamma(t),q)$ to a point not on the geodesic is convex with initial one-sided derivative satisfying
\[
\frac{d}{dt}d(\gamma(t),q)(0^+)=-\cos\angle (\gamma,\gamma_{pq})
\]
for $\gamma_{pq}$ the geodesic connecting $p$ to $q$, \cite[Chap. II.3, Coro. 3.6]{BH}.  We consider simple applications of the variational formula.

The first is Yamada's observed {\em non refraction} of WP geodesics on $\Tbar$: a length minimizing path at most changes strata at its endpoints, \cite{DW2, Wlcomp}.  Specifically consider a pair of unit-speed geodesics $\gamma_0(t),\, \gamma_1(t)$ with initial point $p$ such that the reverse path along $\gamma_0$ followed by $\gamma_1$ is length minimizing.  The Alexandrov angle between the tangents at $p$ is $\pi$.  Specifically the comparison distance $d(\gamma_0(t),\gamma_1(t))$ is at least that of the path from $\gamma_0(t)$ to $p$ to $\gamma_1(t)$ and thus $\lim_{t\rightarrow 0} d(\gamma_0(t),\gamma_1(t))/2t=1$ and the angle is $\pi$.  Elements of $AC_p$ at angle $\pi$ necessarily lie in the subspace $T_p\mathcal T(\sigma)$ and from the Theorem are tangent to paths contained in $\mathcal T(\sigma)$. It further follows that the geodesics are segments on a single geodesic.  The second application concerns projecting geodesics; see Definition \ref{termproj}.  A projecting geodesic $\varsigma(t)$ with initial point $p$ has tangent in $AC_p$ orthogonal to the subspace 
$AC_p(\overline{\mathcal T(\sigma)})$ as follows.  An element of the Alexandrov tangent cone of $\overline{\mathcal T(\sigma)}$ is represented by geodesics $\gamma(t')$ beginning at $p$ and contained in the stratum.  Minimizing the distance to the stratum provides that $d(\varsigma(t),p)\le d(\varsigma(t),\gamma(t'))$ for all $t,\,t'$.  From the variational formula it follows that $\angle(\varsigma,\gamma)\ge\pi/2$.  Orthogonality now follows from two observations.  A pair of vectors in a positive Euclidean orthant form an angle of at most $\pi/2$. The orthogonal complement in $AC_p$ of the Euclidean orthant is the linear subspace of geodesics containing $p$ as an interior point.

The third application concerns limits of geodesics, \cite{Brkwpvs}.    
\begin{example}  
Length-minimizing concatenations of geodesics.
\end{example}

A characterization of limits in terms of concatenations of geodesic segments is provided in \cite[Sec. 7]{Wlcomp}.  The resulting concatenations were found to be the unique length-minimizing paths connecting an initial and terminal point and intersecting a prescribed sequence of closures of strata.   We apply the above variational formula to show that such concatenation paths satisfy a strong form of the classical {\em angle of incidence and reflection equality.}  To this purpose consider a pair of geodesics $\gamma_0$ and $\gamma_1$ each with initial point $p$ on $\overline{\mathcal T(\sigma)}$, $\gamma_0$ with endpoint $q$ and $\gamma_1$ with endpoint $r$.  Consider that the concatenation $\gamma_0+\gamma_1$ is a length-minimizing path connecting $q$ and $r$ to a point of $\overline{\mathcal T(\sigma)}$.  A geodesic $\kappa$ beginning at $p$ contained in $\overline{\mathcal T(\sigma)}$ provides a variation of the configuration.   The initial derivative of distance $d(q,p)+d(r,p)$ along $\kappa$ is $-\cos\angle(\gamma_0,\kappa)-\cos\angle(\gamma_1,\kappa)$.  The geodesics beginning at $p$ contained in $\overline{\mathcal T(\sigma)}$ fill out the Alexandrov tangent cone $AC_p(\overline{\mathcal T(\sigma)})$.  It follows that the sum in $AC_p$ of the initial tangents of $\gamma_0$ and $\gamma_1$ has vanishing projection onto the subcone $AC_p(\overline{\mathcal T(\sigma)})$, the first conclusion. 

Equality at $p$ of the initial derivatives 
$\frac{d\ell_{\alpha}^{1/2}(\gamma_0)}{dt}$ and $\frac{d\ell_{\alpha}^{1/2}(\gamma_1)}{dt}$ is the further property of concatenations resulting as limits of geodesics of $\mathcal T$.   The basic observation is that along a geodesic $\gamma(t)$ the projection of the tangent onto the span of $\{\lambda_{\alpha},J\lambda_{\alpha}\}$ has square norm 
$f(t)=\langle\lambda_{\alpha},\frac{d}{dt}\rangle^2+\langle J\lambda_{\alpha},\frac{d}{dt}\rangle^2$ a Lipschitz function with small constant.  In particular from Theorem \ref{D} the $t$-derivative of $f(t)$ has vanishing principal term and is consequently bounded as $O(\ell_{\alpha}^{3/2})$.  On a segment of $\gamma(t)$ with $\ell_{\alpha}^{1/2}\le c_0$ the function $f(t)$ is Lipschitz with constant $O(c_0^3)$. We consider a sequence of segments of geodesics converging to the concatenation, \cite[Sec. 7]{Wlcomp}.  The segments converge to $\gamma_0+\gamma_1$ in the compact-open topology of $\mathcal T$ with 
$\langle\lambda_{\alpha},\frac{d}{dt}\rangle$ and $\langle J\lambda_{\alpha},\frac{d}{dt}\rangle$ converging correspondingly.  
We combine the expansions $\langle\lambda_{\alpha},\frac{d}{dt}\rangle=a_{\alpha}(\gamma_j)+O(t\ell_{\alpha}^3/2)$, 
$j=1,2$, $\langle J\lambda_{\alpha},\frac{d}{dt}\rangle=O(t\ell_{\alpha}^3/2)$ of Corollary \ref{termbhv} for 
$\gamma_0$ and $\gamma_1$ with the Lipschitz bound for $f(t)$ for $\ell_{\alpha}^{1/2}\le c_0$ to find that $a_{\alpha}(\gamma_0)^2=a_{\alpha}(\gamma_1)^2+O(c_0^3)$.   The conclusion follows. 

We describe an application for homomorphic group actions by isometries.  Farb and Masur \cite{FaMa}, Bestvina and Fujiwara \cite{BeFu}, as well as Hamenstadt \cite{Ham} have important general results on the finiteness of group homomorph images into the mapping class group.  
\begin{example}
Combinatorial harmonic maps.
\end{example}

Wang \cite{WgMT1,WgMT2}, followed by Izeki and Nayatani \cite{IN} have studied equivariant energy minimizing maps from simplicial complexes to complete $CAT(0)$ spaces.  This is the setting for Margulis superrigidity for lattices in semisimple algebraic groups over $p$-adic fields.  Wang formulated a notion of energy for an equivariant map, provided the target is locally compact and the group action is reductive.  A criterion for existence of a energy minimizing map in a homotopy class is presented in terms of a Poincar\'{e}-type inequality for maps of the links of vertices to the target tangent cones.  Izeki and Nayatani extended the approach of Wang to include non locally compat $CAT(0)$ spaces and to remove the hypothesis for a reductive group action.  We now present the Izeki and Nayatani results \cite[Thrms. 1.1 and 1.2]{IN} for the setting of the augmented Teichm\"{u}ller space and describe examples of simplicial complexes with suitable group actions.
\begin{theorem}
Let $X$ be a connected simplicial complex with connected links of vertices.  Let $\Gamma$ be a finitely generated group acting properly discontinuously by automorphisms of $X$ with compact quotient.  Assume for the link of each vertex of $X$ that the first non zero eigenvalue of the Laplacian is greater than $1/2$.  Further let $Y$ be a complete $CAT(0)$ space with each Alexandrov tangent cone isometric to a closed convex cone in a Hilbert space.  Then an isometric action of $\Gamma$ on $Y$ has a  fixed point.
\end{theorem}
The augmented Teichm\"{u}ller space $\Tbar$ provides an example of a suitable target space with the WP isometry group the extended mapping class group \cite{BrMr, MW,Wlcomp}.  The fixed points of WP isometries are 
understood \cite{DW2, Wlcomp}.  A WP isometry fixing a point of $\mathcal T$ is realized as a conformal automorphism of the corresponding Riemann surface.  A WP isometry fixing a point of $\Tbar-\mathcal T$ is reducible pseudoperiodic (fixes a simplex of curves with a power having periodic component mappings) in the sense of Thurston, \cite{Thsurf}.  A group fixing a point of $\Tbar-\mathcal T$ is realized as an extension of the conformal automorphisms of the corresponding Riemann surface with nodes.  The extension is by a group of Dehn twists about the curves in the simplex for the null stratum of the fixed point.  

Examples of source spaces with group actions are provided by suitable Euclidean buildings and 
Ballman-\'{S}wi\c{a}tkowski complexes, \cite[Examples 1 and 2]{IN}.  A building is a simplicial complex that can be expressed as a union of subcomplexes (called apartments) satisfying a certain set of axioms.  A building is Euclidean if its apartments are isomorphic to a Euclidean Coxeter complex.  The quotient space $X=PGL(n;\mathbb Q_p)\slash PGL(n;\mathbb Z_p)$ is the vertex set of a Euclidean building for $p$ a prime, $\mathbb Q_p$ the $p$-adic number field and $\mathbb Z_p$-the $p$-adic integers.  For $n=3$, $X$ is a two-dimensional simplicial complex with links of vertices all isomorphic to a common regular bipartite graph with first non zero eigenvalue of the Laplacian greater than $1/2$.  Examples of suitable groups $\Gamma$ are provided by cocompact lattices in $PGL(3;\mathbb Q_p)$.

Ballman and \'{S}wi\c{a}tkowski described contractible two-dimensional simplicial complexes with links of vertices all isomorphic.  They considered a finite group $H$ with a set $S$ of non trivial generators and the Cayley graph $\mathcal C$ of $H$ with respect to $S$.  Provided the minimal number of edges of a closed circuit in $\mathcal C$ is at least $6$ there exists a two-dimensional simplicial complex $X$ with all links of vertices isomorphic to $\mathcal C$.  Provided $H$ has a presentation $\langle S\mid R\rangle$ then the group $\Gamma$ with presentation $\langle S\cup \{\tau\}\mid R\cup\{\tau^2\}\cup\{(s\tau)^3\mid s\in S\}\rangle$ acts properly discontinuously by automorphisms on $X$ with compact quotient.  Sarnak has described Cayley graphs of finite groups with first non zero eigenvalue of the Laplacian greater than $1/2$, \cite[Chap. 3]{Srbk}.


\end{document}